\documentclass[3p]{elsarticle}
\usepackage{amsmath}
\usepackage{amsfonts}
\usepackage{bm}
\usepackage{mathtools}
\usepackage{subfiles}
\usepackage{algorithm,algorithmic}
\usepackage{graphbox}                   % To vertically align images
\usepackage{amsthm}
\usepackage[hypertexnames=false]{hyperref}
\usepackage{cleveref}                   % Best to load last

\newtheorem{remark}{Remark}
\newtheorem{claim}{Claim}
\newtheorem{lemma}{Lemma}
\newtheorem{corollary}{Corollary}
\newtheorem{theorem}{Theorem}   
\crefname{claim}{claim}{claims}

\Crefname{ALC@unique}{Line}{Lines} % <- Preamble

\DeclareMathOperator*{\minimize}{minimize}
\DeclareMathOperator{\Tr}{Tr}
\DeclareMathOperator{\Rank}{Rank}

\begin{document}

\begin{frontmatter}

    \title{Derivative-based SINDy (DSINDy): Addressing the challenge of discovering governing equations from noisy data}

    % Authors: full names plus addresses.
    \author[1]{Jacqueline Wentz}
    \ead{Jacqueline.Wentz@colorado.edu}

    \author[1]{Alireza Doostan}
    \ead{Alireza.Doostan@colorado.edu}

    \affiliation[1]{
        organization={Smead Aerospace Engineering Sciences, University of Colorado},
        addressline={3775 Discovery Dr, Boulder, CO, 80303},
        country={United States}}

    \begin{abstract}
        Recent advances in the field of data-driven dynamics allow for the discovery of
        ODE systems using state measurements.
        One approach, known as Sparse Identification of Nonlinear Dynamics (SINDy), assumes the dynamics are sparse within a predetermined basis in the states and finds the expansion coefficients through linear regression with sparsity constraints. This approach requires an accurate estimation of the state time derivatives, which is not necessarily possible in the high-noise regime without additional constraints. We present an approach called Derivative-based SINDy (DSINDy) that combines two novel methods to improve ODE recovery at high noise levels. First, we denoise the state variables by applying a projection operator that leverages the assumed basis for the system dynamics. Second, we use a second order cone program (SOCP) to find the derivative and governing equations simultaneously.
        We derive theoretical results for the projection-based denoising step, which allow us to estimate the values of hyperparameters used in the SOCP formulation.
        This underlying theory helps limit the number of required user-specified parameters. We present results demonstrating that our approach leads to improved system recovery for the Van der Pol oscillator, the Duffing oscillator, the R\"{o}ssler attractor, and the Lorenz 96 model.

    \end{abstract}


\end{frontmatter}

\section{Introduction}\label{sec:introduction}

Physical systems are often represented using differential equations to describe how the systems evolve over time. From simple systems, such as the pendulum, to complex systems, such as the interaction of hundreds of proteins within a cell, differential equations provide a way to predict future states and study stability properties within a system. However, often the governing equations of systems are either partially or fully unknown. In these cases, methods for data-driven equation discovery can be applied to learn the governing equations directly from state measurements \cite{Ghadami2022}. These data-driven methods have shown success, for example, in learning chemical laws \cite{Neumann2020}, biological networks \cite{Mangan2016}, ecological systems \cite{Dam2017}, and fluid dynamics \cite{Raissi2020}. However, some fundamental challenges in equation discovery remain, including accurate state time derivative estimation in the presence of noise, robustness with respect to noise, and the need for user-specified hyperparameters.

We focus specifically on the data-driven discovery of ordinary differential equation (ODE) systems. The best method for discovering a given ODE system depends on characteristics of the problem, such as the number of states, prior information, and the goal of the discovery, e.g., interpretability, forecasting, or stability analysis. For instance, although symbolic regression is highly expressive and leads to interpretable governing equations \cite{Bongard2007,Quade2016,Neumann2020}, it is a computationally infeasible approach for high dimensional systems. Similarly, neural-networks allow for high expressivity and require limited information on the form of the nonlinearities \cite{Raissi2018}, but the resulting solution is often uninterpretable and requires many user-specific hyperparameters, e.g., in defining the neural network architecture. Recently developed sparsity-promoting methods, e.g., Sparse Identification of Nonlinear Dynamics (SINDy), assume the governing equations are sparse in some known basis of the state variables \cite{Schmidt2009, Wang2011, Brunton2016}.
This leads to parsimonious and interpretable solutions but requires potentially unavailable knowledge on the form of the basis. In contrast to neural networks and symbolic regression, sparsity-promoting methods allow for the rapid discovery of an ODE system by solving a sparsity regularized least-squares optimization problem.

Here, we will focus on sparsity-promoting techniques as we are interested in interpretable equations. One of the main challenges with these methods is the accurate approximation of the state time derivatives in the presence of measurement noise. A variety of methods have been proposed to circumvent this issue such as \textit{a priori} smoothing \cite{Cortiella2022,Delahunt2022}, using an integral formulation \cite{Schaeffer2017}, and modeling the state variables with a neural network to find derivatives through automatic differentiation \cite{Chen2021}. Several groups have explored simultaneous denoising and recovery of the system dynamics \cite{Tran2017,Sun2021,Hokanson2022,Kaheman2022}. Although promising, the resulting optimization problem is non-convex and often involves the tuning of hyperparameters. One of the most successful approaches uses the weak form of the dynamics \cite{Wang2019, Gurevich2019, Reinbold2020, Messenger2021, Messenger2021a}. This approach, known as Weak SINDy or WSINDy, applies integration by parts and places the derivative on test functions to avoid the derivative calculation from noisy data. One potential downside to WSINDy is there are several hyperparameters that may need tuning in order to obtain reasonable results \cite{Messenger2021}.

Our approach for discovery governing equations, which we call Derivative-based SINDy (DSINDy), seeks to improve upon the aforementioned methods in several ways. Instead of \textit{a priori} smoothing without any physical information \cite{Cortiella2022}, we apply a denoising step that leverages the assumed basis for the system dynamics. Additionally, in contrast to simultaneous denoising and recovery methods \cite{Tran2017,Sun2021,Hokanson2022,Kaheman2022}, our approach for system recovery only involves solving a convex program. Finally, by using theoretical results from the denoising step, we recover the system dynamics without requiring the tuning of hyperparameters. Although a comprehensive comparison of DSINDy with other equation discovery methods in the high noise regime is beyond the scope of this work, we do provide a comparison of DSINDy with both WSINDy and an $\ell_1$-minimization version of SINDy \cite{Cortiella2021}, which we refer to as $\ell_1$-SINDy. Note that in $\ell_1$-SINDy, Tikhonov regularization is used to find smooth time derivatives from the noisy measurements.

The DSINDy algorithm can be broken down into two steps. First we use a novel algorithm, called Projection-based State Denoising (PSDN), which leverages, through a projection operation, the assumed basis for the system dynamics. We next formulate an iteratively reweighted, second order cone program (IRW-SOCP) that allows us to directly find the state time derivatives, as opposed to the states themselves, while enforcing sparsity of the coefficients, hence the name ``Derivative"-based SINDy. A theoretical analysis of PSDN allows us to estimate the values of hyperparameters needed for IRW-SOCP. When used in concert, PSDN and IRW-SOCP lead to the best performance in terms of coefficient estimation and system recovery. Thus, although these two methods can be applied independently, we present them together as a comprehensive approach for learning the governing equations of an ODE system.

Our paper is outlined as follows. In \Cref{sec:statement}, we state the equation discovery problem and give relevant notation. In \Cref{sec:methods}, we describe DSINDy in detail, and, in \Cref{sec:theory}, we present theoretical convergence results. Finally, in \Cref{sec:ode-systems,sec:results}, we introduce four example ODE systems and compare the recovery performance of DSINDy with WSINDy and $\ell_1$-SINDy.

\section{Problem statement and notation}\label{sec:statement}

We consider a system with $m$ state variables that change over time according to
\begin{align}
    \dot{u}_k^*(t) & = F_k(u_1^*(t),u_2^*(t),...,u_m^*(t)) \quad\text{for}\quad k=1,2,...,m,
\end{align}
where we use dot notation to refer to the time derivative.
Here, $F_k$ for $k=1,2,...,m$ are unknown functions of the state variables.
Note we use an asterisk $(*)$ to refer to true values and derivatives of the state variables, in contrast to noisy estimates.

One approach for discovering the ODE system is to assume that $F_k$ can be written as a linear combination of known basis functions \cite{Brunton2016}.
We make this assumption here, but note that for many systems this information may not be available.
We define the known set of basis functions as $\{\theta_j\}_{j=1}^p$, where each function acts on the state space and returns a scalar, i.e., $\theta_j:\mathbb{R}^m \rightarrow \mathbb{R}$ for $j=1,2,...,p$.
By assumption, for each state variable $u_k$ there is a vector $\bm{c}_k = [c_{k,1},...,c_{k,p}]^T \in \mathbb{R}^p$, such that
\begin{equation}
    \dot{u}_k^*(t) = \sum_{j=1}^p c_{k,j}  \theta_j(u_1^*(t),u_2^*(t),...,u_m^*(t)).
\end{equation}
As in SINDy, the goal of this work is to recover $\bm{c}_k$, for $k=1,2,...,m$, using a set of noisy state measurements.

We model noise as an additive term on the true state variable trajectory. That is, for a system with $N$ measurements obtained at times $0=t_1,t_2,...,t_N=t_{end}$, the $k$-th state measurement vector $\bm{u}_k = [u_{k,1}, u_{k,2}, ..., u_{k,N}]^T$ is
\begin{equation}
    \bm{u}_k = \bm{u}^*_k + \bm{\epsilon}_k \in \mathbb{R}^N,
    \quad\text{for}\quad k=1,2,...,m,
\end{equation}
where $\bm{u}^*_k = [u^*_k(t_1),u^*_k(t_2),...,u^*_k(t_N)]^T$ and $\bm{\epsilon}_k=[\epsilon_{k,1},\epsilon_{k,2},...,\epsilon_{k,N}]^T$. We assume the noise is an i.i.d. Gaussian random variable with zero mean and variance $\sigma^2$, i.e., $\epsilon_{k,i} \sim \mathcal{N}(0,\sigma^2)$ for $k=1,2,...,m$ and $i=1,2,...,N$. In this work, we only consider systems with uncorrelated measurement noise.

We define the libraries $\Theta \in \mathbb{R}^{N\times p}$ and $\Theta^* \in \mathbb{R}^{N\times p}$ to contain evaluations of the basis functions at the noisy and true state variable values, respectively. Each row corresponds to a measurement time and each column corresponds to a basis function such that, for $i=1,2,...,N$ and $j=1,2,...,p$,
\begin{equation}\label{eq:Theta}
    \Theta_{ij} = \theta_j(u_{1,i},u_{2,i},...,u_{m,i}) \quad \text{ and } \quad \Theta^*_{ij} = \theta_j(u_{1,i}^*,u_{2,i}^*,...,u_{m,i}^*),
\end{equation}
where $u_{k,i}^* = u^*_k(t_i)$.
Let $\dot{\bm{u}}^*_{k} \in \mathbb{R}^N$ be the vector of true time derivatives of state variable $u_k$ at times $t_1,...,t_N$. It follows that
\begin{equation}\label{eq:Thetastarc_equal_udotstar}
    \Theta^* \bm{c}_k = \dot{\bm{u}}^*_k.
\end{equation}
In practice, the true values and time derivatives of the state variables are unknown. Therefore, we use an approximated version of \Cref{eq:Thetastarc_equal_udotstar}, i.e.,
\begin{equation}\label{eq:Thetac_equal_udot}
    \Theta \bm{c}_k \approx \dot{\bm{u}}_k,
\end{equation}
to find the coefficients. Here, $\dot{\bm{u}}_k$ is an approximation of the true time derivative vector, $\dot{\bm{u}}^*_k$. Noise in the state measurements impacts both sides of \Cref{eq:Thetac_equal_udot}, which may cause inaccurate coefficient recovery. As discussed in \Cref{sec:introduction}, many groups have started with \Cref{eq:Thetac_equal_udot} to recover the governing equations of ODE systems, e.g., \cite{Wang2011, Brunton2016, Tran2017, Schaeffer2017, Sun2021, Reinbold2020, Messenger2021,Cortiella2021}. Note that similar to these works we only consider the case where the system is overdetermined, i.e., $N \ge p$.

In \Cref{sec:methods} we present our novel DSINDy approach for finding the coefficient vector $\bm{c}_k$. In \Cref{sec:results} we compare DSINDy to two other approaches, i.e., WSINDy and $\ell_1$-SINDy, which are described in detail in \Cref{sec-app:algorithms}.

\subsection{Multivariate monomial basis}

We define the basis $\{\theta_j\}_{j=1}^p$ to be the set of multivariate monomials up to total degree $d$. The $j$-th basis function is given as,
\begin{equation}
    \theta_j(u_1,u_2,...,u_m) = \prod_{k=1}^m (u_k)^{\alpha_k^{(j)}}.
\end{equation}
We use $\bm{\alpha}^{(j)} = [\alpha^{(j)}_1,\alpha^{(j)}_2...,\alpha^{(j)}_m] \in \mathbb{N}_{0}^m$ to denote the multi-index vector corresponding to the $j$-th element of the monomial basis such that $\sum_{k=1}^m \alpha_k^{(j)} \le d$. As an example, for a system where $m=2$ and $d=2$, there are $p=6$ basis terms where
\begin{equation}
    \begin{aligned}
        \bm{\alpha}^{(1)} & = [0,0], &
        \bm{\alpha}^{(2)} & = [1,0], &
        \bm{\alpha}^{(3)} & = [0,1]    \\
        \bm{\alpha}^{(4)} & = [2,0], &
        \bm{\alpha}^{(5)} & = [1,1], &
        \bm{\alpha}^{(6)} & = [0,2].
    \end{aligned}
\end{equation}
Although we use the set of multivariate monomials as the basis, the methods presented in this paper could be applied to systems with other types of basis elements.

\subsection{Additional notation}\label{sec:notation}

Throughout the paper we use an asterisk ($*$) and tilde ($\sim$) to refer to true state values and denoised measurements, respectively. When referencing the original noisy data no decorations are used. Analogous decorations are used to refer to functions or matrices evaluated with a given state version. When the same method or logic can be applied to each state in the system, we often drop the subscript $k$ and, e.g., refer to an arbitrary state vector as $\bm{u}\in\mathbb{R}^N$.

In DSINDy we use a discrete integral operator $T \in \mathbb{R}^{N\times N}$, which performs Trapezoidal quadrature, i.e.,
\begin{equation}\label{eq:T}
    T := \frac{\Delta t}{2}
    \begin{bmatrix}
        0      & 0      & 0 & ...    & 0 & 0      \\
        1      & 1      & 0 & ...    & 0 & 0      \\
        1      & 2      & 1 &        & 0 & 0      \\
        \vdots & \vdots &   & \ddots &   & \vdots \\
        1      & 2      & 2 & ...    & 1 & 0      \\
        1      & 2      & 2 & ...    & 2 & 1
    \end{bmatrix},
\end{equation}
where $\Delta t = t_{end}/(N-1)$. If $\bm{u}^*, \dot{\bm{u}}^* \in \mathbb{R}^N$ are the vectors of the true values and true time derivatives of a state variable, then
\begin{equation}\label{eq:quad-err}
    u^*_0 + T\dot{\bm{u}}^* = \bm{u}^* + \bm{e}_q,
\end{equation}
where $\bm{e}_q \in \mathbb{R}^N$ is error introduced by performing trapezoidal quadrature (note that $e_{q,1}=0$). We additionally use a finite difference matrix $D \in \mathbb{R}^{(3N - 3) \times N}$ to penalizes non-smooth solutions, i.e.,
\begin{equation}
    D := \begin{bmatrix}
        I \\ D_1 \\ D_2
    \end{bmatrix},
\end{equation}
where $D_1 \in \mathbb{R}^{(N - 1) \times N}$ and $D_2 \in \mathbb{R}^{(N - 2) \times N}$ are the first and second order finite difference operators, respectively. That is
\begin{equation}
    D_1 := \frac{1}{\Delta t} \begin{bmatrix}
        -1 & 1  &        &        \\
           & -1 & 1      &        \\
           &    & \ddots & \ddots
    \end{bmatrix}, \quad \quad
    D_2 := \frac{1}{\Delta t^2} \begin{bmatrix}
        1 & -2 & 1      &        &        \\
          & 1  & -2     & 1      &        \\
          &    & \ddots & \ddots & \ddots
    \end{bmatrix}.
\end{equation}

Our denoising strategy involves projecting a vector $\bm{x} \in \mathbb{R}^N$ onto the column space of a matrix $A \in \mathbb{R}^{N \times M}$. The result is given as $P_A \bm{x}$ where
\begin{equation}\label{eq:P_A}
    P_A := A A^{\dagger}.
\end{equation}
Here and in other sections, we use a dagger ($\dagger$) to represent the pseudoinverse, i.e., $A^\dagger:=(A^T A)^{-1} A^T$.

\begin{remark}
    To avoid numerical instabilities, we find $P_A$ using the singular value decomposition, i.e., $P_A := \hat{U} \hat{U}^T$, where $\hat{U}$ contains the left singular vectors of $A$ that correspond to nonzero singular values.
\end{remark}

For an arbitrary vector $\bm{x}$, we use $x_i$ to denote the $i$th element. The expected value and variance of a random variable $X$ are denoted as $\mathbb{E}[X]$ and $\mathbb{V}[X]$, respectively. We use $\bm{1}_N \in \mathbb{R}^N$ to represent a vector of ones and unless otherwise noted $\|\cdot\|=\|\cdot\|_2$.

\section{DSINDy Overview}\label{sec:methods}

The DSINDy algorithm involves first denoising the state measurements (see \Cref{sec:methods-PSDN}) and then using an SOCP formulation to find the coefficients and derivatives (see \Cref{sec:methods-socp}). Here, we drop the subscript $k$ and use, for example, $\bm{u}$ to refer to an arbitrary state in the system.
% DSINDy is implemented in python and all codes are available at the GitHub repository. \todo{I need to add code to GitHub}

\subsection{Projection-based state denoising}\label{sec:methods-PSDN}

To motivate the PSDN approach, we start with the equation
$\Theta^* \bm{c} = \dot{\bm{u}}^*$; see \Cref{sec:statement}. We discretely integrate both sides, using the operator $T$ given by \Cref{eq:T}, to obtain
\begin{equation}\label{eq:T_Thetastar_c}
  T \Theta^* \bm{c} = \bm{u}^* - u_0^* + \bm{\epsilon}_q,
\end{equation}
where $\bm{\epsilon}_q$ is the error introduced by quadrature. Ignoring this error, \Cref{eq:T_Thetastar_c} can be rewritten as
\begin{equation}\label{eq:Phi-b-eq-u}
  \Phi^* \bm{b} \approx \bm{u}^*,
\end{equation}
where
\begin{equation}\label{eq:b_and_Phistar}
  \bm{b} := \begin{bmatrix} u^*_0 \\ \bm{c} \end{bmatrix}
  \quad \text{and} \quad
  \Phi^* := \begin{bmatrix} \bm{1}_N & T \Theta^* \end{bmatrix}.
\end{equation}
From \Cref{eq:Phi-b-eq-u} it is clear that the true state vector $\bm{u}^*$ approximately lies in the column space of $\Phi^*$. With this in mind, our denoising strategy is to project the measurements, $\bm{u}$, onto the column space of an approximation of $\Phi^*$.

If $\Phi^*$ were known, we could directly perform this projection to find
\begin{equation}\label{eq:tildeustar}
  \tilde{\bm{u}}^* := P_{\Phi^*} \bm{u},
\end{equation}
where $P_{\Phi^*}$ is the projection operator as defined in \Cref{eq:P_A}.
As $N \rightarrow \infty$ and $t_{end}$ is held constant, we can show that $\tilde{\bm{u}}^*$ approaches $\bm{u}^*$ (see \Cref{lemma:err_bound_1}). In practice, $\Phi^*$ must be approximated from the noisy measurements as follows,
\begin{equation}\label{eq:Phi}
  \hat{\Phi} := \begin{bmatrix} \bm{1}_N & T \hat{\Theta}\end{bmatrix},
\end{equation}
where $\hat{\Theta}$ is an unbiased estimator of $\Theta^*$ such that $\mathbb{E}[\hat{\Theta}]=\Theta^*$ (see \Cref{remark:theta-hat}). For details on how to construct $\hat{\Theta}$ for the monomial basis see \Cref{sec-app:centered-library}. The PSDN approach is then summarized as a projection of the data onto the column space of $\hat{\Phi}$,
\begin{equation}\label{eq:tildeu}
  \tilde{\bm{u}} := P_{\hat{\Phi}} \bm{u} \quad \quad \text{(PSDN)}.
\end{equation}
Even though we are introducing additional error due to the noise in $\hat{\Phi}$, the result given by \Cref{eq:tildeu} has similar asymptotic properties as $\tilde{\bm{u}}^*$. Specifically, as $N\rightarrow\infty$ and $t_{end}$ is held constant, $\tilde{\bm{u}}$ approaches $\bm{u}^*$ (see \Cref{thm:main}).

\begin{remark}\label{remark:theta-hat}
  Using $\hat{\Theta}$ to find $\tilde{\bm{u}}$ helps simplify the proofs in \Cref{sec:theory}. The calculation of $\hat{\Theta}$ may not be possible for alternative basis functions, but, in practice, we did not observe a significant change in performance when the projection was performed onto the column space of an uncentered library.
\end{remark}

For some systems, PSDN led to near optimal denoising, but this was not always the case (see \Cref{sec:results-PSDN-theory}). Therefore, we also consider an iterative approach, i.e., IterPSDN, where we gradually project the data onto the column space of a sequence of $\Phi^*$ estimates (see \Cref{alg:PSDN}).
At each iteration, we first find an estimate of $\Phi^*$ based on the current state vector estimates and then perform a partial projection onto the column space of the $\Phi^*$ estimate as shown on \Cref{alg-line:proj} of \Cref{alg:PSDN}.

\begin{algorithm}
  \caption{IterPSDN$\left(\{\bm{u}_k\}_{k=1}^m; \alpha, \left\{\sigma_k\right\}_{k=1}^m,\text{CheckDiverg}\right)$}
  \label{alg:PSDN}
  \begin{algorithmic}[1]
    \STATE Set $\bm{u}_k^{(0)} = \bm{u}_k$ for $k=1,2,...,m$
    \FOR{$i=0,1,2,...$}
    \STATE Let $\Theta^{(i)}$ be monomial library evaluated at $\{\bm{u}^{(i)}_k\}_{k=1}^m$ \;
    \STATE Set $\Phi^{(i)} = \begin{bmatrix} \bm{1}_N & T \Theta^{(i)}\end{bmatrix}$ \;
    \FOR{$k=1,2,...,m$}
    \STATE Set $\bm{u}^{(i+1)}_k = \alpha P_{\Phi^{(i)}}  \bm{u}^{(i)}_k + (\alpha - 1)  \bm{u}^{(i)}_k$ \label{alg-line:proj}
    \COMMENT{Perform partial projection}
    \IF{CheckDiverg \AND $\frac{1}{\sqrt{N}}\|\bm{u}_k^{(i+1)}-\bm{u}_k^{(0)}\| > \sigma_k$}
    \STATE Set $\bm{u}^{(i+1)}_k = \bm{u}^{(i)}_k$
    \COMMENT{If diverging, revert state values}
    \ENDIF
    \ENDFOR
    \IF{$\max_k \left(\frac{\|\bm{u}_k^{(i+1)} - \bm{u}_k^{(i)}\|}{\|\bm{u}_k^{(i)}\|}\right) < 10^{-8}$}
    \STATE Break
    \ENDIF
    \ENDFOR
    \RETURN{$\{\bm{u}_k^{(i)}\}_{k=1}^m$}
  \end{algorithmic}
\end{algorithm}

\Cref{alg:PSDN} requires a few additional input parameters. First, we introduce $\alpha \in [0,1]$. If $\alpha=1$, a full projection is performed and, if $\alpha=0$, the state vector estimates are not changed. For $\alpha < 1$, a weighted average of the projection result and current state vector estimate is used to generate the new state vector estimate. In practice, we set $\alpha=0.1$ but note that for $\alpha < 0.1$ similar results were obtained. The `CheckDiverg' flag is included as \Cref{alg:PSDN} might lead to diverging estimates. When `CheckDiverg=True', an estimate of the standard deviation of the noise, i.e., $\sigma_k$ for $k=1,2,...,m$, must be specified.

\subsection{Finding derivatives and coefficients simultaneously}\label{sec:methods-socp}

To find the state time derivatives while enforcing sparsity of the coefficients, we write the coefficients as a function of the derivative. We start with the equation of the system dynamics, $\Theta^* \bm{c} = \dot{\bm{u}}^*$, and multiply both sides of this equation by the transpose of $\tilde{\Theta}$ (the monomial library evaluated at the smoothed states) to obtain the oblique projection equation, $\tilde{\Theta}^T \Theta^* \bm{c} = \tilde{\Theta}^T \dot{\bm{u}}^*$.
Then, assuming $\tilde{\Theta}^T \Theta^*$ is invertible, we have that
\begin{equation}\label{eq:c_eq_fudot}
  \bm{c} = G^{-1} \tilde{\Theta}^T \dot{\bm{u}}^* \quad \text{where} \quad G := \tilde{\Theta}^T \Theta^*.
\end{equation}
As $G$ cannot be directly evaluated, we explored two methods for estimating this matrix. First, we derived an estimator $\hat{G}$ such that $\hat{G}/N$ is consistent under certain assumptions (see \Cref{sec-app:unbiased-G-hat}). Second, we approximated $\Theta^*$ as $\tilde{\Theta}$ and $G$ as
\begin{equation}
  \tilde{G} := \tilde{\Theta}^T \tilde{\Theta}.
\end{equation}
Since we did not observe a significant difference between these two approaches,
we only show results for the second method in \Cref{sec:results}.  However, for alternative systems or sample sizes the choice of the estimator of $G$ may play a more significant role.

\begin{remark}
  The matrix $\tilde{\Theta}$ could be any approximated version of $\Theta^*$, e.g., found using state variable estimates from Gaussian process regression. However, in practice we obtain significantly better results when $\tilde{\Theta}$ is obtained using PSDN or IterPSDN.
\end{remark}

To find the state time derivatives we solve a convex optimization problem, where the objective enforces sparsity of the coefficients and the constraints require the derivative be smooth and agree with the data. Similar to iteratively reweighted Lasso (IRW-Lasso), see \Cref{sec-app:bg-reg}, we solve this problem multiple times, where, at each iteration, we weight the coefficients based on their previous values. Because this problem can be formulated as a SOCP, we refer to the algorithm as IRW-SOCP.

More specifically, for $i=1,2,...$, we find $\dot{\bm{u}}^{(i)}$ and the corresponding coefficient estimate $\bm{c}^{(i)} = \tilde{G}^{-1} \tilde{\Theta}^T\dot{\bm{u}}^{(i)}$ by solving,
\begin{equation}\label{eq:SOCP}
  \begin{aligned}
    \text{minimize}_{u_0,\dot{\bm{u}}}                           &  &                                   &
    \|W^{(i-1)} \tilde{G}^{-1} \tilde{\Theta}^T \dot{\bm{u}}\|_1 &  & \text{(sparsity of coefficients)}   \\
    \text{subject to}                                            &  &                                   &
    \|D\dot{\bm{u}}\| \le C                                      &  & \text{(smooth derivative)}          \\
                                                                 &  &                                   &
    \left\|\begin{bmatrix} \bm{1} & T \end{bmatrix}
    \begin{bmatrix} u_0 \\ \dot{\bm{u}} \end{bmatrix}
    - P_{\tilde{\Phi}}\tilde{\bm{u}}
    \right\| \le \gamma_i                                        &  & \text{(match smoothed data).}       \\
  \end{aligned}
\end{equation}
For the first iteration, $W^{(0)} = I$ and for $i > 0$, $W^{(i)}$ is a diagonal matrix where
\begin{equation}
  W^{(i)}_{jj} =
  \frac{1}{|c^{(i)}_j|+ \varepsilon \max_j |c^{(i)}_j|}.
\end{equation}
In practice, we set $\varepsilon=10^{-4}$. The finite difference matrix $D$ and discrete integral operator $T$ are as presented in \Cref{sec:notation}. In the second constraint, we include a projection of $\tilde{u}$ onto the column space of $\tilde{\Phi}$ in order to guarantee feasibility (see \Cref{sec:methods-C}).

There are two hyperparameters in \Cref{eq:SOCP} that need to be determined at each iteration, $C$ and $\gamma_i$. In practice, we set $C$ to guarantee feasibility at small values of $\gamma_i$ (see \Cref{sec:methods-C}). The remaining parameter, $\gamma_i$, is determined either using the Pareto corner criteria or by directly using a theoretical value (see \Cref{sec:gamma}). Both approaches for finding $\gamma_i$ rely on theoretical results given in \Cref{sec:theory}. We include a subscript on $\gamma$ to emphasize that its value may change at each iteration.

\subsubsection{Setting the smoothing hyperparameter \texorpdfstring{$C$}{} to guarantee feasibility}
\label{sec:methods-C}

We estimate $C$ in  \Cref{eq:SOCP} using the PSDN results. Specifically, we find an initial estimate of the derivative,
\begin{equation}\label{eq:tildeudot}
  \dot{\tilde{\bm{u}}} := \tilde{\Theta}\left(
  (\tilde{\Phi}^T \tilde{\Phi})^{-1} \tilde{\Phi}^T
  \right)_{2:,:}\tilde{\bm{u}},
\end{equation}
where the `$2:,:$' subscript implies that we are removing the top row of the matrix.
We then use this estimate to set the value of $C$ as,
\begin{equation}\label{eq:C}
  C := \|D \dot{\tilde{\bm{u}}}\|.
\end{equation}
At this $C$, \Cref{eq:SOCP} is feasible for all $\gamma > 0$ (see \Cref{sec-app:smoothing-hyper-param}).

\begin{remark}
  In practice, the $C$ found with \Cref{eq:C} was often larger than the optimal value. However, generally the results were not very sensitive to this parameter.
\end{remark}

\subsubsection{Setting the data-matching hyperparameter \texorpdfstring{$\gamma$}{}}\label{sec:gamma}

Consider the data-matching constraint in \Cref{eq:SOCP} where $\dot{\bm{u}}$ is replaced with the true time derivative vector, i.e., $\dot{\bm{u}}^*$. The constraint is then approximately
\begin{equation}
  \|\bm{u}^* - P_{\tilde{\Phi}} \tilde{\bm{u}}\| \le \gamma,
\end{equation}
where here we are ignoring quadrature error.
In \Cref{sec:theory}, we show that at large N the best possible expected value of this data-matching error is given as,
\begin{equation}\label{eq:gamma_exp}
  \mathbb{E}\left[\|\bm{u}^* - P_{\Phi^*} \bm{u}\|\right] \approx
  \sigma \sqrt{p+1} =: \gamma_{exp},
\end{equation}
where $\sigma$ is the standard deviation of the measurement noise. Although this optimal value will not necessarily be obtained by PSDN, in \Cref{sec:results-PSDN-theory}, we show that IterPSDN does obtain this value for the systems considered.
Therefore, we predict the optimal value of $\gamma_i$, for $i=1,2,...$, will be close to this estimate, i.e., close to $\gamma_{exp}$.

With this in mind, we consider two methods for finding the data-matching hyperparameter. First, we use the corner point of the Pareto curve, see \Cref{sec-app:Pareto} for details, but restrict the search space to lie within an order of magnitude of $\gamma_{exp}$, i.e.,
\begin{equation}\label{eq:pareto-limits}
  \gamma_i \in [0.1 \cdot \gamma_{exp}, 10 \cdot \gamma_{exp}]
  \quad\quad \text{for} \quad\quad i=1,2,....
\end{equation}
In practice, we find that this improves the performance of the corner finding algorithm and a clear corner point exists within this range. This is the approach used in the main results section, i.e. \Cref{sec:results}. We also explore the option of setting $\gamma_i = \gamma_{exp}$ for all $i$. In \Cref{sec-app:gamma-method}, we compare the performance of these two methods.

\begin{remark}
  By using the results from PSDN or IterPSDN in \Cref{eq:SOCP}, we predict a hyperparameter $\gamma_{exp}$ that is independent of $N$. If we had instead used the original measurements, the data-matching error $\|\bm{u}^* - \bm{u}\|$ would grow with the square root of $N$, leading to a less restrictive constraint.
\end{remark}

\section{Theoretical work}\label{sec:theory}

In this section we derive error bounds for PSDN as given by \Cref{eq:tildeu}. We show that, in the limit as the number of measurements goes to infinity and the training time is held constant, the projected state variables approach their true values with respect to the relative $\ell_2$ error. For notational simplicity we introduce the following definitions:
\begin{equation}
    \begin{aligned}
        \Psi^*      & := \frac{1}{\sqrt{N}} \Phi^*, & \Delta \Theta       & := \hat{\Theta} - \Theta^*,     \\
        \Delta \Phi & := \hat{\Phi} - \Phi^*,       & \quad\quad \Delta P & := P_{\hat{\Phi}} - P_{\Phi^*},
    \end{aligned}
\end{equation}
where $\Theta^*$, $\Phi^*$, and $\hat{\Phi}$ are as given by \Cref{eq:Theta,eq:b_and_Phistar,eq:Phi}. The matrix $\hat{\Theta}$ is an unbiased estimate of the noisy monomial library (see \Cref{sec-app:centered-library}), and the projection operator $P$ is as defined in \Cref{eq:P_A}. Throughout this section we drop the subscript on $\bm{u}_k$ as the results can be applied to any of the $m$ states.

\subsection{An optimal error estimate for PSDN}

We obtain an optimal error estimate by calculating the expected squared error of projecting the noisy data onto the column space of $\Phi^*$. This result informs upon the best possible performance of PSDN, see \Cref{eq:tildeu}. Additionally, these results can be used to find the value of the data-matching hyperparameter $\gamma$ in IRW-SOCP, i.e., \Cref{eq:SOCP}. We first derive an upper bound on the expected squared error of $\tilde{u}^*$ as given by \Cref{eq:tildeustar}.

\begin{lemma}
    \label{lemma:err_bound_1}
    Suppose the function $u^*(t)$ is three times continuously differentiable for all $t \in [0,t_{end}]$, then
    \begin{equation}
        \mathbb{E}\left[\|P_{\Phi^*} \bm{u} - \bm{u}^*\|^2
        \right]
        \le
        \sigma^2 (p+1)  + \frac{C_1^2}{(N-1)^3},
    \end{equation}
    where
    \begin{equation}
        C_1 = \frac{t_{end}^3}{12}\max_{t\in[0,t_{end}]} \left|\frac{d^3 u^*}{dt^3}\right|.
    \end{equation}
\end{lemma}

\begin{proof}
    First let
    \begin{equation*}
        \bm{b} := \begin{bmatrix} u_0^* \\ \bm{c} \end{bmatrix}, \quad
        \bm{e}_q := \Phi^* \bm{b} - \bm{u}^*, \quad
        P^{\perp}_{\Phi^*} := I - P_{\Phi^*},
    \end{equation*}
    where $\bm{e}_q \in \mathbb{R}^N$ is the trapezoidal quadrature error, see \Cref{eq:quad-err}. Then the error of projecting $\bm{u}^*$ onto the column space of $\Phi^*$ is given as,
    \begin{equation}\label{eq-app:err1}
        \begin{aligned}
            P_{\Phi^*} \bm{u}^* - \bm{u}^*
            = P_{\Phi^*} \left( \Phi^* \bm{b}  - \bm{e}_q \right) - \bm{u}^* & = \Phi^* \bm{b} - P_{\Phi^*} \bm{e}_q - \bm{u}^*                                                   \\
                                                                             & = \left(\bm{e}_q + \bm{u}^*\right) - P_{\Phi^*} \bm{e}_q - \bm{u}^* = P_{\Phi^*}^{\perp} \bm{e}_q.
        \end{aligned}
    \end{equation}
    In turn, the error of projecting $\bm{u} = \bm{u}^* + \bm{\epsilon}$ onto the column space of $\Phi^*$ is
    \begin{equation}\label{eq-app:err2}
        P_{\Phi^*} \bm{u} - \bm{u}^* = P_{\Phi^*} \bm{u}^*  + P_{\Phi^*} \bm{\epsilon} - \bm{u}^* = P_{\Phi^*} \bm{\epsilon} + P_{\Phi^*}^{\perp} \bm{e}_q,
    \end{equation}
    where we use \Cref{eq-app:err1} to obtain the second equality. Since $P_{\Phi^*}$ is an orthogonal projection, the norm of the right hand side of \Cref{eq-app:err2} is
    \begin{equation}\label{eq-app:err3}
        \|P_{\Phi^*} \bm{\epsilon} + P_{\Phi^*}^{\perp} \bm{e}_q\|^2
        = \|P_{\Phi^*} \bm{\epsilon} \|^2 + \|P_{\Phi^*}^{\perp} \bm{e}_q\|^2.
    \end{equation}
    We combine \Cref{eq-app:err2,eq-app:err3} and take the expectation to obtain,
    \begin{equation}
        \begin{aligned}
            \mathbb{E}\left[
            \|P_{\Phi^*} \bm{u} - \bm{u}^*\|^2
            \right]
            = \mathbb{E}\left[\bm{\epsilon}^T P_{\Phi^*} \bm{\epsilon}\right] + \|P_{\Phi^*}^{\perp} \bm{e}_q\|^2
            = \sum_{i,j} \mathbb{E}[\epsilon_{i} \epsilon_{j}] (P_{\Phi^*})_{ij} + \|P_{\Phi^*}^{\perp} \bm{e}_q\|^2 \\
            = \sum_{i} \mathbb{E}[\epsilon_{i}^2] (P_{\Phi^*})_{ii} + \|P_{\Phi^*}^{\perp} \bm{e}_q\|^2
            = \sigma^2 \Tr (P_{\Phi^*}) + \|P_{\Phi^*}^{\perp} \bm{e}_q\|^2.
        \end{aligned}
    \end{equation}
    The trace of a projection matrix is equal to the dimension of the target space, i.e., $\Tr (P_{\Phi^*}) = \Rank (\Phi^*)$, which implies
    \begin{equation}\label{eq-app:opt-err}
        \mathbb{E}\left[
        \|P_{\Phi^*} \bm{u} - \bm{u}^*\|^2
        \right] = \sigma^2 \Rank (\Phi^*) + \| P_{\Phi^*}^{\perp}\bm{e}_q\|^2 \le \sigma^2 (p+1)  + \|\bm{e}_q\|^2.
    \end{equation}
    The lemma is then proven by applying a bound (see \Cref{lemma:quad_error} in \Cref{sec-app:theory}) on the norm of the quadrature error, i.e.,
    \begin{equation}
        \|\bm{e}_q\| \le C_1 (N-1)^{-3/2}.
    \end{equation}
    This bound follows almost immediately from the well known trapezoidal rule error bound.

\end{proof}

\Cref{lemma:err_bound_1} implies that the expected squared relative error of projecting the noisy measurements $\bm{u}$ onto the column space of $\Phi^*$, i.e., $\mathbb{E}[\|P_{\Phi^*} \bm{u} - \bm{u}^*\|^2]/\|\bm{u}^*\|^2$, goes to zero as we increase the number of nonzero measurements. Since $C_1$ is proportional to $t_{end}^3$, this result only applies when the training time is held constant.

Next, we present a lower bound on the expected squared error. This bound follows immediately by considering \Cref{eq-app:opt-err} in the proof of \Cref{lemma:err_bound_1}.
\begin{corollary}\label{cor:err_bound_lower}
    If the assumptions of \Cref{lemma:err_bound_1} hold and $\Phi^*$ has full column rank, i.e., $\text{rank}(\Phi^*)=p+1$, then
    \begin{equation}
        \mathbb{E}\left[
        \|P_{\Phi^*} \bm{u} - \bm{u}^*\|^2
        \right]
        \ge
        \sigma^2 (p+1).
    \end{equation}
\end{corollary}
Taken together \Cref{lemma:err_bound_1} and \Cref{cor:err_bound_lower} imply,
\begin{equation}
    \lim_{N \rightarrow \infty}
    \mathbb{E}\left[
    \|P_{\Phi^*} \bm{u} - \bm{u}^*\|^2
    \right] = \sigma^2 (p+1).
\end{equation}

\begin{remark}
    Since $\mathbb{E}\left[X\right]^2 < \mathbb{E}\left[X^2\right]$,
    \Cref{lemma:err_bound_1} provides an upper bound on the expected relative error.
    Although, we have not provided a lower bound on this error, we still find that $\sigma \sqrt{p+1}$ provides a reasonable estimate (see \Cref{sec:results-PSDN-theory}).
\end{remark}

\subsection{An upper bound on the expected error of PSDN}

Our main result provides a bound on expected value of the PSDN error. Under certain assumptions, this result shows that as the sample size increases the relative error of PSDN goes to zero.
\begin{theorem}\label{thm:main}
    Let $\hat{\Phi}$ be given by \Cref{eq:Phi} where $T$ is the discrete integration matrix that performs trapezoidal quadrature. Suppose $\textnormal{rank}(\hat{\Phi}) = \textnormal{rank}(\Phi^*)$, $\|(\Phi^*)^\dagger\| \|\Delta \Phi\| < 1/4$, and $u^*(t)$ is three times continuously differentiable for $t \in [0,t_{end}]$. Then the expected relative error of PSDN is bounded as follows
    \begin{equation}\label{eq:main-bound}
        \mathbb{E}\left[
        \frac{\|P_{\hat{\Phi}} \bm{u} - \bm{u}^*\|}{\|\bm{u}^*\|}
        \right] \le \frac{\sigma \left(\sqrt{p+1} + C_2\right)}{\|\bm{u}^*\|}  + \frac{C_1}{(N-1)^{3/2}\|\bm{u}^*\|} + \frac{C_2}{\sqrt{N}},
    \end{equation}
    where $C_1$ is as defined in \Cref{lemma:err_bound_1} and
    \begin{equation}
        C_2 = t_{end} \|(\Psi^*)^\dagger\| \sum_{j=1}^p \max_k  \mathbb{V}\left[ \Delta \Theta_{jk}\right]^{1/2}.
    \end{equation}
\end{theorem}

A few comments on \Cref{thm:main} are warranted. In \Cref{sec:psi-pinv-norm-bound} we show $\|(\Psi^*)^\dagger\|$ is uniformly bounded over $N$, which implies the bound given by \Cref{eq:main-bound} goes to zero as the measurement density increases. The first assumption of \Cref{thm:main}, i.e., $\textnormal{rank}(\hat{\Phi}) = \textnormal{rank}(\Phi^*)$, is valid if the library $\Theta^*$ has full column rank. This is true for the systems examined here, see \Cref{sec:ode-systems}, and we  leave it as future work to examine the performance of PSDN when this is not the case.
In \Cref{cor:assumption-2}, given in \Cref{sec:thm-proof}, we show that the second assumption, i.e., $\|(\Phi^*)^\dagger\| \|\Delta \Phi\| < 1/4$, is satisfied in expectation at large $N$.

We provide the proof of \Cref{thm:main} in \Cref{sec:thm-proof}. In the proof we use the triangle and Cauchy-Schwarz inequalities to show that,
\begin{equation}\label{eq:exp-bound}
    \mathbb{E}\left[\|P_{\hat{\Phi}} \bm{u} - \bm{u}^*\|\right]
    \le \mathbb{E}\left[\|P_{\Phi^*} \bm{u} - \bm{u}^*\|\right]
    + \left(\mathbb{E}\left[\|\Delta P\|^2\right]\right)^{1/2}
    \left(\|\bm{u}^*\| + \sqrt{N}\sigma\right).
\end{equation}
Since \Cref{lemma:err_bound_1} provides a bound for the first term, it remains to show the expected value of $\|\Delta P\|^2$ goes to zero as $N \rightarrow \infty$, which we do next in \Cref{sec:DeltaP-bound}.

\subsubsection{Bounding the expected projection error \texorpdfstring{$\|\Delta P\|$}{}}\label{sec:DeltaP-bound}

We obtain a bound for the expected projection error using the results of two lemmas.
In \Cref{lemma:DeltaPhi_bound}, we derive a bound for $\mathbb{E}\left[\|\Delta \Phi\|\right]$ that depends on the elementwise variance of $\Delta \Theta = \hat{\Theta} - \Theta^*$. We then, in \Cref{lemma:DeltaP_bound}, derive a bound for  $\mathbb{E}\left[\|\Delta P\|\right]$ that depends on $\mathbb{E}\left[\|\Delta \Phi\|\right]$.

\begin{lemma}\label{lemma:DeltaPhi_bound}
    The expected error of $\hat{\Phi}$ is bounded as follows:
    \begin{equation}
        \mathbb{E}\left[\|\Delta \Phi\|^2\right] \le \frac{t_{end}^2}{2} \sum_{j=1}^p \max_k \left( \mathbb{V}\left[ \Delta \Theta_{kj}\right]\right).
    \end{equation}
\end{lemma}

\begin{proof}
    Using $\Delta \Phi = \hat{\Phi} - \Phi^* = \begin{bmatrix}
            \bm{0} & T\Delta \Theta
        \end{bmatrix}$
    where $T$ is given by \Cref{eq:T}, we have that
    \begin{equation}\label{eq:DeltaPhi_F}
        \|\Delta \Phi\|_F^2 = \|T\Delta \Theta \|_F^2 = \sum_{j=1}^p \sum_{i=2}^N  \left(\sum_{k=1}^i T_{ik} \Delta \Theta_{kj}\right)^2.
    \end{equation}
    We use this result, i.e., \Cref{eq:DeltaPhi_F}, to bound the expected value of $\|\Delta \Phi\|_F^2$,
    \begin{equation*}
        \begin{aligned}
            \mathbb{E}\left[\|\Delta \Phi\|_F^2\right]
             & = \sum_{j=1}^p \sum_{i=2}^N \mathbb{E}\left[\left( \sum_{k=1}^i  T_{ik}
                \Delta \Theta_{kj}\right)^2\right]
            = \sum_{j=1}^p \sum_{i=2}^N \mathbb{V}\left[\sum_{k=1}^i  T_{ik}
            \Delta \Theta_{kj}\right]                                                  \\
             & = \sum_{j=1}^p \sum_{i=2}^N \sum_{k=1}^i  T^2_{ik}
            \mathbb{V}\left[\Delta \Theta_{kj} \right]
            \le \Delta t^2 \sum_{j=1}^p \sum_{i=2}^N \left(i - \frac{3}{2}\right) \max_k \left(\mathbb{V}\left[\Delta \Theta_{kj}\right]\right).
        \end{aligned}
    \end{equation*}
    The inequality follows because the sum of the squared elements in the $i$th row of $T$ is
    \begin{equation*}
        \sum_{k=1}^i T_{ik}^2 = \frac{\Delta t^2}{4}(4i - 6) = \Delta t^2 \left(i - \frac{3}{2}\right).
    \end{equation*}
    Continuing, we have that
    \begin{equation*}
        \begin{aligned}
            \mathbb{E}\left[\|\Delta \Phi\|_F^2\right]
             & \le \Delta t^2 \sum_{j=1}^p \sum_{i=1}^{N-1} \left(i-\frac{1}{2}\right) \max_k \left(\mathbb{V}\left[ \Delta \Theta_{kj}\right]\right) \\
             & = \Delta t^2 \frac{(N-1)^2}{2} \sum_{j=1}^p \max_k \left(\mathbb{V}\left[ \Delta \Theta_{kj}\right]\right)
            = \frac{t_{end}^2}{2} \sum_{j=1}^p \max_k \left(\mathbb{V}\left[ \Delta \Theta_{kj}\right]\right).
        \end{aligned}
    \end{equation*}
    The proof is completed by noting that the 2-norm of a matrix is always less than the Frobenius norm.
\end{proof}

The following lemma provides a bound on $\|\Delta P\|$ that depends on $\|\Delta \Phi\|$. In order to obtain this bound we use results from \cite{Stewart1977} as presented in \Cref{sec-app:P-err}.

\begin{lemma}\label{lemma:DeltaP_bound}
    Suppose $\textnormal{rank}(\hat{\Phi}) = \textnormal{rank}(\Phi^*)$ and $\|(\Phi^*)^\dagger\| \|\Delta \Phi\| < 1/4$, then
    \begin{equation}
        \|\Delta P\|^2
        \le 2 \|(\Phi^*)^\dagger\|^2 \|\Delta \Phi\|^2.
    \end{equation}
\end{lemma}

\begin{proof}
    Let $E_{11}$, $E_{21}$, and $\beta$ be as defined in \Cref{thm:stewart} (a modified version of Theorem~4.1 from \cite{Stewart1977}) and define the following scalars,
    \begin{equation*}
        X := \|E_{11}\|\|(\Phi^*)^\dagger\|, \quad \quad
        Y := \|E_{21}\|\|(\Phi^*)^\dagger\|.
    \end{equation*}
    First note that,
    \begin{equation*}
        \beta \|E_{21}\|/\|\Phi^*\| = \frac{\|E_{21}\| \|(\Phi^*)^\dagger\|}{1-\|E_{11}\|\|(\Phi^*)^\dagger\|} = \frac{Y}{1-X}.
    \end{equation*}
    Using $\|E_{11}\| < \|\Delta \Phi\|$ and the lemma assumptions, we have that,
    \begin{equation*}
        X = \|E_{11}\|\|(\Phi^*)^\dagger\| < \|\Delta \Phi\|\|(\Phi^*)^\dagger\| < 1/4.
    \end{equation*}
    We next square the bound given by \Cref{eq:bound_stewart} and rewrite it as,
    \begin{equation}\label{eq:DeltaP_bound}
        \|\Delta P\|^2 \le \frac{Y^2/(1-X)^2}{1 + Y^2/(1-X)^2}
        = \frac{Y^2}
        {(1-X)^2+Y^2} \le 2 Y^2.
    \end{equation}
    The second inequality follows because, for $X<1/4$,
    \begin{equation*}
        2 Y^2((1-X)^2 + Y^2) = 2 Y^2(1 - 2X + X^2 + Y^2) \ge  2 Y^2(1 - 2X) \ge Y^2.
    \end{equation*}
    The lemma statement follows from \Cref{eq:DeltaP_bound} since $Y \le \|(\Phi^*)^\dagger\| \|\Delta \Phi\|$.
\end{proof}

\subsubsection{Proof of the main result} \label{sec:thm-proof}

We first state a corollary that shows an assumption of \Cref{thm:main} is satisfied at large $N$. This corollary follows immediately from \Cref{lemma:DeltaPhi_bound} and \Cref{cor:Psi_pseudo_bound}, which shows the value of $\| (\Psi^*)^{\dagger} \|$ is bounded from above for all $N$.
\begin{corollary}\label{cor:assumption-2}
    Suppose the assumptions of \Cref{cor:Psi_pseudo_bound} in \Cref{sec-app:theory} hold, i.e., $u_k(t)$ for $k=1,2,..,m$ is twice continuously differentiable for $t \in [0,t_{end}]$ and in the limit as $N \rightarrow \infty$ the rank of $(\Psi^*)^T \Psi^*$ does not change. Then there exists $N_0$ such that for $N>N_0$
    \begin{equation}
        \|(\Phi^*)^\dagger\| \mathbb{E}\left[\|\Delta \Phi\|\right] < 1/4.
    \end{equation}
\end{corollary}
We next provide the proof of \Cref{thm:main} using the results of
\Cref{lemma:err_bound_1,lemma:DeltaPhi_bound,lemma:DeltaP_bound}.
\begin{proof}[Proof of \Cref{thm:main}]
    We first bound the value of $\|P_{\hat{\Phi}} \bm{u} - \bm{u}^*\|$ using the triangle inequality
    \begin{equation}\label{eq-app:pf-bound1}
        \begin{aligned}
            \|P_{\hat{\Phi}} \bm{u} - \bm{u}^*\|
             & \le \|P_{\Phi^*} \bm{u} - \bm{u}^*\| + \|\Delta P\|\|\bm{u}^* + \bm{\epsilon}\|                   \\
             & \le \|P_{\Phi^*} \bm{u} - \bm{u}^*\| + \|\Delta P\|\left(\|\bm{u}^*\| + \|\bm{\epsilon}\|\right),
        \end{aligned}
    \end{equation}
    where $\Delta P = P_{\hat{\Phi}} - P_{\Phi^*}$. Taking expectations of both sides of \Cref{eq-app:pf-bound1} and using Jensen's inequality, we have that
    \begin{equation}\label{eq:allterms}
        \mathbb{E}\left[\|P_{\hat{\Phi}} \bm{u} - \bm{u}^*\|\right]
        \le T_1 + T_2 + T_3,
    \end{equation}
    where
    \begin{equation*}
        T_1 := \sqrt{\mathbb{E}\left[\|P_{\Phi^*} \bm{u} - \bm{u}^*\|^2\right]}, \quad\quad
        T_2 := \|\bm{u}^*\|\sqrt{\mathbb{E}\left[\|\Delta P\|^2\right]}, \quad\quad
        T_3 := \sqrt{\mathbb{E}[\|\Delta P\|^2 \|\bm{\epsilon}\|^2]}.
    \end{equation*}

    We first obtain a bound for $T_1$ using \Cref{lemma:err_bound_1} and the triangle inequality,
    \begin{equation}\label{eq:term1}
        T_1 = \sqrt{\mathbb{E}\left[\|P_{\Phi^*} \bm{u} - \bm{u}^*\|^2\right]}
        \le \sigma \sqrt{p+1}  + \frac{\sqrt{C}_1}{(N-1)^{3/2}}.
    \end{equation}
    To bound $T_2$ and $T_3$ we first use \Cref{lemma:DeltaP_bound,lemma:DeltaPhi_bound} to bound $\mathbb{E}\left[\|\Delta P\|^2\right]$,
    \begin{equation}\label{eq:term3}
        \mathbb{E}\left[\|\Delta P\|^2\right]
        \le 2 \|(\Phi^*)^\dagger\|^2 \mathbb{E} \left[\|\Delta \Phi\|^2\right]
        \le \frac{\|(\Psi^*)^\dagger\|^2}{N} t_{end}^2 \sum_{j=1}^p \max_k \left( \mathbb{V}\left[ \Delta \Theta_{kj}\right]\right).
    \end{equation}
    \Cref{eq:term3} can immediately be used to bound $T_2$. For $T_3$ we first use Cauchy Schwarz to show
    \begin{equation}\label{eq:term2}
        T_3^2 = \mathbb{E}\left[\|\Delta P\|^2 \|\bm{\epsilon}\|^2\right]
        \le
        \mathbb{E}\left[\|\Delta P\|^2 \right]
        \mathbb{E}\left[\|\bm{\epsilon}\|^2\right] = \sigma^2N \mathbb{E}\left[\|\Delta P\|^2 \right].
    \end{equation}
    Putting \Cref{eq:allterms,eq:term1,eq:term2,eq:term3} together, we obtain the final result.
\end{proof}

\section{Methods for obtaining and comparing numerical results}\label{sec:ode-systems}

\subsection{Implementation details for equation discovery algorithms}

We compared the performance of DSINDy with multiple alternative methods, including Gaussian process (GP) regression \cite{Schulz2018}, $\ell_1$-SINDy \cite{Brunton2016,Cortiella2021}, and WSINDy \cite{Messenger2021}.
To perform WSINDy we used the most recent MATLAB code available on GitHub\footnote{\url{https://github.com/MathBioCU/WSINDy_ODE} (accessed 9/28/2022)}. In general, when running WSINDy we used the default parameters given in the MATLAB code. The one exception to this was for the Lorenz 96 model, where WSINDy did poorly at low noise levels. To fix this we set the parameter `alpha\_loss' = 0.4 (default: 0.8) and `overlap\_frac' = 0.1 (default: 1.0).
The python code for all other computations, i.e., GP regression, $\ell_1$-SINDy and DSINDy, is available on our group's GitHub page\footnote{\url{https://github.com/CU-UQ/DSINDy}}.

\Cref{tab:methods} provides a summary and comparison of the major steps behind each of the equation discovery algorithms. For all three methods we include \textit{a priori} denoising of the data. For $\ell_1$-SINDy and DSINDy the denoising is done with IterPSDN, and for WSINDy we used the built-in smoothing algorithm in the WSINDy MATLAB code. The $\ell_1$-SINDy algorithm obtains the state time derivatives using Tikhonov regularization and finds the basis coefficients using IRW-Lasso. In contrast, WSINDy does not require a derivative estimation and finds the coefficients using sequential-thresholding least squares (STLS). A complete mathematical description of these alternative methods is given in \Cref{sec-app:algorithms}

Our results rely on several open-source python packages. To perform GP regression, we used the gaussian\_process package from the scikit-learn python library \cite{Pedregosa2011} and constructed a GP with a radial basis function plus white noise kernel. To implement $\ell_1$-SINDy we used the scikit-learn implementation of Lasso, and for DSINDy the optimization problem given by \Cref{eq:SOCP} was solved using the python CVXOPT package \cite{Anderson2015} coupled with MOSEK Solver \cite{MosekApS2019}.

\begin{remark}
  A review of both local and global smoothing methods, in the context of equation discovery, is provided by \cite{Cortiella2022}. For the present work, we compare our denoising approach with GP regression, as a global denoising method. GP regression offers flexibility in fitting nonlinear functions and has available python implementations for finding the GP hyperparameters \cite{Pedregosa2011}.
\end{remark}
%
%%%%%%%%%% 
\begin{table}
  \centering
  \footnotesize
  \caption{Implementation details for $\ell_1$-SINDy, WSINDy, and DSINDy. In parentheses, we specify how hyperparameters were found. For more details see \Cref{sec-app:algorithms}.}
  \label{tab:methods}
  \begingroup
  \setlength{\tabcolsep}{8pt} % Default value: 6pt
  \renewcommand{\arraystretch}{1.2} % Default value: 1
  \begin{tabular}{l|l|l|l}
    \textbf{Method} & \textbf{Denoising} & \textbf{Derivative Estimation}                  & \textbf{Coefficient Estimation}           \\ \hline
    $\ell_1$-SINDy  & IterPSDN           & Tikhonov Reg. (Pareto)                          & IRW-Lasso (Pareto)                        \\ \hline
    WSINDy          & Weighted Avg.      & Not Applicable                                  & Modified STLS (see \cite{Messenger2021a}) \\ \hline
    DSINDy          & IterPSDN           & \multicolumn{2}{c}{IRW-SOCP (Pareto or theory)}
  \end{tabular}
  \endgroup
\end{table}
%%%%%%%%%%

\begin{table}
  \centering
  \footnotesize
  \caption{Summary of dynamic systems. For the Duffing oscillator, two dynamic parameter sets were considered (PS1 and PS2). For both sets, the same $t_{end}$, $d$ and initial conditions $\bm{u}_{IC}$ were used.}
  \label{tab}
  \begingroup
  \setlength{\tabcolsep}{8pt} % Default value: 6pt
  \begin{tabular}{lll}
    System                                                                                                      & Equations & Parameters \\ \hline \\[-1em]
    \begin{tabular}{@{}l} Duffing  \\ $(m=2)$ \end{tabular}                           &
    $\begin{aligned}
         \dot{u}_1 & = u_2                                       \\
         \dot{u}_2 & = -\kappa u_1 - \gamma u_2 - \epsilon u_1^3
       \end{aligned}$                                                  &
    $\begin{aligned}
          & \text{PS1: } \kappa=1, \gamma=0.1, \epsilon=5  \\
          & \text{PS2: }\kappa=0.2, \gamma=0.2, \epsilon=1 \\
          & t_{end}=10, d=4                                \\
          & \bm{u}_{IC}=[0, 1]
       \end{aligned}$                                                                                   \\ \hline \\[-1em]
    \begin{tabular}{@{}l} Van der Pol                        \\ $(m=2)$ \end{tabular} &
    $\begin{aligned}
         \dot{u}_1 & = u_2                                  \\
         \dot{u}_2 & = -u_1 + \gamma u_2 - \gamma u_2^3 u_2
       \end{aligned}$                                                       &
    $\begin{aligned}
          & \gamma=2           \\
          & t_{end}=10, d=4    \\
          & \bm{u}_{IC}=[0, 1]
       \end{aligned}$                                                                                                               \\ \hline \\[-1em]
    \begin{tabular}{@{}l} R\"{o}ssler               \\ $(m=3)$ \end{tabular}          &
    $\begin{aligned}
         \dot{u}_1 & = -u_2-u_3                     \\
         \dot{u}_2 & = u_1 + \alpha u_2             \\
         \dot{u}_3 & = \beta - \kappa u_3 + u_1 u_3
       \end{aligned}$                                                               &
    $\begin{aligned}
          & \alpha=0.2, \beta=0.2, \kappa=5.7 \\
          & t_{end}=10, d=2                   \\
          & \bm{u}_{IC}=[0,-5,0]
       \end{aligned}$                                                                                                \\ \hline \\[-1em]
    \begin{tabular}{@{}l} Lorenz 96 \\ $(m=6)$ \end{tabular}                          &
    $\begin{aligned}
         \dot{u}_i & = \left(u_{i+1} - u_{i-2}\right)u_{i-1} - u_i + F \\
                   & \text{for } i=1,2,...,m                           \\
                   & \text{where } u_{j-m}=u_{j+m}=u_j
       \end{aligned}$                                            &
    $\begin{aligned}
          & F=8                         \\
          & t_{end}=5, d=3              \\
          & \bm{u}_{IC} = [1,8,8,8,8,8]
       \end{aligned}$
  \end{tabular}
  \endgroup
\end{table}

\subsection{Example ODE systems}

The equations and parameters for the ODE systems we considered are given in \Cref{tab}. For both the Duffing and Van der Pol oscillator, the two state variables, $u_1$ and $u_2$, represent the displacement and velocity, respectively. For the Duffing oscillator, we considered two parameter sets, which both result in a damped system with a stable equilibrium. The first parameter set (PS1) explores the ability of DSINDy to recovery coefficients that vary significantly in magnitude, and the second parameter set (PS2) was used in the WSINDy work \cite{Messenger2021}.
We next considered the Van der Pol oscillator where, by setting $\gamma > 0$, we are guaranteed the system enters a stable limit cycle. This example demonstrates the ability of DSINDy to capture this cyclical behavior.
Finally, we considered two chaotic systems, the R\"{o}ssler attractor and the Lorenz 96 model. Using the given parameter set, the R\"{o}ssler attractor has a chaotic attractor and the third state $u_3$ often has values close to zero. We studied this system to explore how DSINDy performs in this more complex scenario. We additionally studied the Lorenz 96 system with $m=6$ states to examine how DSINDy performs on larger dimensional chaotic systems.

\subsection{Error metrics for comparing results}

In the results section, we report on the error of denoising and system recovery. Typically, the error is reported as a relative $\ell_2$-error. For example, the relative error of a coefficient vector estimate $\bm{c}$, given a true vector $\bm{c}^*$, is given as
\begin{equation}\label{eq:relerr}
  \mathcal{E}(\bm{c}; \bm{c}^*) := \frac{\|\bm{c}-\bm{c}^*\|}{\|\bm{c}^*\|}.
\end{equation}

For the Duffing oscillator, the Van der Pol oscillator, and the R\"{o}ssler attractor, we report the reconstruction error. This error is calculated by simulating the system with the same initial conditions as the training data but for double the training time, i.e., $t_{test} = 2 t_{end}$. For the Lorenz 96 model, we instead explore the ability of DSINDy to predict, new, short term dynamics. To do this, we use the same initial conditions as the training data and find the true system value at $t=t_{end}=5$. We use this as the initial condition for testing, i.e. $u_{IC,test} = u^*(t_{end})$, and simulate the learned dynamics using the initial condition $u_{IC,test}$. We record the time over which the simulated system has a large predictive ability, i.e., the time for which the relative reconstruction error is less than 10\% for all state variables. For all simulations, we record the state estimates using a time interval of $\Delta t = 0.01$ and compare the result to the true solution using \Cref{eq:relerr}.

If the relative reconstruction error exceeds 1.0 for any state in the system, then this is classified as a failed simulation, and we record the final error as 1.0. A system that becomes unbounded is also classified as a failure, and the reconstruction error is set to 1.0 for all state variables. This approach allows us to compare average reconstruction results for multiple methods and noise realizations.

In general, in figures showing the denoising and derivative estimation errors we provide the mean $\pm$ standard deviation (std) of 30 noise replications. However, in the figures showing the coefficient and reconstruction errors, we show the standard error of the mean (sem) instead of the standard deviation in order to make the plots more readable.

\section{Results}\label{sec:results}

We first compare numerical results for the projection-based denoising methods, i.e., PSDN and IterPSDN, with theoretical predictions (\Cref{sec:results-PSDN-theory}). We then examine the results of the complete DSINDy algorithm and compare the error of the denoising step with GP regression (\Cref{sec:results-PSDN-GP}) and the error of the IRW-SOCP derivative estimation with Tikhonov regularization (\Cref{sec:results-deriv}).
Ultimately, DSINDy improves upon coefficient recovery and state variable estimation compared with $\ell_1$-SINDy and WSINDy (\Cref{sec:results-system}).

\subsection{Projection denoising results agree with theoretical predictions}\label{sec:results-PSDN-theory}

\begin{figure}
    \centering
    \includegraphics{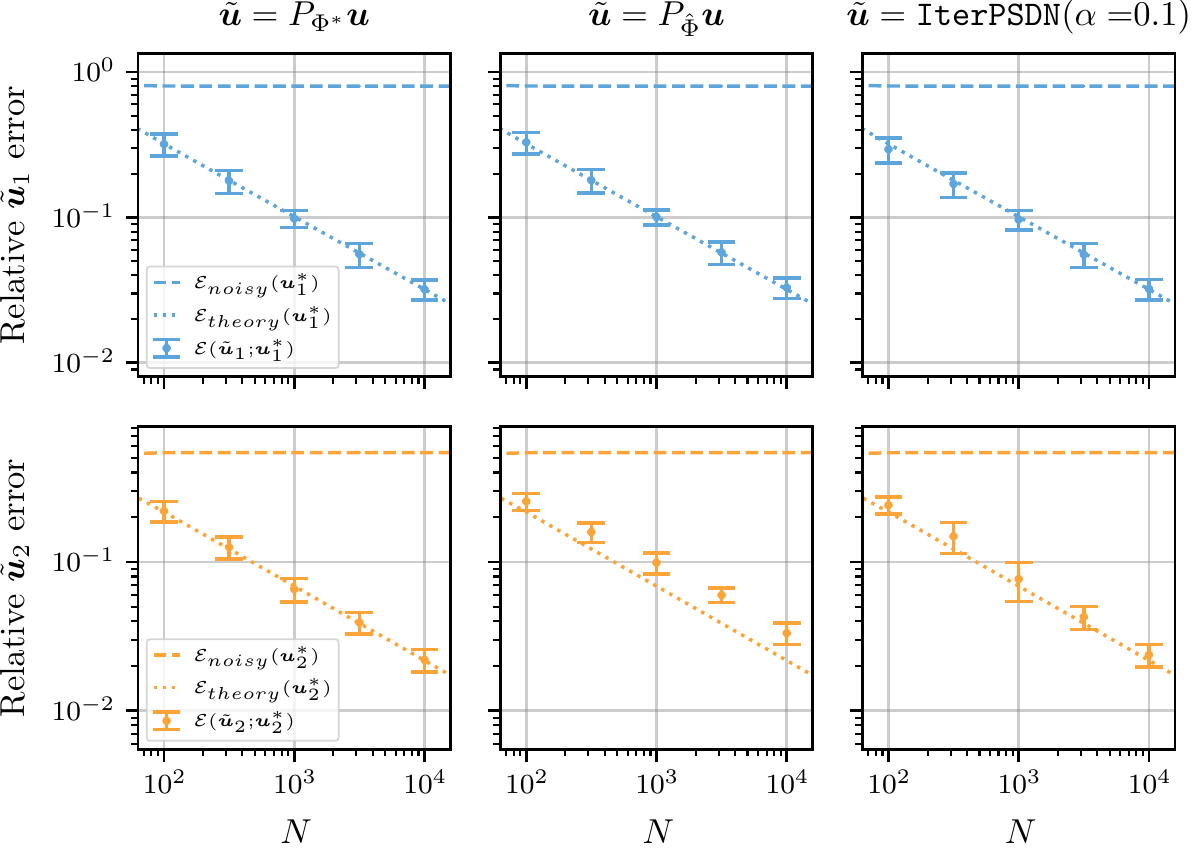}
    \caption{Comparing denoising error with theoretical predictions for the Duffing oscillator (PS1).  Results show mean $\pm$ std for 50 sample replications at noise level $\sigma^2=0.1$ when the integrated library $\Phi^*$ is known (left), when using PSDN (middle), and when using IterPSDN, \Cref{alg:PSDN} (right). Each row represents one of the system states. For system parameters see \Cref{tab}.}
    \label{fig:PSDN-Theory1}
\end{figure}

\begin{figure}
    \centering
    \includegraphics{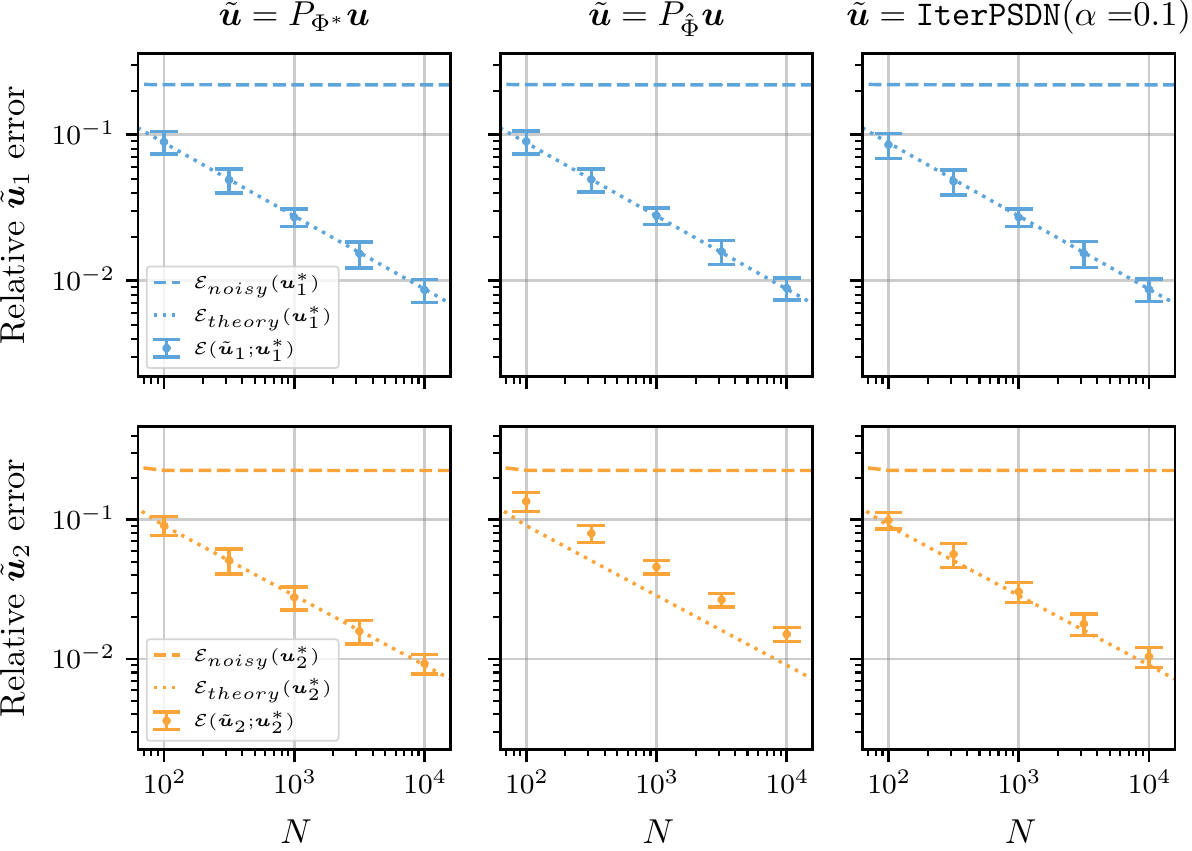}
    \caption{Comparing denoising error with theoretical predictions for the Van der Pol oscillator. Results show mean $\pm$ std for 50 sample replications at noise level $\sigma^2=0.1$ when the integrated library $\Phi^*$ is known (left), when using PSDN (middle), and when using IterPSDN, \Cref{alg:PSDN} (right). Each row represents one of the system states. For system parameters see \Cref{tab}.}
    \label{fig:PSDN-Theory2}
\end{figure}

\begin{figure}[!t]
    \centering
    \includegraphics{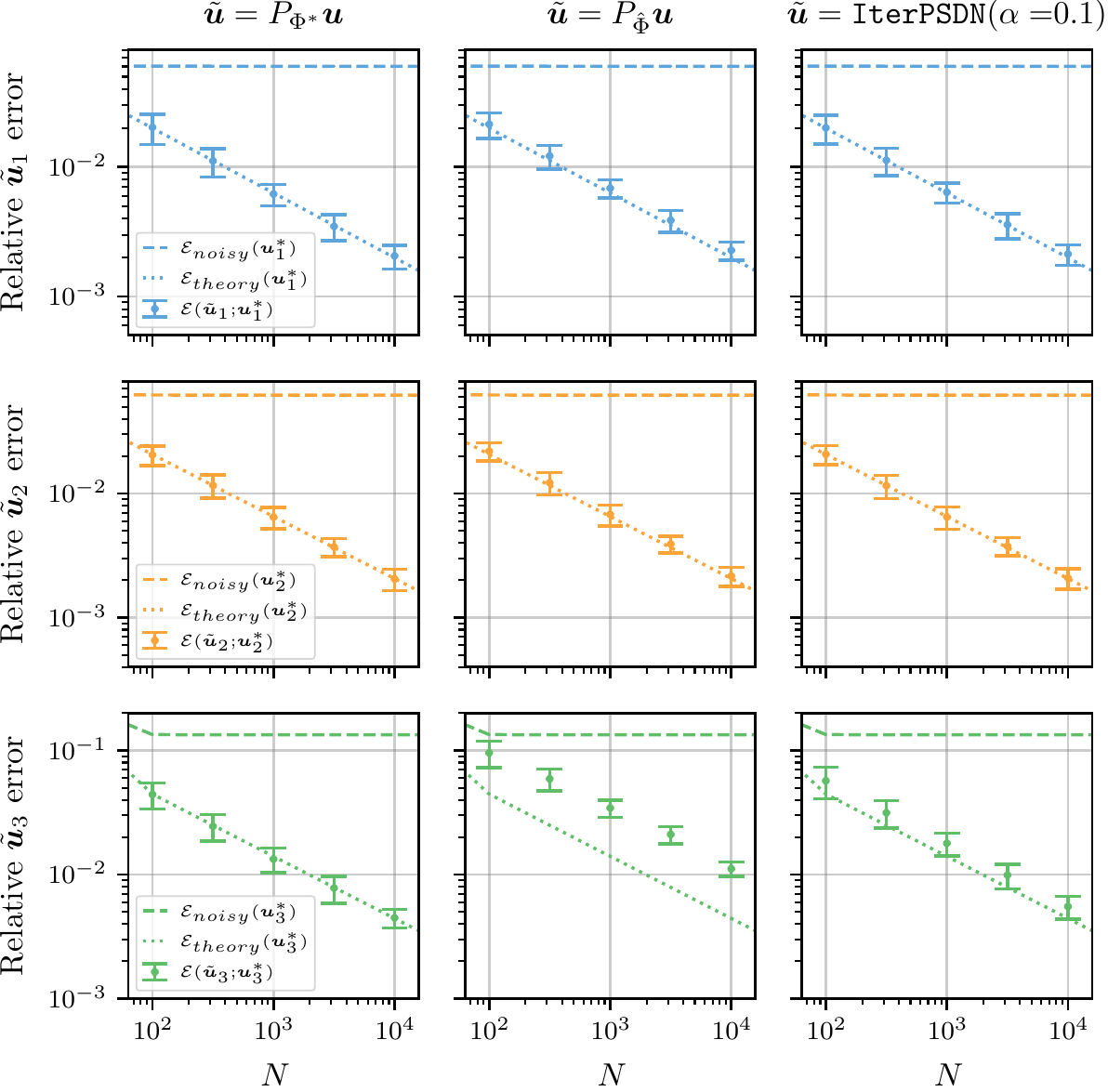}
    \caption{Comparing denoising error with theoretical predictions for the R\"{o}ssler attractor. Results show mean $\pm$ std for 50 sample replications at noise level $\sigma^2=0.1$ when the integrated library $\Phi^*$ is known (left), when using PSDN (middle), and when using IterPSDN, \Cref{alg:PSDN} (right). Each row represents one of the system states. For system parameters see \Cref{tab}.}
    \label{fig:PSDN-Theory3}
\end{figure}
%%%%

Here, we compare denoising results with the theoretical bounds presented in \Cref{sec:theory}. We compare the relative error of a denoised state vector, i.e., $\mathcal{E}(\tilde{\bm{u}}; \bm{u}^*)$, see \Cref{eq:relerr}, with the optimal error estimate derived in \Cref{sec:theory}, i.e.,
\begin{equation}
    \mathcal{E}_{theory}(\bm{u}^*) := \frac{\sigma\sqrt{p+1}}{\|\bm{u}^*\|}.
\end{equation}
This value represents the minimal error we can expect to obtain in PSDN or IterPSDN. For comparison, we also show the expected error of the original noisy data, i.e.,
\begin{equation}
    \mathcal{E}_{noisy}(\bm{u}^*) := \frac{\sigma\sqrt{N}}{\|\bm{u}^*\|}.
\end{equation}

We compare $\mathcal{E}$, $\mathcal{E}_{theory}$, $\mathcal{E}_{noisy}$ at a range of sample sizes $N$ for the four ODE systems presented in \Cref{tab} (\Cref{fig:PSDN-Theory1,fig:PSDN-Theory2,fig:PSDN-Theory3}). Results for the Lorenz 96 model are given in \Cref{sec-app:Lorenz96}. At each $N$, we performed 50 sample replications to approximate the mean and standard deviation of the state prediction error following denoising. Results are shown for three variations of the projection-state denoising approach: denoising when $\Phi^*$ is known (first column of \Cref{fig:PSDN-Theory1,fig:PSDN-Theory2,fig:PSDN-Theory3}), PSDN (second column), and IterPSDN (third column).
When applying \Cref{alg:PSDN}, we set the `CheckDiverg' flag equal to `True' and used the true standard deviation of the measurement noise. For the Duffing and Van der Pol oscillators and the Lorenz 96 model, divergence did not occur. However, for the R\"{o}ssler attractor, checking for divergence led to improved solutions.

\begin{remark}
    The standard deviation of the noise $\sigma$ is often known based on the accuracy of the measurement system or can be estimated. For example, in \Cref{sec:results-PSDN-GP}, $\sigma$ is estimated using the result of GP regression.
\end{remark}

For all four systems, we find that as expected (see \Cref{lemma:err_bound_1} and \Cref{cor:err_bound_lower}) the projection error when $\Phi^*$ is known matches the theoretical value. For PSDN, some state variables have an increased error compared to this optimal value, i.e., see the second state in the Duffing oscillator system (middle column, bottom row of \Cref{fig:PSDN-Theory1}), the second state of the Van der Pol oscillator (middle column, bottom row of \Cref{fig:PSDN-Theory2}), and the third state of the R\"{o}ssler attractor (middle column, bottom row of \Cref{fig:PSDN-Theory3}). However, IterPSDN improves upon PSDN and obtains the optimal value for the aforementioned states. The one exception to this is the third state of the R\"{o}ssler attractor where IterPSDN improves upon PSDN, but the error is still slightly larger than the optimal value. These results suggest that IterPSDN helps remove the error introduced in the estimation of $\Phi^*$.

\begin{remark}
    Although IterPSDN improves upon the performance of PSDN for the systems considered, this error reduction is not guaranteed. We leave a more in depth analysis and comparison of PSDN and IterPSDN to future work.
\end{remark}

\subsection{DSINDy Results}

We applied the complete DSINDy algorithm (see \Cref{tab:methods}) to the example systems in \Cref{sec:ode-systems} using the parameters given in \Cref{tab}. Unless otherwise noted, for all systems we consider a sample size of $N=1000$ and noise levels ranging from $\sigma=10^{-3}$ to $\sigma=1$, where $\sigma$ represents the standard deviation of the Gaussian measurement noise. Since the magnitude of each state is different, the signal-to-noise level also varies for each state and system.
The DSINDy results in this section were obtained by setting the hyperparameter $\gamma$ in \Cref{eq:SOCP} using the corner point of the Pareto curve where the search space was limited as shown in \Cref{eq:pareto-limits}. In general when applying \Cref{alg:PSDN}, we set the `CheckDiverg' flag equal to `True' and set the standard deviations $\sigma_k$ for $k=1,...,m$ using the estimate obtained from GP regression. For the Duffing oscillator results at $N=8000$, we ran IterPSDN with the `CheckDiverg=False'. Divergence was not an issue, and this avoided performing GP regression to find the standard deviation.

\subsubsection{The IterPSDN step of DSINDy improves state estimation}
\label{sec:results-PSDN-GP}

\begin{figure}[t]
    \centering
    \includegraphics[trim=0 25 0 0, clip]{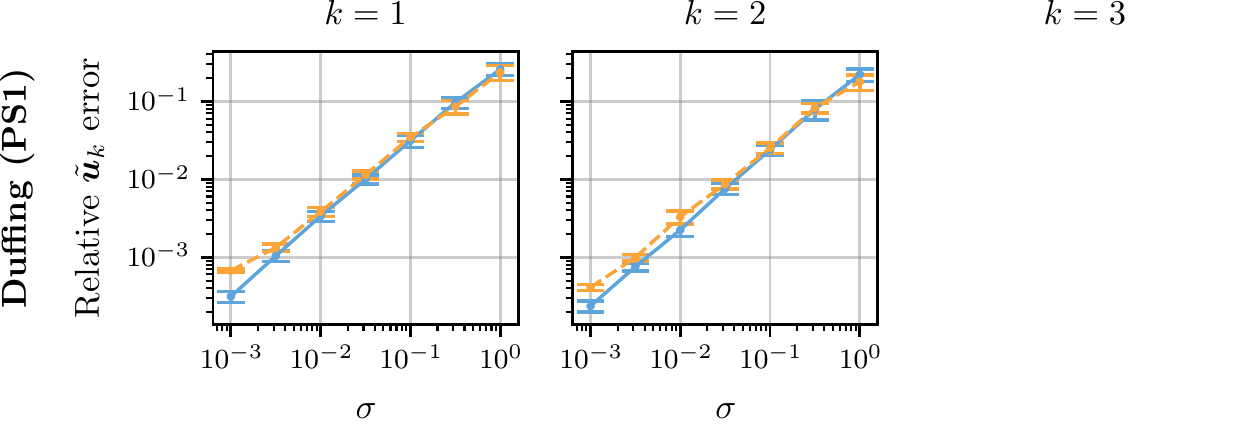}
    \includegraphics[trim=0 25 0 10, clip]{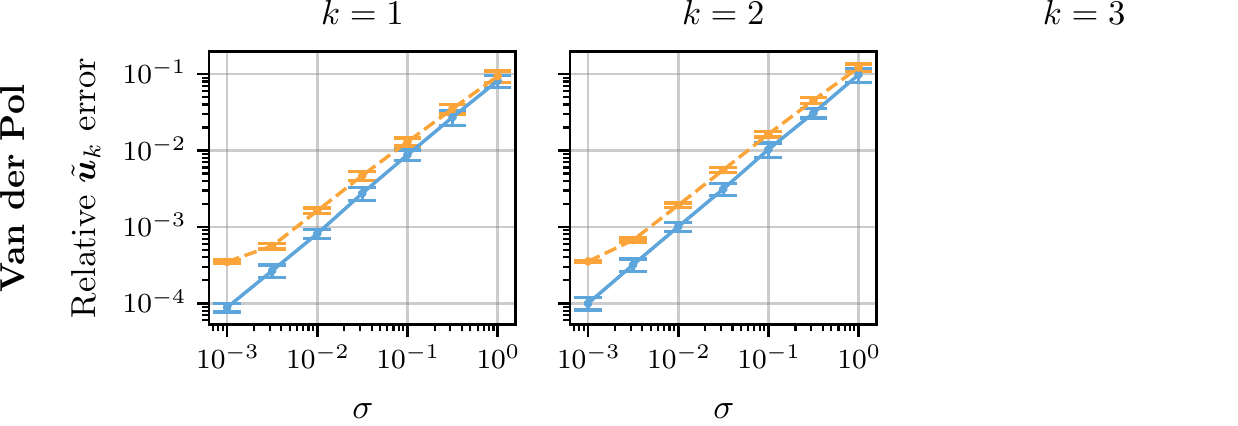}
    \includegraphics[trim=0 0 0 10, clip]{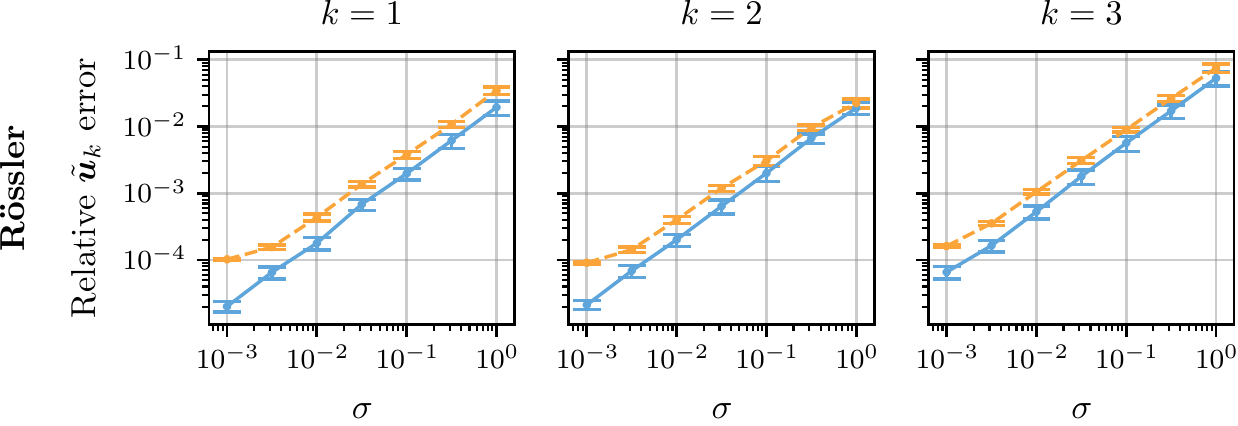}
    \includegraphics{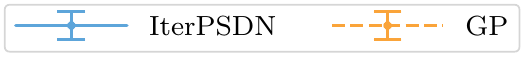}

    \caption{IterPSDN consistently outperforms or is equivalent to GP regression. Showing the mean $\pm$ std of 30 realizations with $N=1000$ at each noise level. For system parameters see \Cref{tab}.}
    \label{fig:PSDN-vs-GP}
\end{figure}

For all systems, IterPSDN (\Cref{alg:PSDN}) outperforms or is equivalent to GP regression (\Cref{fig:PSDN-vs-GP}).
IterPSDN and GP regression performed equally well at denoising the state measurement in the Duffing oscillator (PS1) (top row of \Cref{fig:PSDN-vs-GP}) as well as the Duffing oscillator (PS2) and the Lorenz 96 model (results not shown).
For the Van der Pol oscillator and R\"{o}ssler attractor, IterPSDN improves the state estimation (middle and bottom rows of \Cref{fig:PSDN-vs-GP}). In addition to improving state estimates, IterPSDN is faster than GP regression (given the standard deviation of the noise can be estimated quickly), and using projection state denoising allows us to estimate the relevant hyperparameters in the IRW-SOCP formulation (see \Cref{sec:methods}).

\subsubsection{The IRW-SOCP step of DSINDy improves state time derivative estimation}
\label{sec:results-deriv}

\begin{figure}
    \centering
    \includegraphics[trim=0 25 0 0, clip]{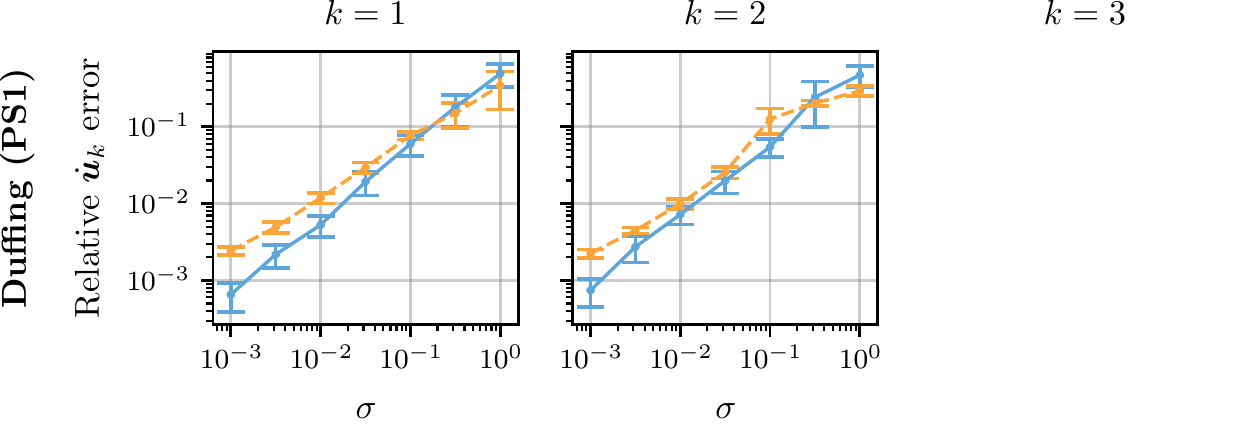}
    \includegraphics[trim=0 25 0 10, clip]{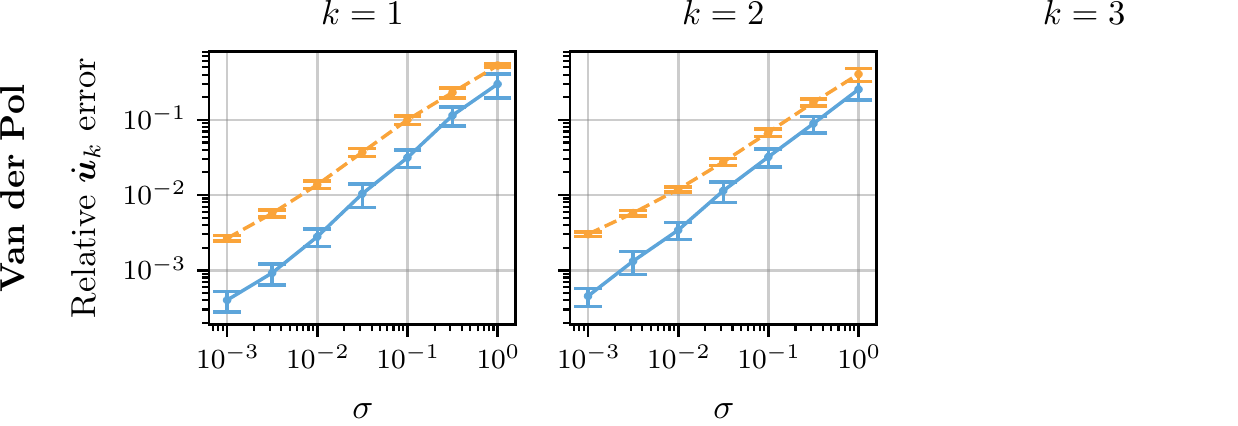}
    \includegraphics[trim=0 0 0 10, clip]{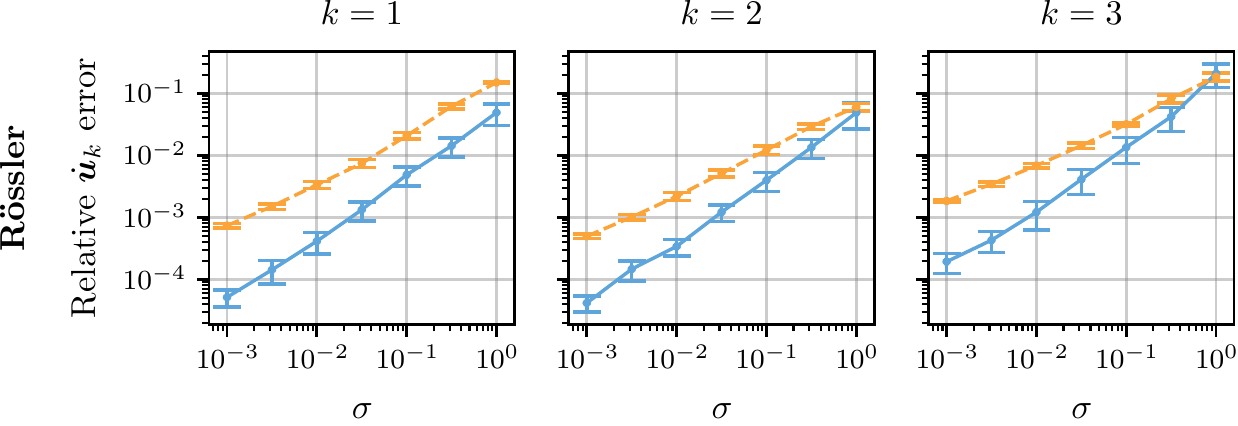}
    \includegraphics{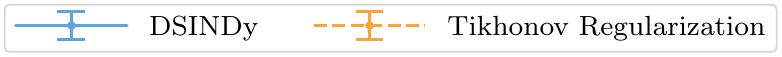}

    \caption{DSINDy improves derivative estimation compared with Tikhonov regularization. Showing the mean $\pm$ std of 30 realizations with $N=1000$ at each noise level. For system parameters see \Cref{tab}.}
    \label{fig:deriv}
\end{figure}

The IRW-SOCP step of DSINDy, see \Cref{eq:SOCP}, led to improved derivative estimation compared with Tikhonov regularization (see \Cref{fig:deriv}, results for the Lorenz 96 model are given in \Cref{sec-app:Lorenz96}). Although we observe improvement across multiple noise levels, this improvement becomes less significant as the standard deviation of the noise increases. Indeed, in some cases Tikhonov regularization slightly outperforms DSINDy at the highest noise level considered, e.g., see second state of Duffing oscillator (PS1) (\Cref{fig:deriv}, top row). However, as described below,  increasing the sample size may alleviate this issue.

\begin{figure}
    \centering
    \includegraphics[trim=0 32 0 0, clip]{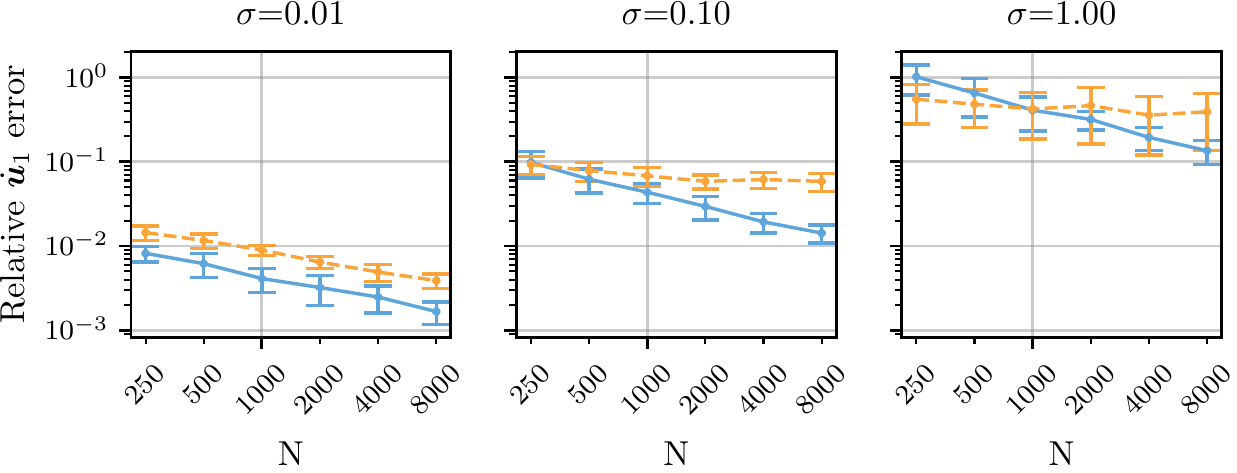}
    \includegraphics[trim=0 0 0 10, clip]{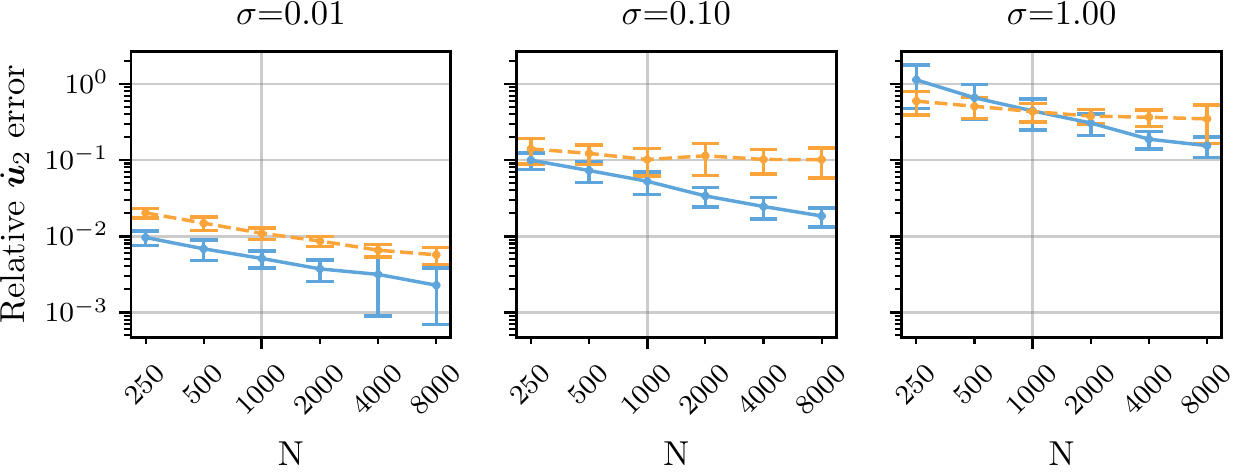}
    \includegraphics{figures/legends/legend_socp_sm_tikreg_h}
    \caption{Increasing the number of training samples improves DSINDy derivative estimation for the Duffing oscillator system (PS2). Each column corresponds to a different noise level. Showing the mean $\pm$ std of 30 noise realizations at each sample size. For system parameters see \Cref{tab}.}
    \label{fig:deriv-vs-N}
\end{figure}

We examined whether increasing the sample size improved the prediction of the state time derivative. Specifically, we ran DSINDy using $N=250,500,1000,2000,4000,8000$ measurements from the Duffing oscillator (PS2) using the same parameters given in \Cref{tab}. In DSINDy, the error of the derivative estimation consistently decreases with the square root of $N$ (\Cref{fig:deriv-vs-N}). However, for Tikhonov regularization this relationship is not always evident, e.g., see middle and right column of \Cref{fig:deriv-vs-N}. Additionally, we find that as the sample size increases, DSINDy starts to outperform Tikhonov regularization at $\sigma=1.0$ (see right column of \Cref{fig:deriv-vs-N}).

\subsubsection{DSINDy improves coefficient recovery and system reconstruction}
\label{sec:results-system}
\begin{figure}
    \hspace{1.8cm}
    \includegraphics[trim=0 25 0 0, clip, align=c]{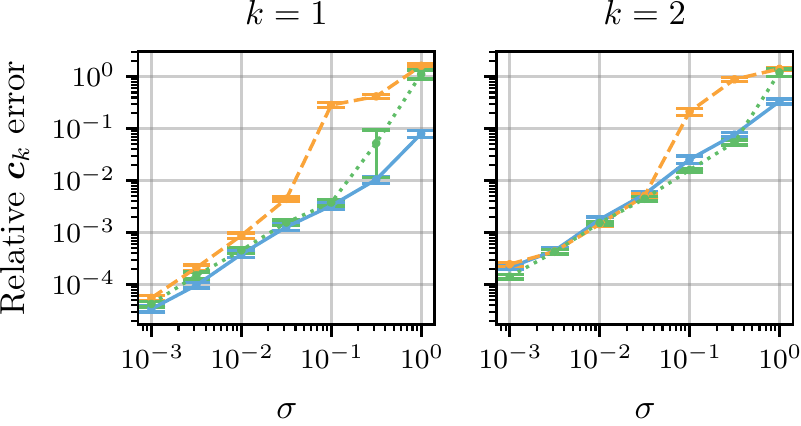}
    \hspace{.7cm}
    \includegraphics[align=c]{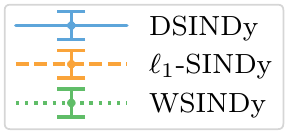}

    \hspace{1.8cm}
    \includegraphics[trim=0 0 0 10, clip]{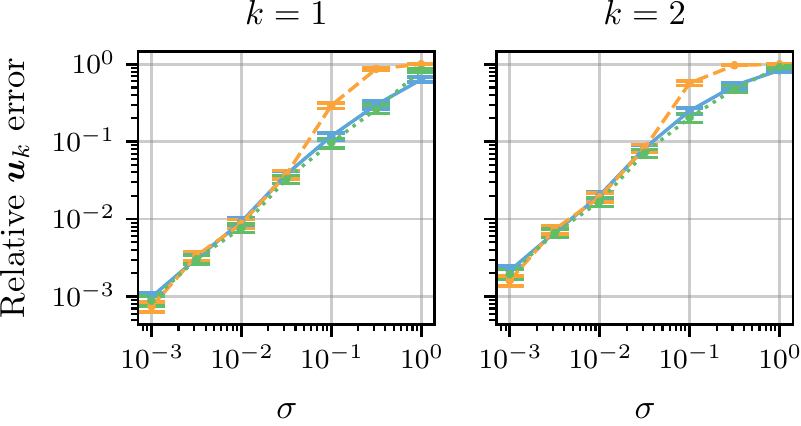}
    \hspace{.2cm}
    \includegraphics{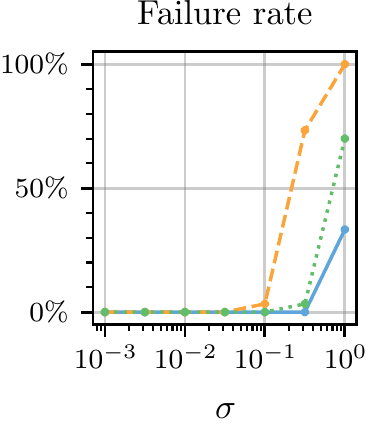}

    \caption{Coefficient and reconstruction error for Van der Pol system.  Showing the mean $\pm$ sem of 30 realizations with $N=1000$ at each noise level. The failure rate represents the fraction of simulations out of 30 that failed. For system parameters see \Cref{tab}.}
    \label{fig:DSINDy2}
\end{figure}

\begin{figure}
    \centering
    \includegraphics[trim=0 24 0 0, clip]{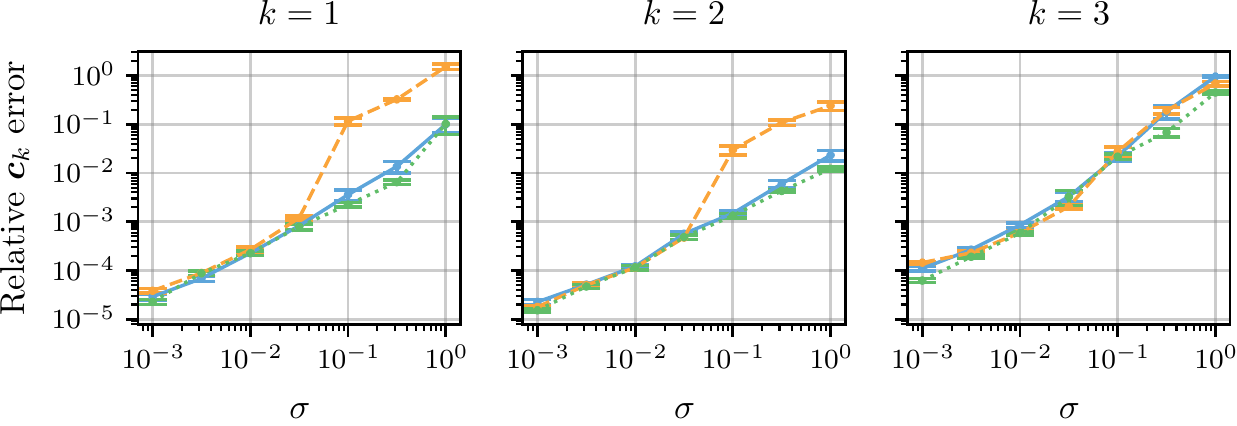}
    \includegraphics[trim=0 0 0 10, clip]{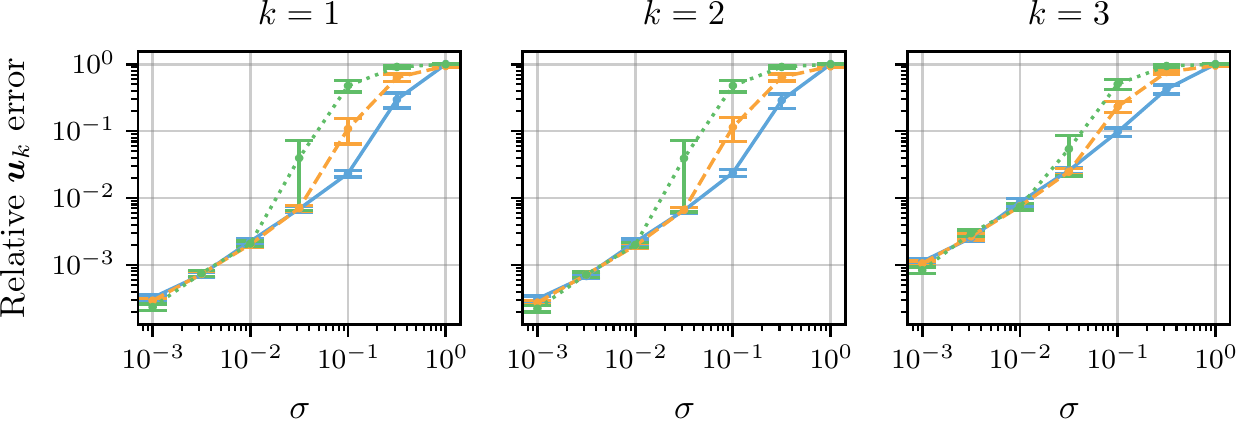}

    \includegraphics[align=c]{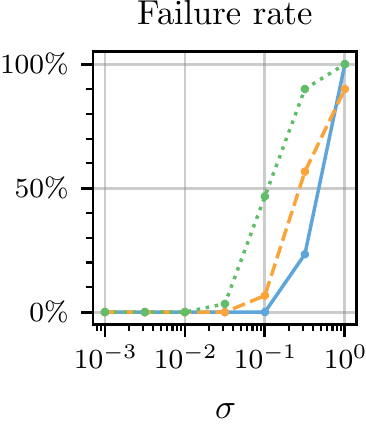}
    \hspace{.5cm}
    \includegraphics[align=c]{figures/legends/legend_socp_sm_lasso_WSINDY}

    \caption{Coefficient and reconstruction error for the R\"{o}ssler attractor. Showing the mean $\pm$ sem of 30 realizations with $N=1000$ at each noise level. The failure rate represents the fraction of simulations out of 30 that failed. For system parameters see \Cref{tab}.}
    \label{fig:DSINDy3}
\end{figure}

\begin{figure}
    \centering
    \includegraphics[trim=0 24 0 0, clip]{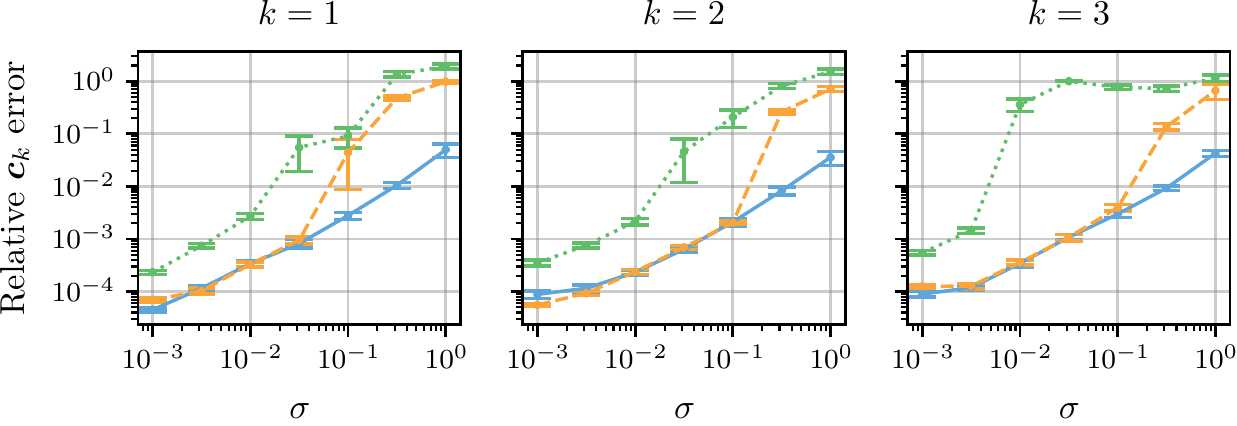}

    \vspace{.2cm}

    \includegraphics[trim=0 0 0 0, clip]{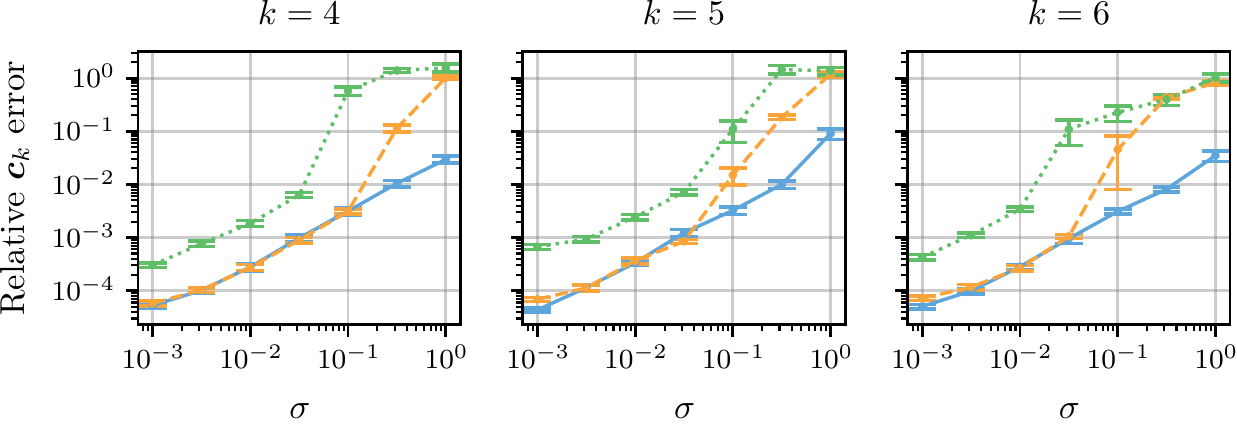}

    \includegraphics[align=c, trim=0 0 0 10, clip]{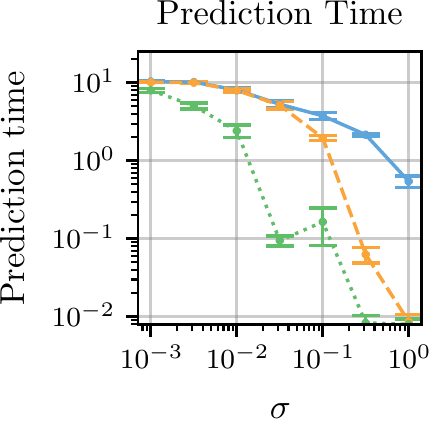}
    \hspace{.5cm}
    \includegraphics[align=c]{figures/legends/legend_socp_sm_lasso_WSINDY}

    \caption{Coefficient errors for the 6 states of the Lorenz 96 system. Showing the mean $\pm$ sem of 30 realizations with $N=2000$ at each noise level. The bottom plot shows the time at which the relative reconstruction error of any state    exceeds 10\% when the learned dynamics are initialized at $t=t_{end}=5$. Note that a prediction time less than $10^{-2}$ implies zero predictive ability as this is time step of the simulation. For system parameters see \Cref{tab}. An example simulation result with the prediction time labeled is given in \Cref{fig:Lorenz96-pred}.}
    \label{fig:DSINDy4}
\end{figure}

\begin{figure}
    \includegraphics[]{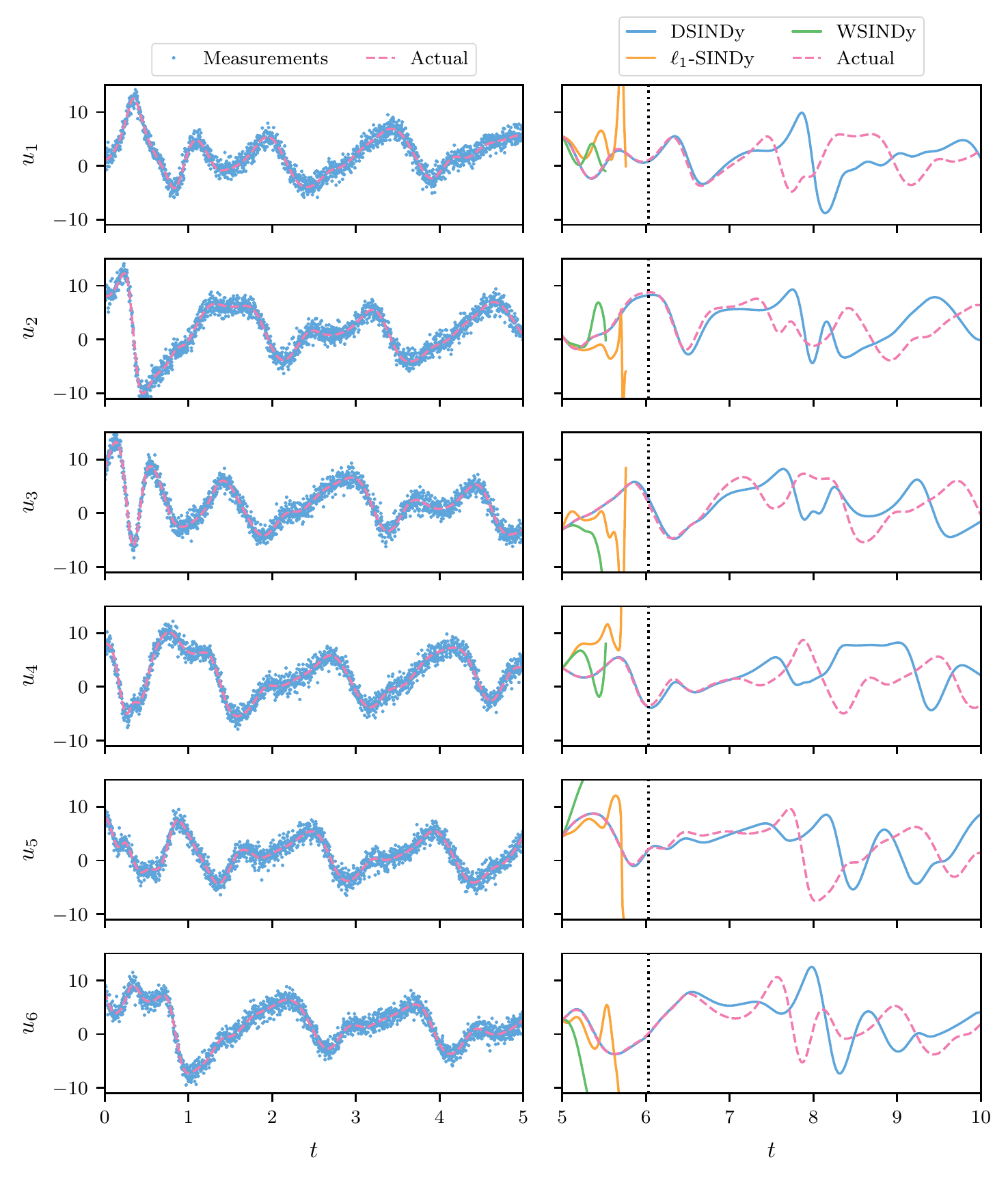}
    \caption{DSINDy improves predictions of the Lorenz 96 model compared with $\ell_1$-SINDy and WSINDy. Left column: Example set of measurements obtained for $t\in[0,5]$ ($\sigma = 1$, $N=2000$). Right column: Predicted system dynamics using example measurements starting at $t_{IC,test}=5$. The vertical black dotted line represents the time at which the relative DSINDy prediction error for one of the states, i.e., the third state, exceeds 10\%.}
    \label{fig:Lorenz96-pred}
\end{figure}

\begin{figure}
    \centering
    \includegraphics[trim=0 32 0 0, clip]{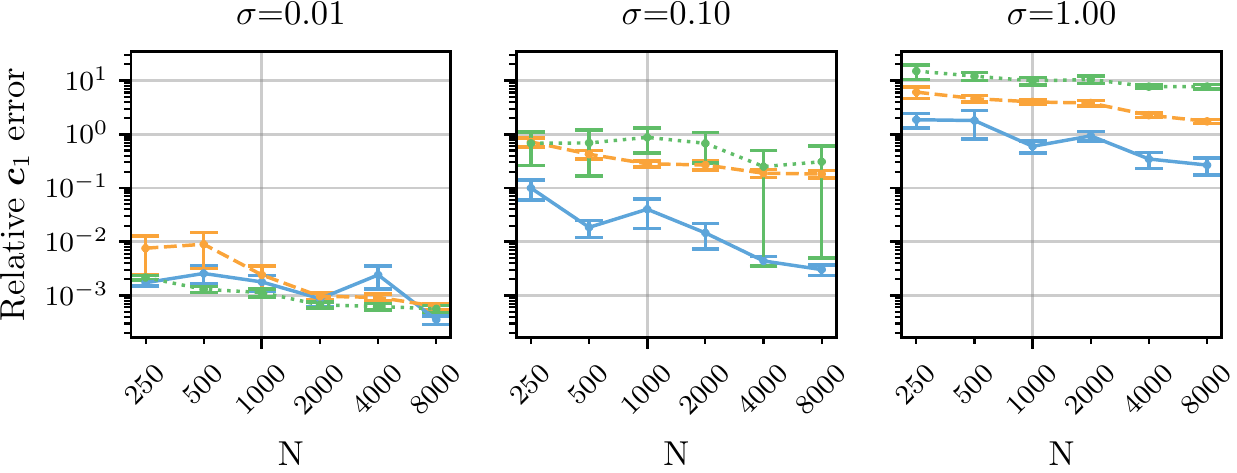}
    \includegraphics[trim=0 32 0 10, clip]{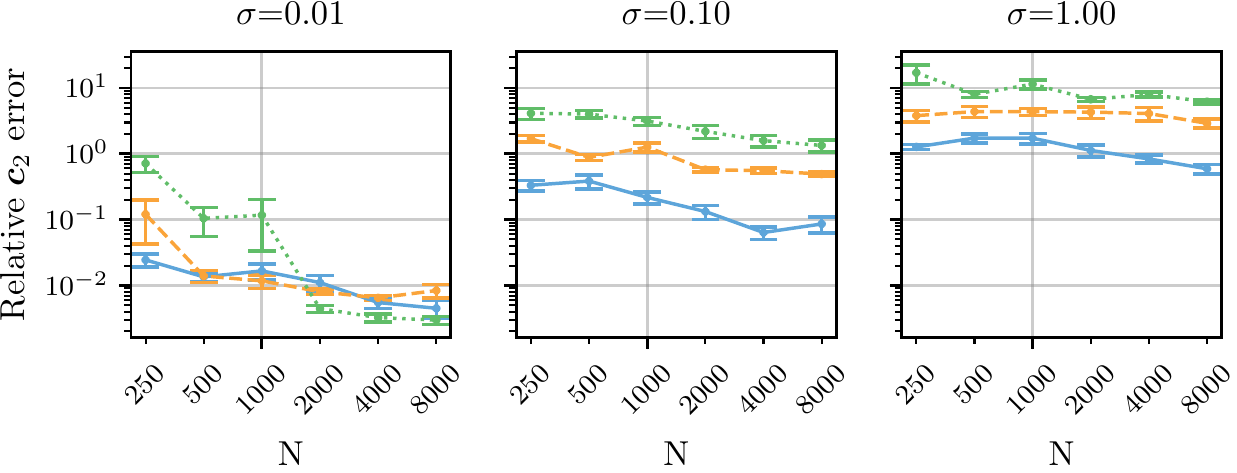}
    \includegraphics[trim=0 32 0 10, clip]{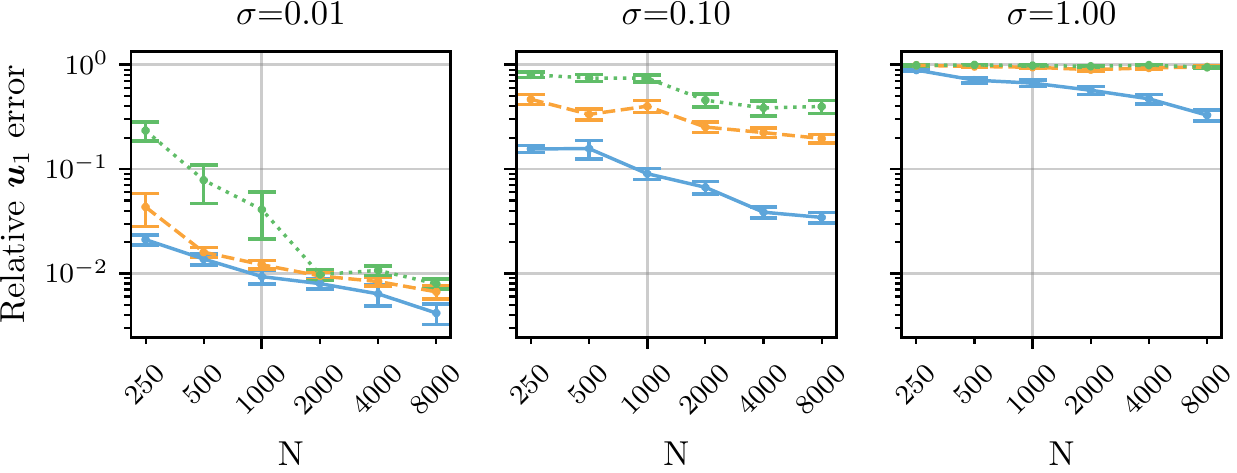}
    \includegraphics[trim=0 32 0 10, clip]{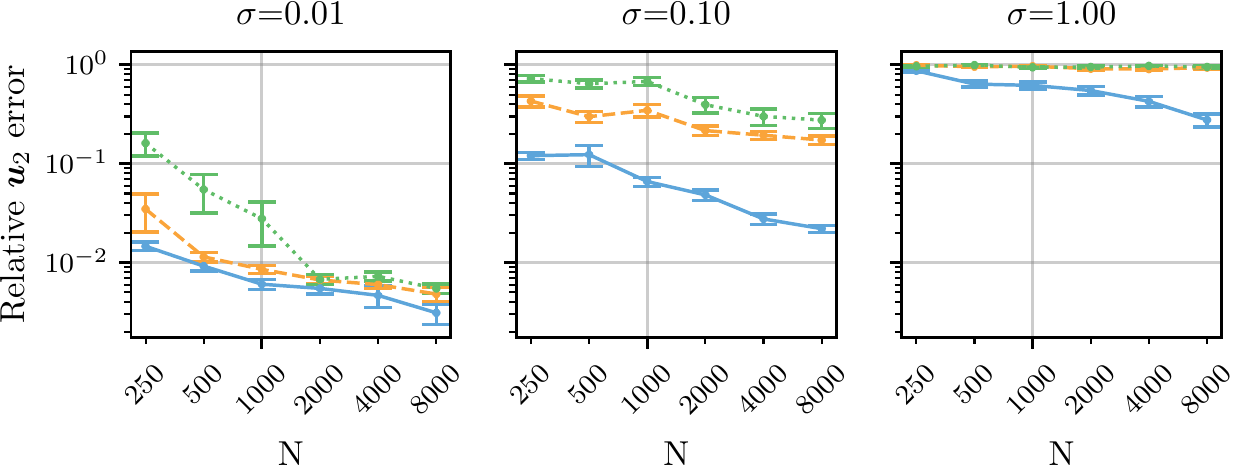}
    \includegraphics[trim=0 0 0 10, clip]{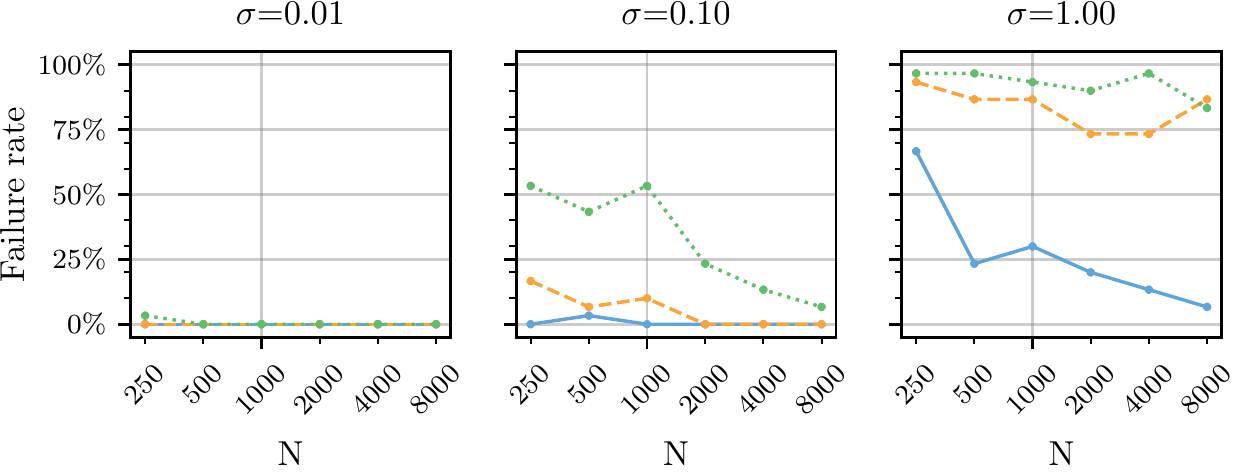}
    \includegraphics{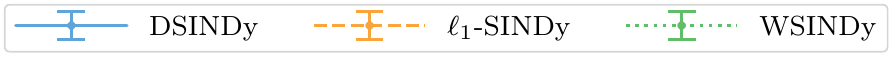}
    \caption{Coefficient and reconstruction error for Duffing oscillator system (PS2). Each column represents a different noise level. Showing the mean $\pm$ sem of 30 realizations at each noise level. The failure rate represents the fraction of simulations out of 30 that failed. For system parameters see \Cref{tab}. An example simulation result for this system is given in \Cref{fig:DSINDy-example}.}
    \label{fig:DSINDy-vs-N}
\end{figure}

DSINDy reduces the error of coefficient recovery and/or system reconstruction compared with $\ell_1$-SINDy and WSINDy (\Cref{fig:DSINDy2,fig:DSINDy3,fig:DSINDy4,fig:Lorenz96-pred,fig:DSINDy-vs-N}). Note that results for Duffing Oscillator (PS1) are given in \Cref{sec-app:Duffing}. In terms of coefficient recovery, DSINDy generally outperforms $\ell_1$-SINDy and WSINDy, e.g., see top row of \Cref{fig:DSINDy2,fig:DSINDy3} and top two rows of \Cref{fig:DSINDy4}. WSINDy only obtains better coefficients, compared with DSINDy, for the R\"{o}ssler attractor system (\Cref{fig:DSINDy3}, top row). However, this does not translate to a better system recovery (\Cref{fig:DSINDy3}, bottom row). This is likely because WSINDy fails to identify one of the small nonzero coefficients. Indeed, at $\sigma=0.1$, in 14 out of 30 replications, WSINDy fails to identify the coefficient that corresponds with $\beta$ in the R\"{o}ssler attractor system (see \Cref{tab}).

\begin{figure}
    \includegraphics{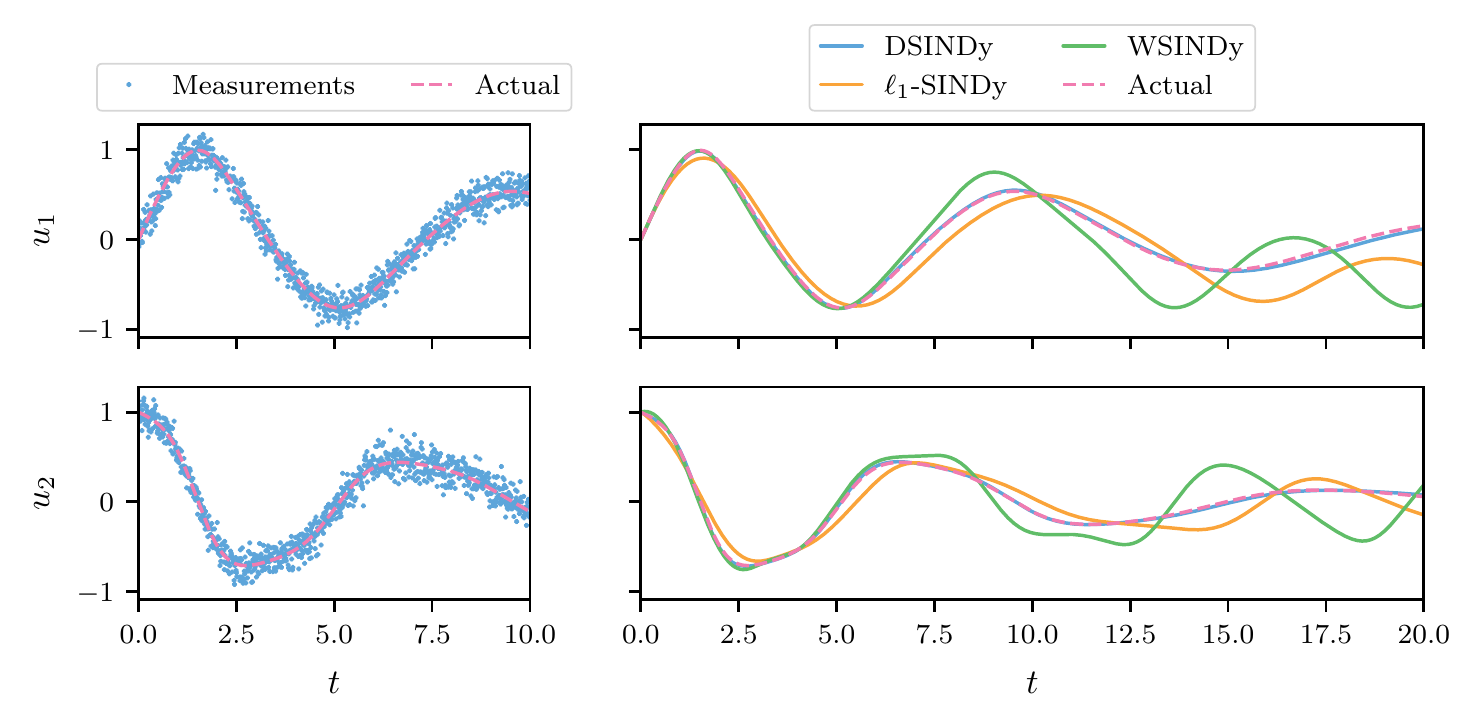}
    \caption{Example results for the Duffing oscillator (PS2) system. Left column: Example set of measurements obtain for $t \in [0,10], \sigma=0.1, N=1000$. Right column: Predicted system dynamics using example measurements over double the training time.}
    \label{fig:DSINDy-example}
\end{figure}

In terms of system reconstruction, DSINDy consistently outperforms both $\ell_1$-SINDy and WSINDy, e.g., see \Cref{fig:DSINDy3,fig:DSINDy-example,fig:DSINDy-vs-N}. The one exception to this is the Van der Pol oscillator where the reconstruction error of DSINDy is equivalent to WSINDy (\Cref{fig:DSINDy2}, bottom row). DSINDy led to accurate predictions of short term dynamics in the chaotic Lorenz 96 model, even when the measurements were highly corrupted by noise (see \Cref{fig:DSINDy4}, bottom panel). \Cref{fig:Lorenz96-pred} shows an example of this predictive capability, alongside state measurements obtained at $\sigma = 1$. In contrast to DSINDy, neither $\ell_1$-SINDy nor WSINDy were able to predict short term dynamics of the Lorenz 96 model at the highest noise level considered.

The benefits of DSINDy are also apparent for the Duffing oscillator (PS2) system, where we consider the effect of sample size (\Cref{fig:DSINDy-vs-N}). At low noise levels, DSINDy and $\ell_1$-SINDy performed similarly (\Cref{fig:DSINDy-vs-N}, first column), whereas WSINDy struggled at the low sample sizes, i.e., at $N=250, 500$. At the two highest noise levels, $\sigma=0.1$ and $\sigma=1.0$, DSINDy improves upon both $\ell_1$-SINDy and WSINDy in terms of both the coefficient and reconstruction errors (\Cref{fig:DSINDy-vs-N}, second and third columns). In particular, at the highest noise level considered, increasing the sample size improves the performance of DSINDy but does not improve the performance of the other two methods.

\subsubsection{Method for selecting data-matching hyperparameter \texorpdfstring{$\gamma$}{}}
\label{sec:results-dm-hyper}
We did not observe a significant affect of the data-matching hyperparameter selection method on the performance of DSINDy in any of the ODE systems (see \Cref{sec-app:results} for complete results). Recall, that we either found $\gamma$ in \Cref{eq:SOCP} using the corner point of the Pareto curve (as done in the previous section) or using the estimate given in \Cref{eq:gamma_exp}. Although both methods rely on our theoretical results, the second approach is much faster as we do not have to explore the hyperparameter space. These results suggest that directly using the theoretical estimate may improve algorithmic speed without compromising performance.

\section{Discussion}\label{sec:discussion}

We developed DSINDy to improve ODE discovery when state measurements are highly corrupted by noise. The method leverages the assumed (overcomplete) basis for the dynamics to denoise the state variable measurements, and finds the state variable time derivatives while enforcing sparsity of the coefficients (\Cref{sec:methods}). Our goal was to improve derivative estimation in the presence of noise, and, in turn improve overall system recovery. We found that DSINDy successfully accomplished this goal for the Duffing oscillator, the Van der Pol oscillator, the R\"{o}ssler attractor, and the higher dimensional Lorenz 96 model (\Cref{sec:results}).

The performance of DSINDy improves as the sample size increases. Theoretical results for the \textit{a priori} denoising strategy PSDN show that the denoised state variables asymptotically approach their true values. Although theoretical results are not available for iterative projections (IterPSDN, \Cref{alg:PSDN}), we numerically found the error of IterPSDN also decreased with the number of samples (see \Cref{fig:PSDN-Theory1,fig:PSDN-Theory2,fig:PSDN-Theory3}).  Additionally, for the Duffing oscillator we see improvements in system recovery as the sample size increases (see \Cref{fig:deriv-vs-N,fig:DSINDy-vs-N}). In particular, at high noise levels, this improvement is exclusive to DSINDy.

There are multiple areas of future research that may improve the DSINDy approach. For example, currently we set the value of the smoothing hyperparameter in IRW-SOCP to guarantee feasibility of \Cref{eq:SOCP}. In practice this value is often larger than necessary, and therefore, future work may involve exploring alternative ways to learn this hyperparameter. This may also provide insight into whether the Pareto curve or theoretical results should be used to set the value of the data-matching hyperparameter, as current results are inconclusive (\Cref{sec:results-dm-hyper}). For the smoothing step, theoretical guarantees do not exist for IterPSDN (\Cref{alg:PSDN}). Therefore, a more thorough analysis and future modifications to this algorithm are warranted.

On top of the SINDy framework additional approaches have been developed to help alleviate the challenge of noise. These include approaches that use Bayesian inference \cite{Niven2020,Galioto2020,Hirsh2022} and ensemble-based/bootstrapping techniques \cite{Fasel2022, Abdullah2022} to obtain estimates of coefficient uncertainty. In \cite{Delahunt2022} a toolkit was proposed that provides extensions to the SINDy framework to help deal with high noise measurements. For example, the toolkit involves modifying the basis after comparing the obtained dynamics to the original measurements. Similar to SINDy, DSINDy could be integrated into these larger frameworks.

One of the biggest challenges with equation discovery is finding a suitable basis for the dynamics. The PSDN approach may potentially be modified to help with basis selection by observing how the inclusion of different basis elements impacts the results. Currently, the theory for PSDN does not hold if the library does not have full column rank, i.e., when there are redundant terms in the basis. This can be a challenging scenario to detect at high noise levels \cite{Cortiella2021}. A topic of future work is to examine the performance of PSDN in this scenario and determine if identification of redundant basis elements is possible in the high-noise regime.

\section*{Acknowledgments}
The authors would like to thank Professor David Bortz and Daniel Messenger for helpful discussions regarding this work. This work was supported by the Department of Energy, National Nuclear Security Administration, Predictive Science Academic Alliance Program (PSAAP) [Award Number DE-NA0003962].

\appendix
% This fixes appendix labels when using cleveref
\gdef\thesection{\Alph{section}} % corrected redefinition of "\thesection"
\makeatletter
\renewcommand\@seccntformat[1]{\appendixname\ \csname the#1\endcsname.\hspace{0.5em}}
\makeatother

\setcounter{theorem}{0}
\setcounter{lemma}{0}
\setcounter{corollary}{0}
\setcounter{claim}{0}
\setcounter{remark}{0}
\renewcommand{\thetheorem}{\Alph{section}\arabic{theorem}}
\renewcommand{\thelemma}{\Alph{section}\arabic{lemma}}
\renewcommand{\thecorollary}{\Alph{section}\arabic{corollary}}
\renewcommand{\theclaim}{\Alph{section}\arabic{claim}}
\renewcommand{\theremark}{\Alph{section}\arabic{remark}}

\section{Description of equation discovery methods that were compared with DSINDy}\label{sec-app:algorithms}

Here we present existing methods that are used for comparison purposes in the main manuscript. These include sparsity-promoting regularization techniques (\Cref{sec-app:bg-reg}) and the Pareto corner criteria for finding the hyperparameter in regularized least squares problems (\Cref{sec-app:Pareto}). Note that the Pareto corner criteria is also used within the DSINDy algorithm. We then provide a description $\ell_1$-SINDy and WSINDy (\Cref{sec:bg-find-dyn1,sec:bg-find-dyn2}). Both equation discovery methods involve the use of sparsity-promoting regularization techniques discussed in \Cref{sec-app:bg-reg}, and $\ell_1$-SINDy uses the Pareto corner criteria discussed in \Cref{sec-app:Pareto}.

\subsection{Sparsity-promoting regularization}
\label{sec-app:bg-reg}

Here we discuss two sparsity-promoting regularization techniques: sequential-thresholding least squares (STLS) \cite{Brunton2016} and iteratively reweighted (IRW) Lasso \cite{Zou2006,Candes2008,Cortiella2021}\footnote{Variations of IRW-Lasso are referred to as adaptive lasso in \cite{Zou2006} and WBPDN in \cite{Cortiella2021}.}. Given a measurement matrix $A \in \mathbb{R}^{N \times p}$ and vector $\bm{b} \in \mathbb{R}^N$, the goal of both methods is to reconstruct a sparse signal (or coefficient vector) $\bm{c} \in \mathbb{R}^p$ such that $A \bm{c} \approx \bm{b}$. In the paper we only considered overdetermined systems, i.e., $N \ge p$.

STLS is a sparsity-promoting method that iteratively sets coefficients to zero \cite{Brunton2016}. Briefly, STLS involves finding the least squares solution, $\bm{c}^{LS} = A^{\dagger} \bm{b}$. Then, all coefficients with magnitudes less than a threshold, $\lambda_{STLS}$, are set equal to zero and the corresponding columns of $A$ are removed. These steps are repeated iteratively until the coefficient vector converges. For the complete algorithm see \cite{Brunton2016}.

IRW-Lasso iteratively solves the Lasso $\ell_1$-regularization problem \cite{Tibshirani1996}, i.e.,
\begin{equation}\label{eq:Lasso}
    \minimize_{\bm{c}} \quad \|A \bm{c} - \bm{b}\|^2 + \lambda_{Lasso} \|W \bm{c}\|_1.
\end{equation}
In standard Lasso as well as in the first iteration of IRW-Lasso, $W$ is the identity matrix. In subsequent iterations of IRW-Lasso, $W$ is a diagonal matrix such that
\begin{equation}
    W_{ii} = \frac{1}{|c^{prior}_i| + \epsilon}, \quad \quad i=1,2,...,p,
\end{equation}
where $\bm{c}^{prior}$ is the coefficient vector obtained in the previous iteration and $\epsilon$ is a small positive number to avoid dividing by zero. For more details on IRW-Lasso and its implementation see \cite{Zou2006,Cortiella2021}. In DSINDy we use a similar iterative process as IRW-Lasso to find a sparse vector $\bm{c}$. However, the formulation differs from \Cref{eq:Lasso} (see \Cref{sec:methods-socp} for details).

\begin{remark}
    Many additional strategies exist for finding a sparse vector $\bm{c}$. These include Orthogonal Matching Pursuit (an $\ell_0$ minimization approach, \cite{Pati1993}) and Basis Pursuit Denoising (a version of \Cref{eq:Lasso} where the $\| A\bm{c}-\bm{b}\|$ loss is incorporated as a constraint, \cite{Chen2001}). However, these approaches are not considered in the present work.
\end{remark}

\subsection{Pareto corner finding algorithm}\label{sec-app:Pareto}
Let $A \in \mathbb{R}^{N \times p}$ and $\bm{b} \in \mathbb{R}^N$, represent a measurement matrix and vector, respectively. Suppose we want to find a coefficient vector $\bm{c} \in \mathbb{R}^p$, such that $A \bm{c} \approx \bm{b}$, by solving
\begin{equation}\label{eq-app:Lasso}
    \minimize_{\bm{c}} \quad \|A \bm{c} - \bm{b}\|^2 + \lambda \|\bm{c}\|_s,
\end{equation}
where $s$ could be 1 or 2 depending on the type of regularization desired. Here, we describe a method that uses the Pareto curve to find the value of the hyperparameter $\lambda$. The Pareto curve displays the trade-off between minimizing the solution residual, $\|A\bm{c}-\bm{b}\|$, and the regularization residual $\|\bm{c}\|_s$ \cite{Hansen1992,VanEwoutBerg2008,Cortiella2021}.

To generate the Pareto curve, we solve \Cref{eq-app:Lasso} at multiple values of the hyperparameter, i.e. set $\lambda=\lambda_1,\lambda_2,...,\lambda_n$, where $\lambda_i \in [0,\infty)$, to obtain $\bm{c}(\lambda_i)$ for $i=1,2,...,n$. We then plot the resulting solution residual verse the regularization residual, i.e., points contained in the set
\begin{equation}
    \{\left(\|\bm{c}(\lambda)\|_s, \|A\bm{c}(\lambda)-\bm{b}\|\right) \mid \lambda \in [0,\infty)\}.
\end{equation}
When plotted on a log scale, the resulting curve often has an L-shape characteristic. The value of $\lambda$ at the corner of this L-curve, i.e., $\lambda_c$, is thought to represent a balance between the regularization and solution error \cite{Hansen1992}. Therefore, we choose $\bm{c}(\lambda_c)$ as the final coefficient vector.

Various methods have been proposed for finding $\lambda_c$ efficiently \cite{Singh2011,Cultrera2020,Cuate2020}. One algorithm iteratively performs a golden section search to narrow in on a point of the Pareto curve with large curvature \cite{Cultrera2020}. This is the algorithm that is used in the present work.

\subsection{\texorpdfstring{$\ell_1$}{l1}-SINDy}\label{sec:bg-find-dyn1}

The $\ell_1$-SINDy algorithm is a variation of SINDy \cite{Brunton2016} that uses $\ell_1$-minimization to find the coefficient vector. Specifically, the $\ell_1$-SINDy algorithm incorporates the modifications to SINDy presented in \cite{Cortiella2021}. In $\ell_1$-SINDy we solve \Cref{eq:Thetac_equal_udot} directly by approximating the state time derivatives and assuming that the coefficient vector is sparse. The recovery of the ODE system then reduces to a regularized linear regression problem. The original SINDy paper \cite{Brunton2016} used STLS to enforce sparsity of the coefficient vector. Later, \cite{Cortiella2021} found that IRW-Lasso led to improved coefficient prediction. Therefore, in $\ell_1$-SINDy we also use IRW-Lasso (see \Cref{sec-app:bg-reg}). Additionally, we will follow \cite{Cortiella2021} and estimate the sparsity-hyperparameter $\lambda_{lasso}$ using the corner point of the Pareto curve (see \Cref{sec-app:Pareto}). We refer to this algorithm as $\ell_1$-SINDy to emphasize this departure from the original SINDy framework.

The SINDy algorithm, and hence $\ell_1$-SINDy, requires that the derivative be estimated from the noisy measurements, $\bm{u} \in \mathbb{R}^N$. Again, we follow the approach of \cite{Cortiella2021} and use Tikhonov regularization to find the derivative. We  solve
\begin{equation}\label{eq:Tik}
    \minimize_{\dot{\bm{u}}} \quad \|T \dot{\bm{u}} - \bm{u}\|^2 + \lambda_{Tik} \|D \dot{\bm{u}}\|^2,
\end{equation}
where $T$ and $D$ are as defined in \Cref{sec:notation}.
The first term of \Cref{eq:Tik} requires that the integrated time derivative predict the measurements and the second term promotes smoothness of the derivative. When implementing this method, we used the Pareto corner finding algorithm to find the hyperparameter $\lambda_{Tik}$ (see \Cref{sec-app:Pareto}) \cite{Cultrera2020}.

One modification we use here in $\ell_1$-SINDy that was not used in \cite{Cortiella2021} is that we perform \textit{a priori} denoising of the state measurements. We use the smoothed measurements from IterPSDN as input to $\ell_1$-SINDy. We do this to more accurately compare $\ell_1$-SINDy with the IRW-SOCP step of DSINDy.

\subsection{WSINDy}\label{sec:bg-find-dyn2}

WSINDy uses the weak formulation to learn the governing equations \cite{Messenger2021}. In WSINDy, we construct a set of test functions that are compactly supported within the time domain of the training data, i.e., $\phi_\ell(t) \in C_c([0,t_{end}])$ for $\ell=1,...,L$. The inner product, i.e., $\langle u,v \rangle = \int_0^{t_{end}} u(t)v(t)dt$, of the dynamics for state $k$ with each test function is given as, i.e.,
\begin{equation}
    \langle \phi_{\ell}, \dot{u}_k \rangle = \langle \phi_\ell, F_k(u_1,u_2,...,u_m) \rangle,
\end{equation}
and modified using integration by parts to obtain,
\begin{equation}\label{eq:wsindy-ibp}
    - \langle \dot{\phi}_\ell, u_k \rangle = \langle \phi_\ell, F_k(u_1,u_2,...,u_m) \rangle.
\end{equation}
Note that since the test functions are compactly supported on $[0,t_{end}]$, the boundary terms go to zero. A discrete approximation to the integral in \Cref{eq:wsindy-ibp}, based on the measurement acquisition times, is given as
\begin{equation}\label{eq:WSINDy-eq}
    - \sum_{i=1}^N \dot{\phi}_\ell(t_i) u_k(t_i) \approx \sum_{i=1}^N \phi_\ell(t_i) \Theta_{i} \bm{c}_k,
\end{equation}
where $\Theta_i$ represents the $i$th row of the basis function library $\Theta$, defined in \Cref{eq:Theta}. Note that \Cref{eq:WSINDy-eq} is not exact as the state variables are approximated from noisy measurements and error is introduced in discretization.

Finally, we consider \Cref{eq:WSINDy-eq} for all test functions to obtain a linear system,
\begin{equation}\label{eq:WSINDy-eq2}
    \bm{b}_k \approx H \bm{c}_k,
\end{equation}
where
$b_{k,\ell} = - \sum_{i=1}^N \dot{\phi}_\ell(t_i) u_k(t_i)$ and the $\ell$th row of $H$ is $ H_{\ell} = \sum_{i=1}^N \phi_\ell(t_i) \Theta_{i}$.  One key benefit of WSINDy is that it avoids the state derivative calculation through the use of integration by parts.

To find the sparse coefficient vector $\bm{c}_k \in \mathbb{R}^p$ for $k=1,2,...,m$, the most recent WSINDy code available uses a modified version of STLS, i.e., MSTLS. This version differs from STLS in that an upper bound on the coefficient magnitude is incorporated. To find the optimal hyperparameter $\lambda$,
MSTLS is performed at a range of $\lambda$ values, and the chosen $\lambda$ minimizes the following loss function,
\begin{equation}
    \mathcal{L}(\lambda) = \frac{\|H(\bm{c}(\lambda) -  \bm{c}^{LS})\|}
    {\|H\bm{c}^{LS}\|}
    + \frac{\|\bm{c}(\lambda)\|_0}{p},
\end{equation}
where $\bm{c}^{LS}$ is the least squares solution.
The first term of this loss function encourages the coefficient vector to match the least squares solution, while the second term encourages sparsity (for additional details see \cite{Messenger2021a}). Note that MSTLS was not used in the original WSINDy paper for discovering ODE systems \cite{Messenger2021}.

The current implementation of WSINDy on GitHub also performs prior smoothing of the measurements. The method smooths the data by finding a weighted average of each point with its nearest neighbors. The weighting and size of the smoothing window is calculated automatically.

\setcounter{theorem}{0}
\setcounter{lemma}{0}
\setcounter{corollary}{0}
\setcounter{claim}{0}
\setcounter{remark}{0}
\section{Additional theoretical results and proofs}\label{sec-app:theory}

\setcounter{theorem}{0}
\setcounter{lemma}{0}
\setcounter{corollary}{0}
\setcounter{claim}{0}
\setcounter{remark}{0}

In this section we provide supplemental theoretical results. We first present an error bound on the $\ell_2$-error of discrete integration which follows almost immediately from the well known trapezoidal rule error bound (\Cref{sec-app:trap-err}). We then present existing matrix perturbation results that provide a bound on the projection operator error (\Cref{sec-app:P-err}), and we then provide a bound for the pseudoinverse of $\Psi^*$ (\Cref{sec:psi-pinv-norm-bound}). Next, we explicitly present the unbiased monomial library $\hat{\Theta}$ and derive an estimator of the matrix $(\Psi^*)^T\Psi^*$ (\Cref{sec-app:theory-estimator}). Finally, we show that the IRW-SOCP formulation given by \Cref{eq:SOCP} is feasible for our chosen value of the derivative-smoothing hyperparameter (\Cref{sec-app:smoothing-hyper-param}).

\subsection{Trapezoidal quadrature error}
\label{sec-app:trap-err}

The proof of \Cref{lemma:err_bound_1} requires a bound on the error introduced by discrete integration in \Cref{eq:quad-err}. This error bound, specific to the trapezoidal rule, is given in the following lemma.

\begin{lemma}\label{lemma:quad_error}
    Let $T \in \mathbb{R}^{N \times N}$ be the discrete integral operator that performs trapezoidal quadrature as defined in \Cref{eq:T}. Let $0=t_1,t_2,...,t_N=t_{end}$ represent a set of $N$ equispaced time points in the interval $[0,t_{end}]$ and let $u(t)$ be defined by the following integral equation
    \begin{equation}
        u(t) = u(0) + \int_0^{t} \dot{u}(t')dt',
    \end{equation}
    where $\dot{u} = du/dt$ and assume that the constant
    \begin{equation}
        M = \max_{t \in [0,t_{end}]} \left| \frac{d^3 u}{dt^3} \right|
    \end{equation}
    exists and is finite.

    Suppose we estimate $\bm{u}=[u(t_1),...,u(t_{end})]$ as
    \begin{equation}
        \hat{\bm{u}} = u(0) + T\dot{\bm{u}}
    \end{equation}
    where $u(0)$ and $\dot{\bm{u}}=[\dot{u}(t_1),...,\dot{u}(t_N)]$ are known.
    Define the quadrature error vector $\bm{e}_q := \hat{\bm{u}} - \bm{u}$. Then
    \begin{equation}
        \|\bm{e}_q\| \le \frac{t_{end}^3 M}{12} (N-1)^{-3/2}.
    \end{equation}
\end{lemma}

\begin{proof}
    Using the error bound for the trapezoidal rule, we have that
    \begin{equation}\label{eq:e_qi}
        \left|e_{q,i}\right| \le \frac{((i-1) \Delta t)^3 M}{12 (i-1)^2}
        = \frac{(i-1) (\Delta t)^3 M}{12} = \frac{(i-1) t_{end}^3 M}{12 (N-1)^3}.
    \end{equation}
    Using \Cref{eq:e_qi}, we next derive a bound on $\|\bm{e}_q\|^2$ noting that $e_{q,1}=0$,
    \begin{equation}\label{eq:quad_proof}
        \begin{aligned}
            \|\bm{e}_{q}\|^2 = \sum_{i=2}^N e_{q,i}^2 & \le \sum_{j=1}^{N-1} \left(\frac{j t_{end}^3 M}{12 (N-1)^3}\right)^2              \\
                                                      & = \left(\frac{t_{end}^3 M}{12 (N-1)^3}\right)^2 \left(\sum_{j=1}^{N-1} j^2\right) \\
                                                      & = \left(\frac{t_{end}^3 M}{12 (N-1)^3}\right)^2 \frac{(N-1)N(2N-1)}{6}            \\
                                                      & = \left(\frac{t_{end}^3 M}{12}\right)^2 \frac{N(2N-1)}{6(N-1)^5}                  \\
                                                      & \le \left(\frac{t_{end}^3 M}{12}\right)^2 \frac{1}{(N-1)^3}.                      \\
        \end{aligned}
    \end{equation}
    The last inequality follows because for $N\ge 2$,
    \begin{equation}
        \frac{N}{N-1} \le 2 \quad \text{and} \quad \frac{2N-1}{N-1} \le 3.
    \end{equation}
    Taking the square root of \Cref{eq:quad_proof} completes the proof.
\end{proof}

\subsection{Existing results on projection operator error}
\label{sec-app:P-err}

In this section we present results from a study that examined how matrix perturbations impact projection operations \cite{Stewart1977}. These existing results are used in the proof of \Cref{lemma:DeltaP_bound} and are applicable to any two matrices $\Phi^*$ and $\hat{\Phi}$, where $\hat{\Phi}$ is a perturbed version of $\Phi^*$.  Note that we use these matrix names to place the results of \cite{Stewart1977} in the context of the equation discovery problem presented in \Cref{sec:statement}. We will use the following definitions,
\begin{equation}
    \Delta \Phi := \hat{\Phi} - \Phi^*, \quad\quad \Delta P := P_{\hat{\Phi}} - P_{\Phi^*},
\end{equation}
where the projection operators ($P_{\hat{\Phi}}$ and $P_{\Phi^*}$) are calculated using \Cref{eq:P_A}.

\begin{theorem}[modified version of Theorem 4.1 in \cite{Stewart1977}] \label{thm:stewart}
    Let the SVD of $\Phi^*$ be given as
    \begin{equation*}
        \Phi^* =
        \begin{bmatrix} \hat{U}_{\Phi^*} & \check{U}_{\Phi^*} \end{bmatrix}    \begin{bmatrix} \Sigma_{\Phi^*} & 0 \\ 0 & 0 \end{bmatrix}
        \begin{bmatrix} \hat{V}_{\Phi^*} & \check{V}_{\Phi^*} \end{bmatrix}^T.
    \end{equation*}
    Define the error matrices
    \begin{equation*}
        E_{11} := \hat{U}_{\Phi^*}^T \Delta\Phi \hat{V}_{\Phi^*}, \quad\quad
        E_{21} := \check{U}_{\Phi^*}^T \Delta\Phi \hat{V}_{\Phi^*},
    \end{equation*}
    and scalar value
    \begin{equation*}
        \beta := \frac{\|\Phi^*\|\|(\Phi^*)^\dagger\|}{1 - \|E_{11}\|\|(\Phi^*)^\dagger\|}.
    \end{equation*}
    If $\Rank(\hat{\Phi})=\Rank(\Phi^*)$ and $\|(\Phi^*)^\dagger\| \|E_{11}\| < 1$, then
    \begin{equation}\label{eq:bound_stewart}
        \|\Delta P\| = \|P_{\hat{\Phi}} - P_{\Phi^*}\| \le \frac{\beta \|E_{21}\|/\|\Phi^*\|}{(1+(\beta\|E_{21}\|/\|\Phi^*\|)^2)^{1/2}}.
    \end{equation}
\end{theorem}

The statement of \Cref{thm:stewart} differs from the statement of Theorem 4.1 in \cite{Stewart1977} in two ways. First, Theorem 4.1 in \cite{Stewart1977} requires that $\hat{\Phi}$ be an acute perturbation of $\Phi^*$. Here, we instead require  that $\Rank(\hat{\Phi})=\Rank(\Phi^*)$ and $\|(\Phi^*)^\dagger\| \|E_{11}\| < 1$. Together, these requirements imply that $\hat{\Phi}$ is an acute perturbation of $\Phi^*$ (see Theorem 2.2 and Theorem 2.5 in \cite{Stewart1977}). Second, we replace the variables $\kappa$ and $\gamma$, referenced in \cite{Stewart1977} with $\beta=\kappa/\gamma$. For additional details and definitions see \cite{Stewart1977}.

\subsection{Bounding the pseudoinverse of \texorpdfstring{$\Psi^*$}{the integrated library}}\label{sec:psi-pinv-norm-bound}

In this section we show that the pseudoinverse of $\Psi^*$ is bounded for all $N$ (see \Cref{cor:Psi_pseudo_bound}). This result relies on the following lemma, which shows $(\Psi^*)^T \Psi$ converges pointwise to a finite matrix.

\begin{lemma}\label{lemma:Psi_limit}
    Suppose $u_k^*(t)$ for $k=1,2,...,m$ is twice continuously differentiable for $t\in[0,t_{end}]$. Then, in the limit as $N\rightarrow\infty$, $(\Psi^*)^T \Psi^*$ converges pointwise to $B \in \mathbb{R}^{p \times p}$, defined elementwise as
    \begin{equation*}\label{eq-app:B}
        B_{\ell j} = \frac{1}{t_{end}} \int_0^{t_{end}} \xi_\ell(t) \xi_j(t) dt
        \quad\quad
        \text{where}
        \quad\quad
        \xi_j(t) =
        \begin{cases}
            1                                & \text{for } j=1  \\
            \int_0^{t} \theta^*_{j-1}(t')dt' & \text{for } j>1.
        \end{cases}
    \end{equation*}
    Here, $\theta^*_j(t) = \theta_j(u^*_1(t),u^*_2(t),...,u^*_m(t))$ represents the $j$th basis function.
\end{lemma}

\begin{proof}
    First consider $T\Theta^*$ and note that for $i=1,2,...,N$ and $j=1,2,..,p$,
    \begin{equation}
        (T\Theta^*)_{ij} = \sum_{k=1}^i T_{ik} \Theta^*_{kj} = \int_0^{t_i} \theta_j^*(t) dt + e_{i,j+1} = \xi_{j+1}(t_i) + e_{i,j+1},
    \end{equation}
    where $e_{i, j+1}$ is error due to the trapezoidal quadrature step, i.e.,
    \begin{equation}
        |e_{i,j+1}| \le \frac{(i-1)t_{end}^3}{12(N-1)^3} \max_{t\in[0,t_{end}]} \left|\frac{d^2 \theta^*_j}{dt^2}\right| = i \mathcal{O}(N^{-3}).
    \end{equation}
    We additionally define $e_{i, 1}=0$ for all $i=1,2,...,N$.

    Define $B^{(N)} := (\Psi^*)^T \Psi^*$, i.e.,
    \begin{equation}
        B^{(N)} = \frac{1}{N} \begin{bmatrix} N                     & \bm{1}^T T \Theta^*      \\
                (T \Theta^*)^T \bm{1} & (T\Theta^*)^T T \Theta^*
        \end{bmatrix}.
    \end{equation}
    Therefore, $B^{(N)}$ is given elementwise as,
    \begin{equation*}
        \begin{aligned}
            B_{\ell j}^{(N)} & = \frac{1}{N}\sum_{i=2}^N \left(\xi_\ell(t_i) + e_{i,\ell}\right)\left(\xi_j(t_i) + e_{i,j}\right)
            = \frac{1}{N}\sum_{i=2}^N \xi_\ell(t_i) \xi_j(t_i) + \mathcal{O}\left(N^{-2}\right)                                              \\
                             & = \frac{1}{t_{end}}\int_0^{t_{end}} \xi_\ell(t) \xi_j(t) dt + \delta_{\ell j}+ \mathcal{O}\left(N^{-2}\right)
            = B_{\ell j} + \mathcal{O}\left(N^{-1}\right).
        \end{aligned}
    \end{equation*}
    The summation following the second equality gives a right Riemann sum, implying, $\delta_{\ell j} \le \mathcal{O}(N^{-1})$. This result gives us the lemma statement.
\end{proof}

\begin{corollary}\label{cor:Psi_pseudo_bound}
    Suppose the assumptions of \Cref{lemma:Psi_limit} hold where $(\Psi^*)^T \Psi^*$ converges to $B$ without changing rank, i.e., there exists $N_0$ such that for $N>N_0$ the rank of $(\Psi^*)^T \Psi^*$ and $B$ are equal. Then, the value of $ \| (\Psi^*)^{\dagger} \|$ is bounded for all $N$.
\end{corollary}

\begin{proof}
    The singular values and pseudoinverse matrix norm are related as follows:
    \begin{equation}
        \|(\Psi^*)^{\dagger}\| = \frac{1}{\sigma_{min}(\Psi^*)},
    \end{equation}
    where $\sigma_{min}(\Psi^*)$ is the smallest nonzero singular value of $\Psi^*$. Since $(\Psi^*)^T \Psi^*$ converges to $B$ without changing rank, $\sigma_{min}(\Psi^*)$ must approach a constant nonzero value. Therefore, there exists a $C$ such that $\|(\Psi^*)^{\dagger}\| < C$ for all $N$.
\end{proof}

\subsection{Deriving estimators\texorpdfstring{ of $\Theta^*$ and $G$}{}}
\label{sec-app:theory-estimator}
In this section we set $\bm{u} \in \mathbb{R}^m$ and $\bm{u}^*  \in \mathbb{R}^m$ to be vectors that contain the $m$ noisy and true state variables at an arbitrary time $\tau$, e.g., $\bm{u}=[u_1(\tau),u_2(\tau),...,u_m(\tau)]$. This is in contrast to other sections of the paper where $\bm{u}$ and $\bm{u}^*$ represents the $N$ measurements and true values of an arbitrary state variable.
We additionally define $\bm{\epsilon} \in \mathbb{R}^m$ to be the noise of the $m$ state variables at the arbitrary time $\tau$, i.e., $\bm{u} = \bm{u}^* + \bm{\epsilon}$.

To present the estimators of $\Theta^*$ and $G$, we define a partial ordering of the basis multi-indices. Specifically, $\bm{\alpha}^{(k)}\le \bm{\alpha}^{(j)}$ if $\alpha_i^{(k)} \le \alpha_i^{(j)}$ for all $i \in \{1,...,m\}$.
If additionally, $\bm{\alpha}^{(k)} \ne \bm{\alpha}^{(j)}$ then a strict inequality holds, i.e., $\bm{\alpha}^{(k)} < \bm{\alpha}^{(j)}$. We will use the following shorthand notation when considering the difference between two multi-indices,
\begin{equation}
    \bm{\alpha}_{j\setminus k} := \bm{\alpha}^{(j)} - \bm{\alpha}^{(k)}.
\end{equation}

To keep the results general, we present the estimators of $\Theta^*$ and $G$ as functions that depend on the moments of the noise distribution and use the following shorthand notation,
\begin{equation}
    \mathcal{M}_{j} = \mathbb{E}[\bm{\epsilon}^{\bm{\alpha}^{(j)}}], \quad\quad
    \mathcal{M}_{j \setminus k} = \mathbb{E}[\bm{\epsilon}^{\bm{\alpha}^{(j)}-\bm{\alpha}^{(k)}}].
\end{equation}
In the specific scenario where the measurement noise is a zero-mean random Gaussian variable with variance $\sigma^2$, we have that for a general multi-index $\bm{\alpha}$,
\begin{equation}
    \mathbb{E}[\bm{\epsilon}^{\bm{\alpha}}] = \prod_{k=2,\alpha_k even}^m\sigma^{\alpha_k}(\alpha_k-1)!!,
\end{equation}
where, for odd $n$, $n!! = \prod_{k=1}^{\frac{n+1}{2}} (2k - 1)$.

For notational simplicity we also use the following definition,
\begin{equation}
    C_{jk} := \binom{\bm{\alpha}^{(j)}}{\bm{\alpha}^{(k)}}.
\end{equation}

\subsubsection{Unbiased estimator of the monomial library \texorpdfstring{$\Theta^*$}{}}\label{sec-app:centered-library}
Recall that the library $\Theta$ contains monomial basis functions evaluated at noisy state variable measurements, i.e., $\Theta_{ij} = \theta_j(u_{1}(t_i),u_{2}(t_i),...,u_{m}(t_i))$.
If $j$ is such that $|\bm{\alpha}^{(j)}| > 1$, then the corresponding monomial basis function $\theta_j(u_{1},u_{2},...,u_{m})$ is a biased estimator of $\theta_j(u^*_{1},u^*_{2},...,u^*_{m})$.
To prove that the PSDN approximation approaches the true state variable values as we increase the density of measurements, i.e., \Cref{thm:main}, it is more convenient to work with the unbiased library, which we denote as $\hat{\Theta}$.

In this section, we present a claim that explicitly gives $\hat{\Theta}$ and shows that it is an unbiased estimator. Although, we use monomial basis functions, a similar approach can be applied for any set of polynomial basis functions, e.g., Legendre polynomials, to obtain an unbiased library.

\begin{claim}\label{app-claim:centered-library}
    Let $\bm{\alpha}^{(j)}$ for $j=1,2,...,p$ represent the set of multi-indices used to define the monomial library of maximum total degree $d$. For $j=1,2,3,...,p$, iteratively define the following basis functions
    \begin{equation}\label{eq-app:unbiased-estimate}
        \hat{\theta}_j(\bm{u}) = \bm{u}^{\bm{\alpha}^{(j)}}
        - \sum_{\mathclap{\substack{
                    k \text{ s.t.} \\ \bm{\alpha}^{(k)}<\bm{\alpha}^{(j)}}}}
        C_{jk}
        \hat{\theta}_k(\bm{u}) \mathcal{M}_{j \setminus k}.
    \end{equation}
    Then $\hat{\theta}_j(\bm{u})$ is an unbiased estimate of $\theta_j(\bm{u}^*)$ for $j=1,2,...,p$.
\end{claim}

\begin{proof}
    First we find the expected value of the $j$th monomial basis function evaluated using the noisy measurements $\bm{u}$, given as $\theta_j(\bm{u}) = \bm{u}^{\bm{\alpha}^{(j)}}$.
    \begin{equation}\label{eq-app:biased-estimate}
        \begin{aligned}
            \mathbb{E}\left[\bm{u}^{\bm{\alpha}^{(j)}}\right]
             & = \mathbb{E}\left[(\bm{u}^* + \bm{\epsilon})^{\bm{\alpha}^{(j)}}\right] = \sum_{\mathclap{\substack{
            k \text{ s.t.}                                                                                          \\ \bm{\alpha}^{(k)}\le\bm{\alpha}^{(j)}}}} C_{jk} (\bm{u}^*)^{\bm{\alpha}^{(k)}} \mathcal{M}_{j \setminus k}
             & = \theta_j(\bm{u}^*) +
            \sum_{\mathclap{\substack{
            k \text{ s.t.}                                                                                          \\ \bm{\alpha}^{(k)}<\bm{\alpha}^{(j)}}}}
            C_{jk} \theta_k(\bm{u}^*) \mathcal{M}_{j \setminus k},
        \end{aligned}
    \end{equation}
    where here we use the multi-binomial theorem.

    We can then show by induction that $\hat{\theta}_j(\bm{u})$ is an unbiased estimator of $\theta_j(\bm{u}^*)$. First for $j=1$, we have that $\bm{\alpha}^{(1)}=[0,0,...,0]^T$ and $\hat{\theta}_1(\bm{u})$ is unbiased as
    \begin{equation}
        \hat{\theta}_1(\bm{u}) = 1 = \theta_1(\bm{u}^*).
    \end{equation}
    Next, pick $j > 1$ and suppose that $\hat{\theta}_k(\bm{u})$ is unbiased for all $k<j$. Note that this immediately implies $\hat{\theta}_k(\bm{u})$ is unbiased for all $k$ such that $\bm{\alpha}^{(k)}<\bm{\alpha}^{(j)}$. By combining \Cref{eq-app:unbiased-estimate,eq-app:biased-estimate}, we have that
    \begin{equation}
        \begin{aligned}
            \mathbb{E}\left[\hat{\theta}_j(\bm{u})\right]
             & = \mathbb{E}\left[\bm{u}^{\bm{\alpha}^{(j)}}\right]
            - \sum_{\mathclap{\substack{
            k \text{ s.t.}                                         \\ \bm{\alpha}^{(k)}<\bm{\alpha}^{(j)}}}}
            C_{jk}
            \mathbb{E}\left[\hat{\theta}_k(\bm{u})\right]
            \mathcal{M}_{j \setminus k}                            \\
             & = \theta_j(\bm{u}^*)
            - \sum_{\mathclap{\substack{
            k \text{ s.t.}                                         \\ \bm{\alpha}^{(k)}<\bm{\alpha}^{(j)}}}}
            C_{jk}
            \left(
            \mathbb{E}\left[\hat{\theta}_k(\bm{u})\right] - \theta_k(\bm{u}^*)
            \right)
            \mathcal{M}_{j \setminus k}
            = \theta_j(\bm{u}^*).
        \end{aligned}
    \end{equation}
    Therefore, $\hat{\theta}_j(\bm{u})$ is unbiased.
\end{proof}

When we refer to the centered monomial library $\hat{\Theta}$, we are referring to the library that contains the $p$ basis functions  $\hat{\theta}_1, \hat{\theta}_2, ..., \hat{\theta}_p$ evaluated at the noisy measurements.

\subsubsection{A consistent estimator of \texorpdfstring{$G$}{Gramian}}\label{sec-app:unbiased-G-hat}

In this section, we present a consistent estimator for $ H := G/N = \tilde{\Theta}^T \Theta^*/N$. We will assume $\tilde{\Theta}$ is independent of the noise in the measurement library $\Theta$. In practice this is not the case, as $\tilde{\Theta}$ is approximated from $\Theta$. However, this work provides a useful starting point for considering estimators for this class of Gramian matrices. Note that a similar approach can be used to find a consistent estimator of $(\Theta^*)^T \Theta^*/N$ that does not make this assumption of independence. We will use the following definitions, for $i=1,2,...,N$,
\begin{equation}
    \bm{u}^*_i := \begin{bmatrix}
        u_1^*(t_i) \\
        \vdots     \\
        u_m^*(t_i)
    \end{bmatrix}, \quad\quad
    \tilde{\bm{u}}_i := \begin{bmatrix}
        \tilde{u}_1(t_i) \\
        \vdots           \\
        \tilde{u}_m(t_i)
    \end{bmatrix}, \quad\quad
    \bm{\epsilon}_i := \begin{bmatrix}
        \epsilon_1(t_i) \\
        \vdots          \\
        \epsilon_m(t_i)
    \end{bmatrix},
\end{equation}
and $\bm{u}_i := \bm{u}^*_i + \bm{\epsilon}_i$.

The consistent estimator of $G/N$, denoted as $\hat{H}$, is given in the following claim.

\begin{claim}\label{app-claim:unbiased-Gramian}
    Let $\bm{\alpha}^{(j)}$ for $j=1,2,...,p$ represent the set of multi-indices used to define the monomial library of maximum total degree $d$ and suppose that $u_1^*(t), u_2^*(t),...,u_m^*(t)$ are bounded for all $t \in [0,t_{end}]$. Iteratively define $\hat{H}$ elementwise such that for $j,k=1,...,p$,
    \begin{equation}
        \hat{H}_{j k} =
        \frac{1}{N} \sum_{i=1}^N
        \tilde{\bm{u}}_i^{\bm{\alpha}^{(j)}}
        \bm{u}_i^{\bm{\alpha}^{(k)}}
        - \sum_{\mathclap{\substack{
                    \ell \text{ s.t.} \\ \bm{\alpha}^{(\ell)}<\bm{\alpha}^{(k)}}}}
        C_{k\ell} \hat{H}_{j \ell}
        \mathcal{M}_{k \setminus \ell}.
    \end{equation}
    Then $\hat{H}$ is a consistent estimator of $G/N$.
\end{claim}

To prove this claim we fill first present and prove a helpful lemma.

\begin{lemma}\label{lemma:expectation-limit}
    Suppose the assumptions of \Cref{app-claim:unbiased-Gramian} hold and pick $\ell,k \in \{1,...,p\}$ such that $\bm{\alpha}^{(\ell)} < \bm{\alpha}^{(k)}$. Then for all $i' \in \{1,2,...,N\}$ and $j \in \{1,2,...,p\}$,
    \begin{equation}\label{app-eq:expectation-limit-1}
        \mathbb{E}[\hat{H}_{j \ell} \bm{\epsilon}_{i'}^{\bm{\alpha}_{k\setminus \ell}}]
        =
        H_{j \ell} \mathcal{M}_{k\setminus \ell} +\mathcal{O}\left(\frac{1}{N}\right).
    \end{equation}
\end{lemma}

\begin{proof}

    We prove the statement using induction on $\ell$. First note that for $\ell=1$ for and $k \in \{2,3,...,p\}$ we have that
    \begin{equation}
        \mathbb{E}\left[\hat{H}_{j\ell} \bm{\epsilon}_{i'}^{\bm{\alpha}_{k \setminus \ell}}\right]
        =  \frac{1}{N}\sum_{i=1}^N
        \tilde{\bm{u}}_i^{\bm{\alpha}^{(j)}}
        \mathcal{M}_{k \setminus \ell}
        = H_{j \ell} \mathcal{M}_{k\setminus \ell}.
    \end{equation}
    Therefore, at $\ell=1$, \Cref{app-eq:expectation-limit-1} is satisfied.

    Then, for $\ell>1$ we assume \Cref{app-eq:expectation-limit-1} holds for $\ell' < \ell$ and note that
    \begin{equation}\label{eq:G-noise}
        \begin{aligned}
            \frac{1}{N} \sum_{i=1}^N
            \tilde{\bm{u}}_i^{\bm{\alpha}^{(j)}}
            \bm{u}_i^{\bm{\alpha}^{(\ell)}}
             & = \frac{1}{N}\sum_{i=1}^N
            \tilde{\bm{u}}_i^{\bm{\alpha}^{(j)}}
            \sum_{\mathclap{\substack{
            \ell' \text{ s.t.}              \\ \bm{\alpha}^{(\ell')}\le\bm{\alpha}^{(\ell)}}}}
            C_{\ell \ell'} (\bm{u}^*_i)^{\bm{\alpha}^{(\ell')}}
            \bm{\epsilon}_i^{\bm{\alpha}_{\ell \setminus \ell'}}
             & = \sum_{\mathclap{\substack{
            \ell' \text{ s.t.}              \\ \bm{\alpha}^{(\ell')}\le\bm{\alpha}^{(\ell)}}}}
            \frac{C_{\ell \ell'}}{N}\sum_{i=1}^N
            \tilde{\bm{u}}_i^{\bm{\alpha}^{(j)}}
            (\bm{u}^*_i)^{\bm{\alpha}^{(\ell')}}
            \bm{\epsilon}_i^{\bm{\alpha}_{\ell \setminus \ell'}}.
        \end{aligned}
    \end{equation}
    Using this, we then have that
    \begin{align*}
        \mathbb{E}\left[\hat{H}_{j \ell} \bm{\epsilon}_{i'}^{\bm{\alpha}_{k \setminus \ell}}\right]
         & =\begin{multlined}[t]
                \sum_{\mathclap{\substack{
                            \ell' \text{ s.t.} \\ \bm{\alpha}^{(\ell')}\le\bm{\alpha}^{(\ell)}}}}
                \frac{C_{\ell\ell'}}{N}\sum_{i=1}^N
                \tilde{\bm{u}}_i^{\bm{\alpha}^{(j)}}
                (\bm{u}^*_i)^{\bm{\alpha}^{(\ell')}}
                \mathbb{E}\left[
                \bm{\epsilon}_i^{\bm{\alpha}_{\ell \setminus \ell'}}
                \bm{\epsilon}_{i'}^{\bm{\alpha}_{k \setminus \ell}}
                \right]
                - \sum_{\mathclap{\substack{
                            \ell' \text{ s.t.} \\ \bm{\alpha}^{(\ell')}<\bm{\alpha}^{(\ell)}}}}
                C_{\ell\ell'}
                \mathbb{E}\left[\hat{H}_{j \ell'} \bm{\epsilon}_{i'}^{\bm{\alpha}_{k \setminus \ell}}\right]
                \mathcal{M}_{\ell \setminus \ell'}
            \end{multlined}                                                                                                  \\
         & =\begin{multlined}[t]
                \sum_{\mathclap{\substack{
                            \ell' \text{ s.t.} \\ \bm{\alpha}^{(\ell')}\le\bm{\alpha}^{(\ell)}}}}
                C_{\ell\ell'}
                H_{j \ell'}
                \mathcal{M}_{\ell \setminus \ell'}\mathcal{M}_{k \setminus \ell}
                - \sum_{\mathclap{\substack{
                            \ell' \text{ s.t.} \\ \bm{\alpha}^{(\ell')}<\bm{\alpha}^{(\ell)}}}}
                C_{\ell\ell'}
                \mathbb{E}\left[\hat{H}_{j \ell'} \bm{\epsilon}_{i'}^{\bm{\alpha}_{k \setminus \ell}}\right]
                \mathcal{M}_{\ell \setminus \ell'}\\
                + \frac{1}{N}\sum_{\mathclap{\substack{
                            \ell' \text{ s.t.} \\ \bm{\alpha}^{(\ell')}\le\bm{\alpha}^{(\ell)}}}}
                C_{\ell\ell'}
                \tilde{\bm{u}}_{i'}^{\bm{\alpha}^{(j)}}
                (\bm{u}^*_{i'})^{\bm{\alpha}^{(\ell')}}
                \left(\mathcal{M}_{k\setminus\ell'} - \mathcal{M}_{\ell \setminus \ell'}\mathcal{M}_{k \setminus \ell}\right)
            \end{multlined}                                                                                 \\
         & = \begin{multlined}[t]
                 H_{j \ell} \mathcal{M}_{k \setminus \ell}
                 +
                 \sum_{\mathclap{\substack{
                             \ell' \text{ s.t.} \\ \bm{\alpha}^{(\ell')}<\bm{\alpha}^{(\ell)}}}}
                 C_{\ell\ell'}
                 \left(
                 H_{j \ell'} \mathcal{M}_{k \setminus \ell}
                 - \mathbb{E}\left[
                     \hat{H}_{j \ell'} \bm{\epsilon}_{i'}^{\bm{\alpha}_{k \setminus \ell}}\right]
                 \right)
                 \mathcal{M}_{\ell \setminus \ell'} \\
                 +
                 \frac{1}{N}\sum_{\mathclap{\substack{
                             \ell' \text{ s.t.} \\ \bm{\alpha}^{(\ell')}\le\bm{\alpha}^{(\ell)}}}}
                 C_{\ell\ell'}
                 \tilde{\bm{u}}_{i'}^{\bm{\alpha}^{(j)}}
                 (\bm{u}^*_{i'})^{\bm{\alpha}^{(\ell')}}
                 \left(
                 \mathcal{M}_{k\setminus\ell'}-
                 \mathcal{M}_{\ell \setminus \ell'}
                 \mathcal{M}_{k \setminus \ell}
                 \right).
             \end{multlined}
    \end{align*}
    Here, for the second equality we are using the fact that the measurement noise is independent across the $N$ samples. For the third equality, we are pulling out the term from the first summation that corresponds to $G_{j\ell}\mathcal{M}_{k \setminus \ell}$.

    As desired, since \Cref{app-eq:expectation-limit-1} holds for $\ell' < \ell$,
    \begin{equation}
        \mathbb{E}\left[\hat{H}_{j \ell} \bm{\epsilon}_{i'}^{\bm{\alpha}_{k \setminus \ell}}\right]
        = H_{j \ell} \mathcal{M}_{k \setminus \ell} + \mathcal{O}\left(\frac{1}{N}\right).
    \end{equation}
\end{proof}

We use this result to prove \Cref{app-claim:unbiased-Gramian}.

\begin{proof}[Proof of \Cref{app-claim:unbiased-Gramian}]
    To prove the result we need to show that $\hat{H}$ is unbiased and that the variance goes to zero as we increase the sample size.

    First we show that $\hat{H}$ is unbiased. Taking the expectation of \Cref{eq:G-noise}, we immediately have that
    \begin{equation}
        \mathbb{E}\left[
        \frac{1}{N} \sum_{i=1}^N
        \tilde{\bm{u}}_i^{\bm{\alpha}^{(j)}}
        \bm{u}_i^{\bm{\alpha}^{(k)}}
        \right]
        = H_{jk} + \sum_{\mathclap{\substack{
                    \ell \text{ s.t.} \\ \bm{\alpha}^{(\ell)}<\bm{\alpha}^{(k)}}}}
        C_{k \ell}
        H_{j \ell}
        \mathcal{M}_{k \setminus \ell}.
    \end{equation}
    Using the same inductive reasoning as that given in the proof to \Cref{app-claim:centered-library}, we then have that $\hat{H}$ is an unbiased estimator of $H$.

    Next, we show that the variance of each element of $\hat{H}$ goes to zero in the limit as $N \rightarrow \infty$. First we rewrite $\hat{H}_{jk}$ using the definition of $\hat{H}$ and calculation shown in \Cref{eq:G-noise},
    \begin{equation}\label{app-eq:Ghat-2}
        \hat{H}_{jk}
        = H_{jk} +
        \sum_{\mathclap{\substack{
                    \ell \text{ s.t.} \\ \bm{\alpha}^{(\ell)}<\bm{\alpha}^{(k)}}}}
        C_{k \ell}
        \left(
        \frac{1}{N}\sum_{i=1}^N
        \tilde{\bm{u}}_i^{\bm{\alpha}^{(j)}}
        (\bm{u}^*_i)^{\bm{\alpha}^{(\ell)}}
        \bm{\epsilon}_i^{\bm{\alpha}_{k \setminus \ell}}
        - \hat{H}_{j \ell}\mathcal{M}_{k \setminus \ell}
        \right).
    \end{equation}
    Therefore, we can show $\mathbb{V}[\hat{H}_{jk}]$ goes to zero by showing that the variance of the second term on the right hand side of \Cref{app-eq:Ghat-2} goes to zero. We again use a proof by induction.

    For $k=1$ we have that $\hat{H}_{jk} = H_{jk}$ and, therefore $\mathbb{V}[\hat{H}] = 0$. For $k > 0$, we consider the following variance,
    \begin{align*}
        \mathbb{V} \left[
        \frac{1}{N}\sum_{i=1}^N
        \tilde{\bm{u}}_i^{\bm{\alpha}^{(j)}}
        (\bm{u}^*_i)^{\bm{\alpha}^{(\ell)}}
        \bm{\epsilon}_i^{\bm{\alpha}_{k \setminus \ell}}
        - \hat{H}_{j \ell}\mathcal{M}_{k \setminus \ell}
        \right]
         & = \begin{multlined}[t]\frac{1}{N^2}\sum_{i,i'=1}^N
                 (\tilde{\bm{u}}_i\tilde{\bm{u}}_{i'})^{\bm{\alpha}^{(j)}}
                 (\bm{u}^*_i\bm{u}^*_{i'})^{\bm{\alpha}^{(\ell)}}
                 \mathcal{M}_{k \setminus \ell}^2
                 + \mathbb{E}\left[\hat{H}_{j \ell}^2\right] \mathcal{M}_{k \setminus \ell}^2 \\
                 - \frac{2}{N} \sum_{i=1}^N
                 \tilde{\bm{u}}_i^{\bm{\alpha}^{(j)}}
                 (\bm{u}^*_i)^{\bm{\alpha}^{(\ell)}}
                 \mathbb{E}\left[
                 \bm{\epsilon}_i^{\bm{\alpha}_{k \setminus \ell}} \hat{H}_{j \ell}
                 \right]
                 \mathcal{M}_{k \setminus \ell} \\
                 + \frac{1}{N^2} \sum_{i=1}^N
                 (\tilde{\bm{u}}_i)^{2\bm{\alpha}^{(j)}}
                 (\bm{u}^*_i)^{2\bm{\alpha}^{(\ell)}}
                 \left(
                 \mathbb{E}[\bm{\epsilon}^{2\bm{\alpha}_{k \setminus \ell}}]
                 - \mathcal{M}_{k \setminus \ell} ^2
                 \right)
             \end{multlined} \\
         & = H_{j \ell}^2 \mathcal{M}_{k \setminus \ell}^2 + \mathbb{E}\left[\hat{H}_{j \ell}^2\right]
        \mathcal{M}_{k \setminus \ell}^2
        - 2 H_{j \ell}^2 \mathcal{M}_{k \setminus \ell}^2
        + \mathcal{O}\left(\frac{1}{N}\right).
    \end{align*}
    In this derivation we use the result from \Cref{lemma:expectation-limit} to replace $\mathbb{E}\left[\bm{\epsilon}_i^{\bm{\alpha}_{k \setminus \ell}} \hat{H}_{j \ell}\right]$ with $H_{j \ell} \mathcal{M}_{k \setminus \ell}$.
    By induction, we have that for $\ell < k$, $\mathbb{V}[\hat{H}_{j\ell}]$ goes to zero and thus,
    \begin{equation}
        \lim_{N\rightarrow \infty} \mathbb{E}[\hat{H}_{j\ell}^2] = H_{j\ell}^2.
    \end{equation}
    Therefore, in the limit as $N\rightarrow \infty$, $\mathbb{V}[\hat{H}_{jk}]$ goes to zero.
\end{proof}

\subsection{Value of smoothing hyperparameter leads to feasible SOCP}
\label{sec-app:smoothing-hyper-param}
Using the value of $C$ in \Cref{eq:C} guarantees feasibility of the SOCP problem, i.e., \Cref{eq:SOCP}, at small values of $\gamma$. To see this, we will show that the constraints in \Cref{eq:SOCP} are satisfied when $\dot{\bm{u}} = \dot{\tilde{\bm{u}}}$ as defined by \Cref{eq:tildeudot}. We immediately have that $\|D\dot{\tilde{\bm{u}}}\| = C$, implying the first constraint in \Cref{eq:SOCP} is satisfied. For the second constraint, we have that,
\begin{equation*}
    \begin{aligned}
        \left\|
        \begin{bmatrix} \bm{1} & T \end{bmatrix} \begin{bmatrix} u_0 \\ \dot{\tilde{\bm{u}}} \end{bmatrix}- P_{\tilde{\Phi}}\tilde{\bm{u}}
        \right\|
         & = \|u_0 \bm{1} + T \tilde{\Theta}  \left(
        (\tilde{\Phi}^T \tilde{\Phi})^{-1} \tilde{\Phi}^T
        \right)_{2:,:} \tilde{\bm{u}}
        - P_{\tilde{\Phi}} \tilde{\bm{u}}
        \|                                                                            \\
         & =
        \|u_0 \bm{1} + \tilde{\Phi} (\tilde{\Phi}^T \tilde{\Phi})^{-1} \tilde{\Phi}^T \tilde{\bm{u}}
        - \left(\left(
        (\tilde{\Phi}^T \tilde{\Phi})^{-1} \tilde{\Phi}^T
        \right)_{1,:} \bm{\tilde{u}}\right)\bm{1}  - P_{\tilde{\Phi}}\tilde{\bm{u}}\| \\
         & = 0,
    \end{aligned}
\end{equation*}
where the final equality holds by setting $u_0 = \left(\left(
    (\tilde{\Phi}^T \tilde{\Phi})^{-1} \tilde{\Phi}^T
    \right)_{1,:} \bm{\tilde{u}}\right)$. Thus, the problem is feasible for all $\gamma \ge 0$. Recall that the `$2:,:$' subscript implies we are removing the top row of the matrix. Similarly, the `$1,:$' subscript implies we are considering only the top row of the matrix.

\setcounter{theorem}{0}
\setcounter{lemma}{0}
\setcounter{corollary}{0}
\setcounter{claim}{0}
\setcounter{remark}{0}
\section{Additional results}\label{sec-app:results}

\subsection{Additional results for the Lorenz 96 model}
\label{sec-app:Lorenz96}

\begin{figure}
    \centering
    \includegraphics[trim=0 0 0 0,clip]{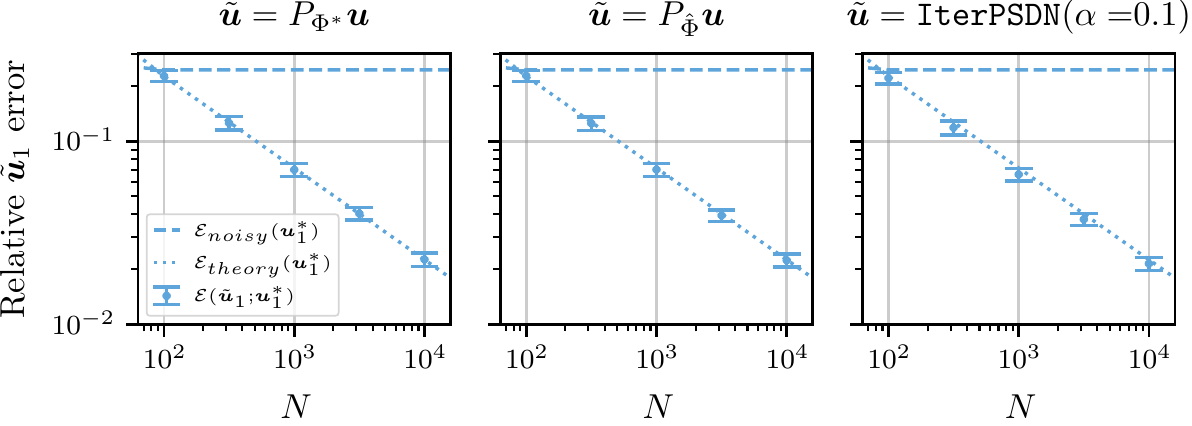}
    \caption{Comparing denoising error with theoretical predictions for the Lorenz 96 model. Results show mean $\pm$ std for 50 sample replications at noise level $\sigma^2=1$ when the integrated library $\Phi^*$ is known (left), when using PSDN (middle), and when using IterPSDN, \Cref{alg:PSDN} (right). Although we only show the results for the first state of the Lorenz 96 model, nearly identical results were observed for the other 5 states. For system parameters see \Cref{tab}.}
    \label{fig:PSDN-Theory4}
\end{figure}

\begin{figure}
    \centering
    \includegraphics[trim=0 25 0 0, clip]{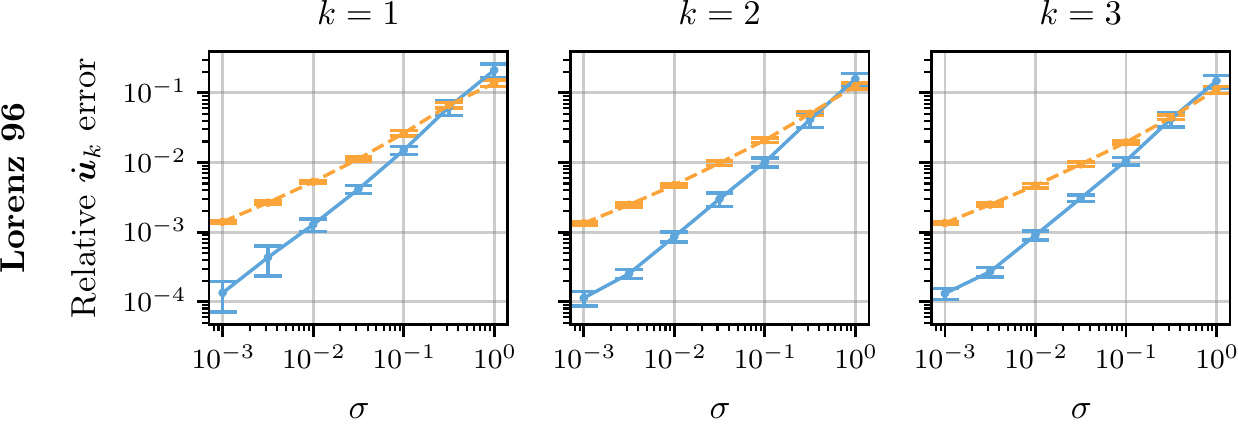}

    \vspace{.2cm}

    \includegraphics[trim=0 0 0 0, clip]{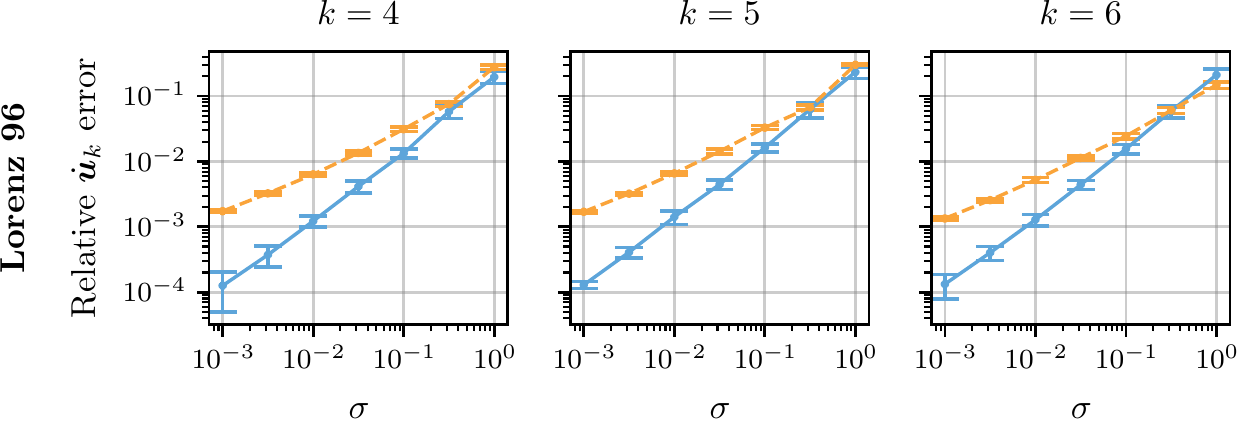}

    \includegraphics{figures/legends/legend_socp_sm_tikreg_h}

    \caption{Derivative estimation errors for the 6 states of the Lorenz 96 system. Showing the mean $\pm$ sem of 30 realizations with $N=2000$ at each noise level. For system parameters see \Cref{tab}.}
    \label{fig:DSINDy-Lorenz2}
\end{figure}

Here, we provide additional results for the Lorenz 96 model. \Cref{fig:PSDN-Theory4} shows results of the PSDN and IterPSDN compared to theoretical estimations (see \Cref{sec:results-PSDN-theory} for additional details). \Cref{fig:DSINDy-Lorenz2} shows the error of the derivation estimation from the IRW-SOCP step of DSINDy.

\subsection{Additional results for Duffing oscillator (PS1)}
\label{sec-app:Duffing}

\begin{figure}
    \centering
    \includegraphics[trim=0 25 0 0, clip, align=c]{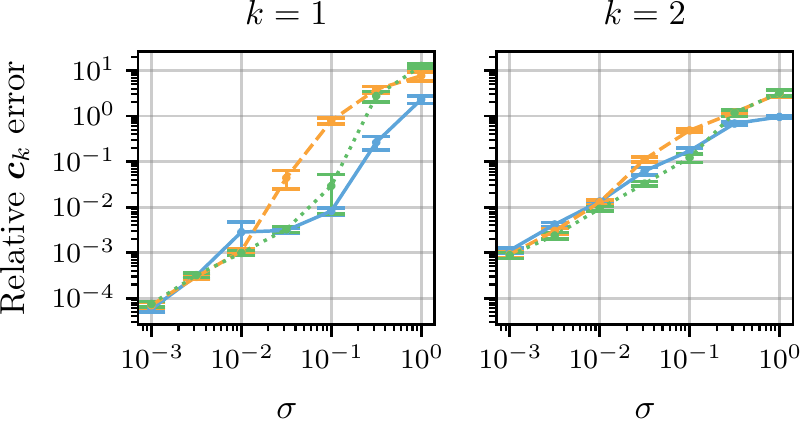}
    \hspace{.7cm}
    \includegraphics[align=c]{figures/legends/legend_socp_sm_lasso_WSINDY}

    \includegraphics[trim=0 0 0 10, clip]{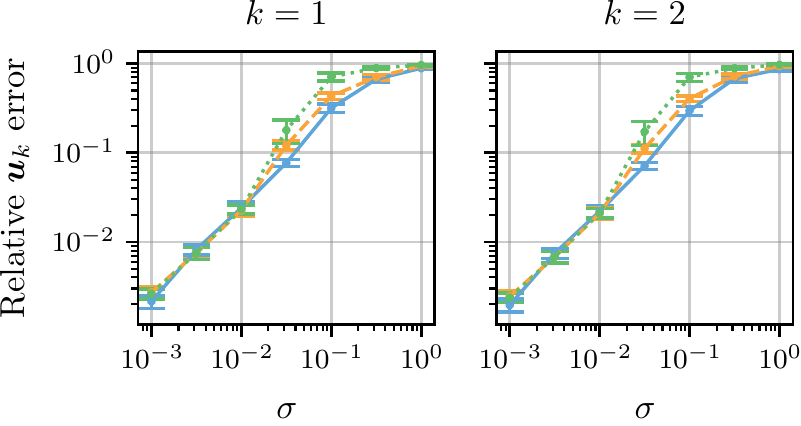}
    \hspace{.2cm}
    \includegraphics{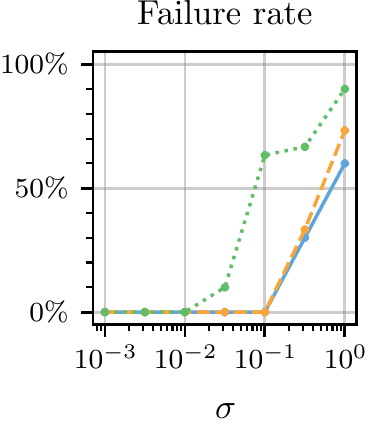}

    \caption{Coefficient and reconstruction error for Duffing oscillator system (PS1). Showing the mean $\pm$ sem of 30 realizations with $N=1000$ at each noise level. The failure rate represents the fraction of simulations out of 30 that failed. For system parameters see \Cref{tab}.}
    \label{fig:DSINDy1}
\end{figure}

Here, we provide results on coefficient error and state reconstruction error for the Duffing oscillator system when using Parameter Set 1 (see \Cref{fig:DSINDy1}).

\subsection{Examining hyperparameter \texorpdfstring{$\gamma$}{} selection method in DSINDy}\label{sec-app:gamma-method}

\begin{figure}
    \includegraphics[trim=0 24 0 0,clip]{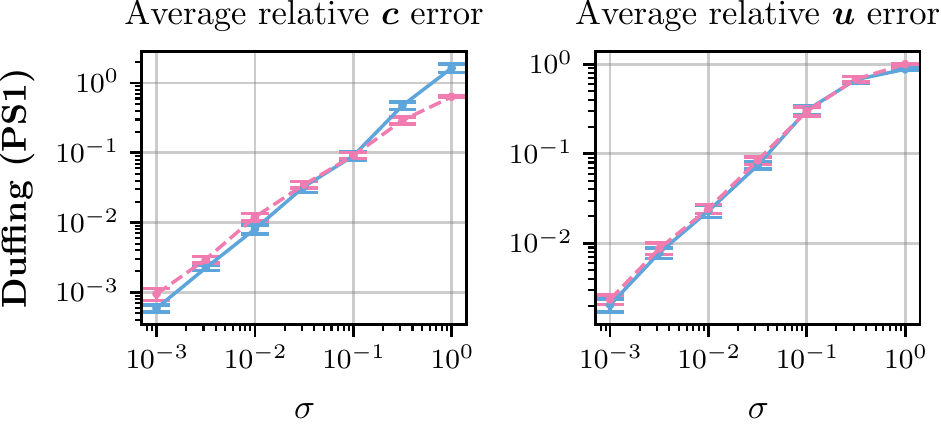}
    \includegraphics[trim=0 24 0 0,clip]{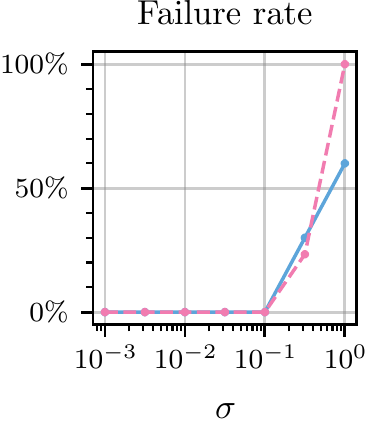}

    \includegraphics[trim=0 24 0 13,clip]{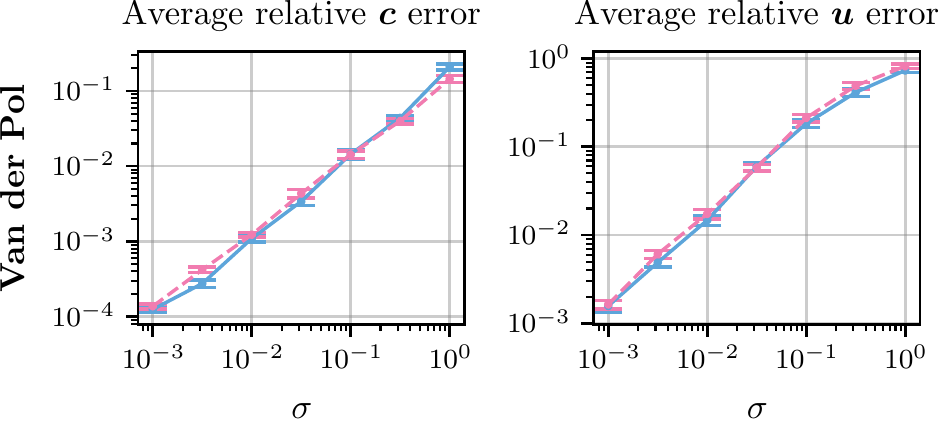}
    \includegraphics[trim=0 24 0 13,clip]{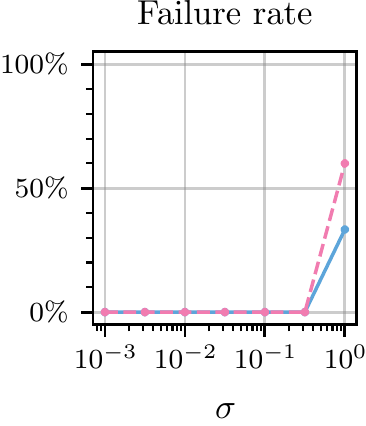}

    \includegraphics[trim=0 0 0 13,clip]{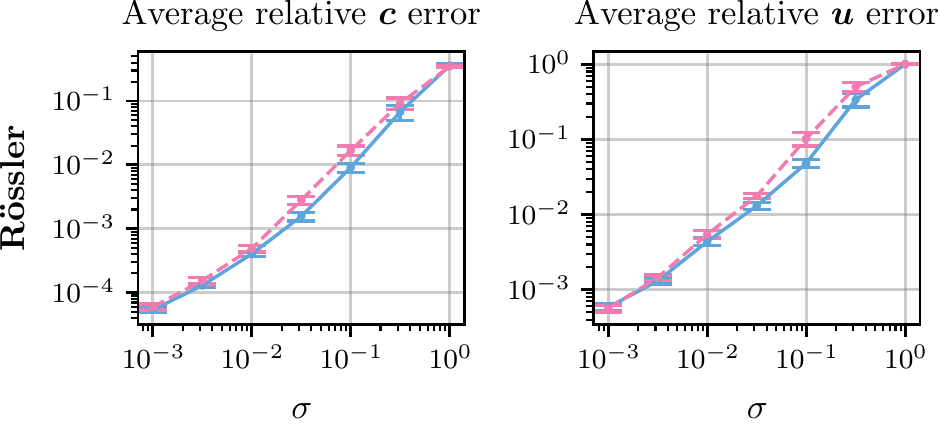}
    \includegraphics[trim=0 0 0 13,clip]{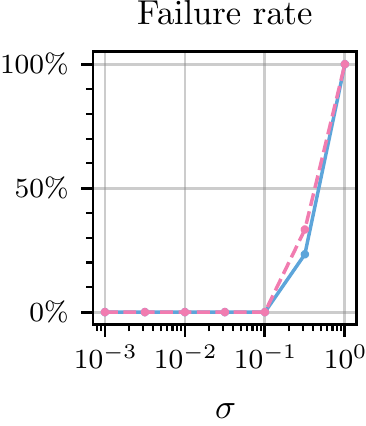}

    \vspace{.5cm}

    \includegraphics[trim=0 0 135 0,clip,align=c]{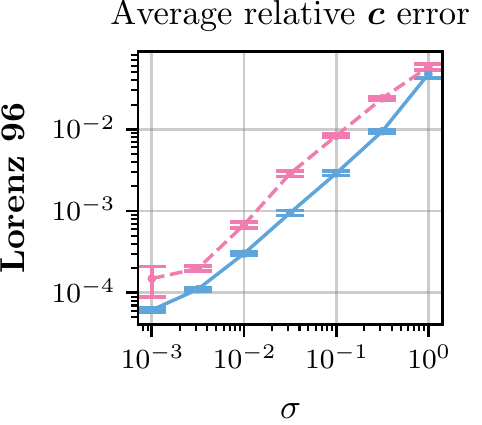}
    \includegraphics[trim=15 0 0 0,clip,align=c]{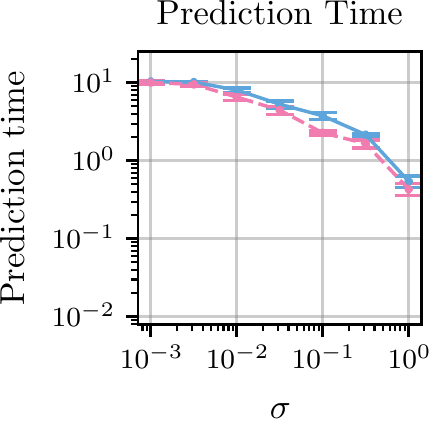}
    \hspace{.4cm}
    \includegraphics[align=c]{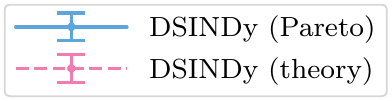}

    \caption{Coefficient and reconstruction error for example systems. Showing the mean $\pm$ sem of 30 realizations with $N=1000$ at each noise level ($N=2000$ for Lorenz 96 model). The failure rate represents the fraction of simulations out of 30 that failed. For system parameters see \Cref{tab}.}
    \label{fig-app:two-gamma-methods-1}
\end{figure}

\begin{figure}[t]
    \centering
    \includegraphics[trim=0 32 0 0, clip]{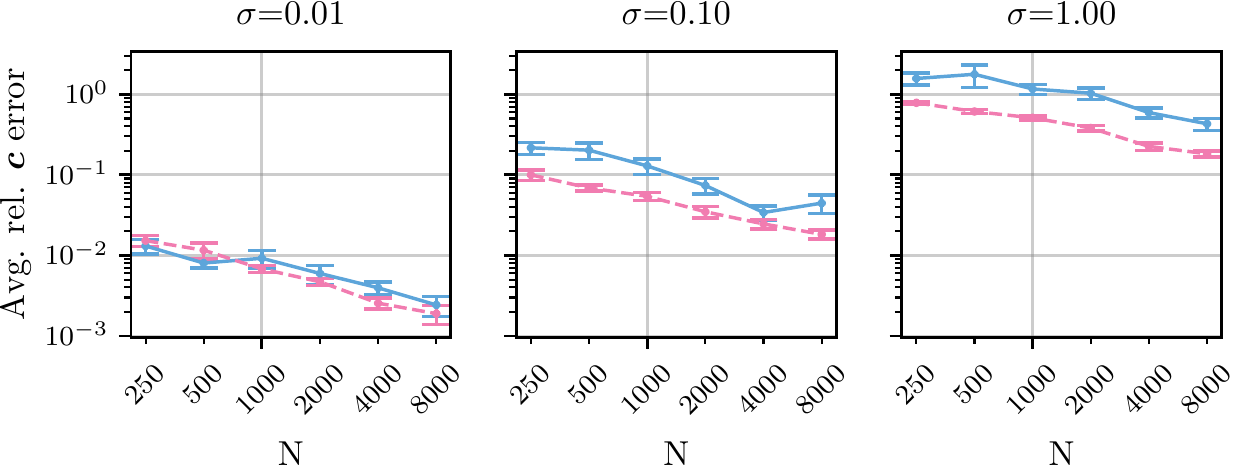}
    \includegraphics[trim=0 32 0 10, clip]{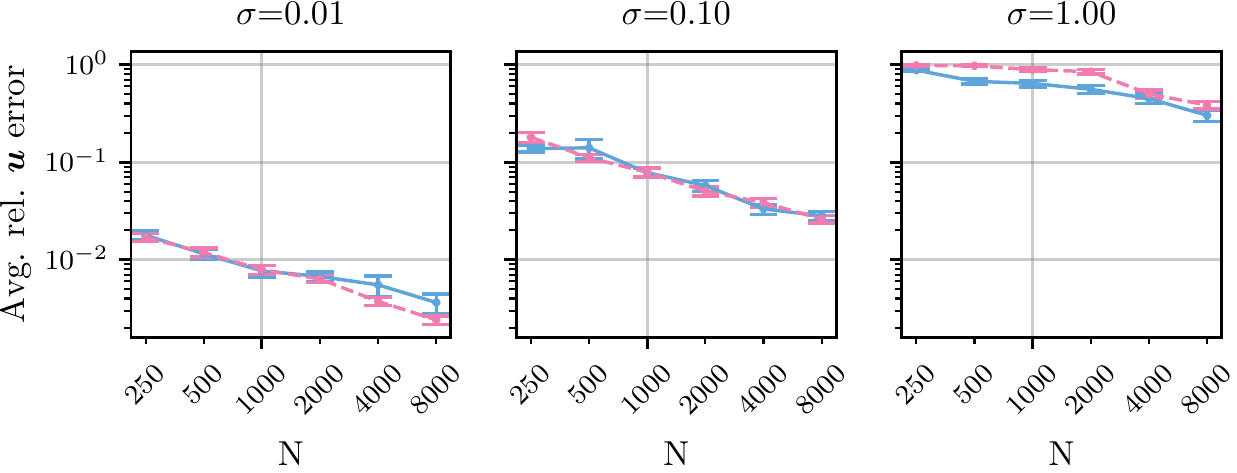}
    \includegraphics[trim=0 0 0 10, clip]{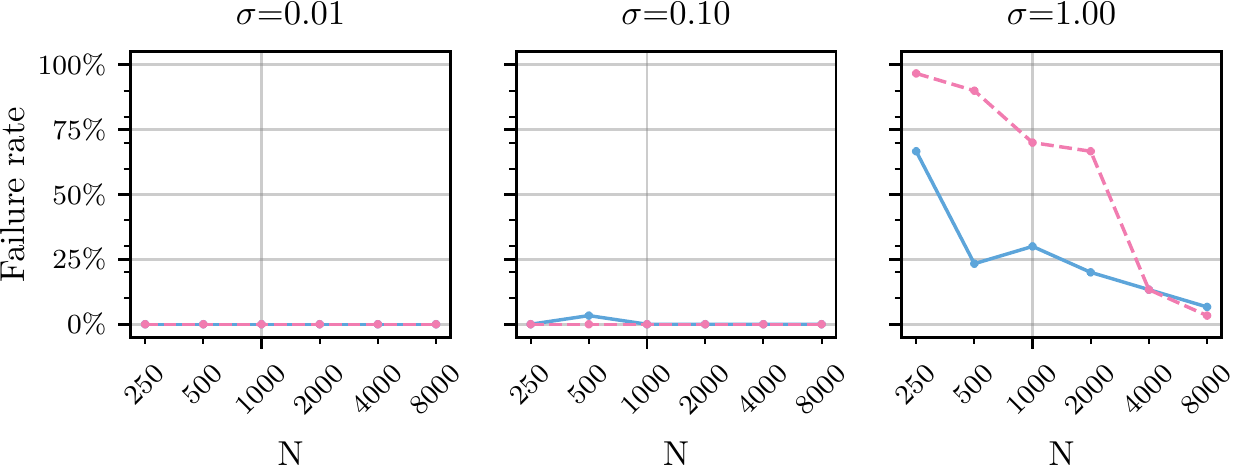}
    \includegraphics{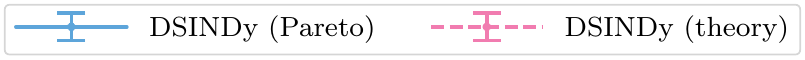}
    \caption{Coefficient and reconstruction error for Duffing oscillator system (PS2). Each column represents a different noise level. Showing the mean $\pm$ sem of 30 realizations at each noise level. The failure rate represents the fraction of simulations out of 30 that failed. For system parameters see \Cref{tab}.}
    \label{fig-app:two-gamma-methods-2}
\end{figure}

Here, we present results showing that the two approaches for selecting the data-matching hyperparameter $\gamma$ in \Cref{eq:SOCP} lead to roughly equivalent results (\Cref{fig-app:two-gamma-methods-1,fig-app:two-gamma-methods-2}).
Note that for conciseness in these figures we show the relative error of the coefficients and system state reconstructions averaged across the states of a given system. For details on these systems sec \Cref{sec:ode-systems}.

\bibliographystyle{elsarticle-num}
\bibliography{EqnDiscov-Paper}

\end{document}